\documentclass[11pt,fleqn,english]{amsart}%article}

\RequirePackage[OT1]{fontenc}
 \RequirePackage{amsthm,amsmath}
\RequirePackage[numbers]{natbib}
\usepackage{amsmath}   % many want amsmath extensions
\usepackage{amsfonts,amssymb, color}
\def\Xint#1{\mathchoice
   {\XXint\displaystyle\textstyle{#1}}%
   {\XXint\textstyle\scriptstyle{#1}}%
   {\XXint\scriptstyle\scriptscriptstyle{#1}}%
   {\XXint\scriptscriptstyle\scriptscriptstyle{#1}}%
   \!\int}
\def\XXint#1#2#3{{\setbox0=\hbox{$#1{#2#3}{\int}$}
     \vcenter{\hbox{$#2#3$}}\kern-.5\wd0}}

\def\dashint{\Xint-}

\usepackage{graphicx} % avoid epsfig or earlier such packages
\usepackage{url}      % for URLs

\usepackage{amsfonts}
\usepackage[english]{babel}
\usepackage{latexsym}
\usepackage{enumerate}
\usepackage{verbatim}
\newcommand\dela[1]{}
\newcommand\delzb[1]{}%This is new command. I use it in May 2013. Zdzislaw
\newcommand\delb[1]{}%03 September 2013
\newcommand\deln[1]{}%%28 November 2013
\newcommand{\bcase}{\begin{cases}}
\newcommand{\ecase}{\end{cases}}

\newcommand\cadlag{c{\`a}dl{\`a}g\,\,}
\newcommand\inte{\mbox{\sl int}}
\newcommand{\BE}{E}
\newcommand{\En}{\mathbb{E}_n}
\newcommand{\eps}{\varepsilon}
\newcommand{\Lve}{\lVert}
\newcommand{\Rve}{\rVert}
\newcommand{\lve}{\lvert}
\newcommand{\rve}{\rvert}
\newcommand\del[1]{}

\del{

\newcommand\del[1]{}
}

\newcounter{zaehler}

\newcounter{gg1}

\newcounter{gr1}
\newenvironment{numlist}
{\begin{list} { (\roman{gr1}):} {\usecounter{gr1}
\setlength{\leftmargin}{0.9cm}
\setlength{\topsep}{0.1cm} \setlength{\itemsep}{0.0cm}
\setlength{\parsep}{0.1cm} \setlength{\itemindent}{0.9cm}
\setlength{\parskip}{0.0cm}}} {\end{list}}
\newcounter{gr1n}
\newenvironment{numlistn}
{\begin{list} { (\roman{gr1n})} {\usecounter{gr1n}
\setlength{\leftmargin}{0.9cm}
\setlength{\topsep}{0.1cm} \setlength{\itemsep}{0.0cm}
\setlength{\parsep}{0.1cm} \setlength{\itemindent}{-0.7cm}
\setlength{\parskip}{0.0cm}}} {\end{list}}
\newcounter{gr11}

\newcounter{lil}
\newenvironment{steps}
{\begin{list} { \bf Step (\Roman{lil})} { \em %\usecounter{lil}
\setlength{\leftmargin}{0.0cm}
\setlength{\topsep}{0.2cm} \setlength{\itemsep}{0.0cm}
\setlength{\parsep}{0.1cm} \setlength{\itemindent}{0.8cm}
\setlength{\parskip}{0.0cm}}} {\end{list}}

\newcounter{lil22}
\newenvironment{steps-2}
{\em \begin{list} {  (\roman{lil22})} {\usecounter{lil22} \em
\setlength{\leftmargin}{1.3cm}
\setlength{\topsep}{0.2cm} \setlength{\itemsep}{0.0cm}
\setlength{\parsep}{0.1cm} \setlength{\itemindent}{0.4cm}
\setlength{\parskip}{0.0cm}}} {\end{list}}
\newcommand{\embed}{\hookrightarrow}

\def\vt{\vartheta}
\def\vp{\varphi}
\def\Re{\mathop{ \rm Re }}
\def\Arg{\mathop{ \rm Arg }}
\def\sgn{\mathop{ \rm sgn }}

\newcommand\lb{\langle}
\newcommand\rb{\rangle}

\newcommand{\la}{{\langle}}
\newcommand{\ra}{{\rangle}}
\newcommand{\BB} {E}

\newcommand{\Law}{{{ \mathcal L }aw}}

\newcommand{\Rb}[1]{{\mathbb{R}_{#1}}}

\DeclareMathOperator{\leb}{Leb}
\DeclareMathOperator{\Leb}{Leb}

\newcommand{\E}{{E}}

\input colordvi

\numberwithin{equation}{section}

\newtheorem{theorem}{Theorem}[section]
\newtheorem{notation}{Notation}[section]
\newtheorem{claim}{Claim}[section]
\newtheorem{lemma}[theorem]{Lemma}%[section]
\newtheorem{corollary}[theorem]{Corollary}%[section]
\newtheorem{assumption}[theorem]{Assumption}

\newtheorem{definition}[theorem]{Definition}%[section]
\newtheorem{remark}[theorem]{Remark}%[section]

\newtheorem{proposition}[theorem]{Proposition}%[section]

\begin{document}

\title[Stochastic Reaction Diffusion Equation of Jump Processes \today]{
Stochastic Reaction Diffusion Equation driven by jump
processes}

\author{Zdzis{\l}aw Brze{\'z}niak}
%\address{}
%\curraddr{}
%\email{}
%\thanks{}\author[Z. Brze{\'z}niak et al.]{,
% and  }
\address{Department of Mathematics \\
 University of York, Heslington, York YO10
5DD, UK} \email{zb500@york.ac.uk}

\author{Erika Hausenblas}

\address{Department of Mathematics and Informationtechnology, Montanuniversity of Leoben,
Franz Josef Strasse 18, 8700 Leoben, Austria} \email{
erika.hausenblas@sbg.ac.at}
\thanks{This work was supported by the FWF-Project
P17273-N12}

\author{Paul Andr\'e  {Razafimandimby}}
\address{Department of Mathematics and Informationtechnology, Montanuniversity of Leoben,
Franz Josef Strasse 18, 8700 Leoben, Austria} \email{
paul.razafimandimby@unileoben.ac.at}%
  %\thankstext{t2}{Footnote to the first author with the `thankstext' command.}

\begin{abstract}:
The  paper studies a  reaction diffusion equation driven by
Poissonian, respectively, L\'evy noise. It first shows  existence
of
a martingale solution for  parabolic SPDEs %of a parabolic type
driven by a Poisson random measure. This result is then
transferred to  parabolic SPDE driven by L\'evy noise. This
result allows us to establish the existence of a martingale
solution of reaction diffusion type equation driven by Poissonian
noise, respectively, L\'evy noise. The result answer positively an
open question about existence of martingale
solutions of nonlinear SPDEs/reaction-diffusion equation driven by genuine L\'evy processes on Banach spaces.
\end{abstract}

%\amssubj{}
\maketitle

\section{Introduction}\label{sec_intro}
The main subject of this paper  is a study of  stochastic reaction
diffusion equation driven by a L\'evy  (or  Poissonian) noise. For
instance, let $L=\{L(t) :t\ge 0  \}$ be a one dimensional L\'evy
process whose characteristic  measure $\nu$ has finite $p$-moment
for some  $p\in(1,2]$,   let $\mathcal{O}\subset {\mathbb{R}}^ d$ be a bounded domain
with smooth boundary and  let $ \Delta$ be the Laplace operator  with
the Dirichlet boundary conditions.
 A typical example of interest
  is the following   stochastic partial differential equation
\begin{eqnarray}\label{spde01}
\\
\nonumber \left\{
\begin{array}{rcl} du(t,\xi) &=& \Delta u(t,\xi) \: dt + \bigl[u(t,\xi)-u(t,\xi)^3\bigr]\, dt  \phantom{\Big|}
\\ & +&
\frac{ \sqrt{|u(t,\xi )|}}{ 1+\sqrt{|u(t,\xi )|} }
\;d{L}(t),\quad t>0,\phantom{\Big|}
\\
u(t,\xi) &=& 0,\quad \xi\in \partial \mathcal{O}, \phantom{\Big|}
\\
u(0,\xi) &=& u_0(\xi),\quad \xi\in \mathcal{O} \phantom{\Big|}.
\end{array}\right.
\end{eqnarray}

One consequence of the main result, i.e.\ Theorem
\ref{Th:general}, is that for any  $u_0\in
C_0(\overline{\mathcal{O}})$ (see )  there exists an
$C_0(\overline{\mathcal{O}})$-valued \del{c\`{a}dl\`{a}g} process
$u=\{ u(t): t\geq 0  \}$ which is a martingale solution  to
problem \eqref{spde01}. The  main result allows  a treatment of
equations with  more general coefficients than those in problem
\eqref{spde01}. For instance, the diffusion coefficient
$g(u)=\sqrt{|u|}/(1+\sqrt{|u|})$  is
just an example of a bounded and continuous function,  the drift term
$f(u)=u-u^3$ is  an example of a dissipative function
$f:\mathbb{R}\to \mathbb{R}$ of polynomial growth, and the Laplace
operator $\Delta$ is a special case of a second order, with
variable coefficients,  uniformly elliptic dissipative operator.
Moreover, the results are applicable to equations with
infinite-dimensional L\'evy processes as well as  more general
classes of  initial data, for instance the $L^q(\mathcal{O})$
spaces with $q\geq p$, and % as well as
 the Sobolev spaces $W^
{\gamma,q}_0(\mathcal{O})$. The above and other examples are
presented in sections \ref{reaction1}, \ref{reaction1-s-t} and \ref{det-set-B}.

In our  paper we adopt an approach used in a recent work by one of the
authors \cite{Brz+Gat_1999}, in which a similar problem with
Wiener processes was treated,  but with major improvements for the
following reasons.

Firstly, the It\^o integral in martingale type $2$ Banach spaces
with respect to cylindrical Wiener processes on which
\cite{Brz+Gat_1999} is so heavily relied, is replaced by the It\^o
integral in martingale type $p$ Banach spaces with respect to
Poisson random measures, see \cite{maxreg}. Secondly, the
compactness argument in \cite{Brz+Gat_1999} depends  on the
H\"older continuity of the trajectories of the stochastic
convolutions driven by a Wiener process. However, since  the
trajectories of the stochastic convolutions driven by a L\'evy
process are not continuous,  the counterpart of the H\"older
continuity, i.e., the \cadlag \ property of the trajectories seems
natural  to be used. Unfortunately, as has been shown by many
counterexamples, see for instance a recent monograph
\cite{Peszat_Z_2007}, as well as even the recent papers
\cite{Brz+Z_2009} and \cite{Brz+GIPPZ_2009},  the trajectories of
the stochastic convolutions driven by a L\'evy process may not
even be \cadlag \ in the space where the L\'evy process lives in,
and, hence, this issue has to be handled with special care. A
third but not least  difference is  that we do not use  a martingale
representation theorem. In \cite{Brz+Gat_1999} a known result by
Dettweiler \cite{Det-88} was used. In the case of a L\'evy
process, to our great surprise, an appropriate martingale
representation Theorem has not been found. In the current paper we
proved  a generalization of the result from Dettweiler
\cite{Det-88}. Finally, instead of using stopping times as in
\cite{Brz+Gat_1999}, we apply  interpolation methods in order to
control certain norms of the solution. On the other hand, our paper confirms an observation already made in
an earlier paper \cite{maxreg} by the first two named authours that the theory of stochastic integration with
respect to a Poisson random measure (with the intensity measure $\nu$) in martingale type $p$
Banach spaces is similar to the theory of stochastic integration with respect to a cylindrical Wiener process
in martingale type $2$ Banach spaces provided the $\gamma$-radonifying norm is replaced by the $L^p(\nu)$ norm.

To the best of our knowledge, the current paper is the first to
give a general result on existence of a solution to a stochastic
reaction diffusion equation (hence with coefficients of polynomial growth) with multiplicative L\'evy (or jump)
noise.

Parabolic  SPDEs driven by additive L\'evy noise were introduced
by Walsh \cite{0471.60083} and Gy{\"o}ngy  \cite{Gyongy_1982}.
Walsh,  whose motivation  came from neurophysiology,   studied   a
particular  example of the cable equation which describes  the
behaviour of voltage potentials of spatially extended neurons, see
also Tuckwell \cite{MR1002192}. Gy{\"o}ngy  considered stochastic
equations with general Hilbert space valued semimartingales
replacing the Wiener process and generalized the existence and
uniqueness theorem of Krylov and B. L. Rozovskii
\cite{Krylov+R_1979}. The question we studied  in this paper is of
similar type to those by Walsh and Tuckwell but more general as we
allow more irregular and non-additive noise. However, it is of
different type than that  by Gy{\"o}ngy: the difference is of the
similar order as between \cite{Brz+Gat_1999} and
\cite{Krylov+R_1979} in the Wiener process case.

Since the early eighties many works have been done  on the topic,
including  Albeverio, Wu, Zhang \cite{934.60055}, Applebaum and Wu
\cite{MR2002e:60099}, Bi\'e \cite{980.39765}, Hausenblas
\cite{levy2,Hausenblas_2007_levy}, Kallianpur and Xiong
\cite{766.60076,859.60050},  and Knoche \cite{MR2103204}. Some of
these papers use the framework of Poisson random measures while
others use the framework of L\'evy processes. Typically,  these
papers  deal with Lipschitz coefficients and/or Hilbert spaces and
hence none of them are applicable to the stochastic reaction
diffusion equations.

SPDEs driven by L\'evy processes in Banach spaces have  not been
intensively studied, apart from a few papers by the second author,
like  \cite{levy2,Hausenblas_2007_levy}, a very recent paper \cite{maxreg} by the first two named authours and \cite{ruediger} by Mandrekar and R\"udiger (who
actually studied  ordinary stochastic differential equations in
martingale type $2$ Banach spaces). On the other hand, Peszat and
Zabczyk \cite{Peszat_Z_2007} formulate their results in the
framework of Hilbert spaces, Zhang and R\"ockner \cite{zhang}
generalized the results of  Gy{\"ongy} \cite{Gyongy_1982} and
Pardoux \cite{Pard_1979}  by firstly, studying the evolution equations
driven by both Wiener and L\'evy processes with global Lipschitz
(and not necessarily linear) coefficients, and secondly studying
the large deviation principle and exponential integrability of the
solutions.

Martingale solution to  SPDEs driven by Levy processes in Hilbert
spaces are not often treated  in the literature. Mytnik
\cite{Mytnik_2002} constructed a weak solution to SPDEs with
non-Lipschitz coefficients driven by space time stable L\'evy
noise. Mueller, \cite{Muller_1998} studied non--Lipschitz SPDEs
driven by non-negative L\'evy noise of index $\alpha\in(0,1)$.
Mueller, Mytnik and Stan \cite{Mueller+Mytnik+Stan_2006} investigated the  heat
equation with one-sided time independent stable noise. One should add that the noise in \cite{Muller_1998} does  not satisfy the hypothesis of the current  work.

\dela{A part of the proof of the main result in \cite{Hausenblas_2007_levy} is missing.}
The proof of the main result in the paper \cite{Hausenblas_2007_levy}   contains a  gap
which has a negative impact on the correctness of the existence result stated in \cite{Hausenblas_2007_levy}. The gap is related to the Step I (see page \ref{step1}) in the current paper. In particular,
in the latter paper \cite{Hausenblas_2007_levy} only the family of laws of $\{u_n:n\in\mathbb{N}\}$ and not the pair $\{(u_n,\eta_n):n\in\mathbb{N}\}$  were considered in order to find the limit.
 In the current paper we treat with care the \dela{omitted} missing part from \cite{Hausenblas_2007_levy}. %We believe that this gap is now closed.

To our best knowledge,   at the time the first version of this paper  was prepared, the literature on nonlinear stochastic partial differential equations driven by  L\'evy processes was  rather limited. For instance Dong and Xu in  \cite{Dong+Xu_2007}
consider Burgers' equations with compound Poisson noise and thus in
fact deals with a deterministic Burgers equation on random
intervals. Some discussion of stochastic Burgers' equations with
additive L\'evy noise is contained in \cite{Brz+Z_2009}, where it
is shown how integrability properties of trajectories of the
corresponding Ornstein-Uhlenbeck process play an important r\^ole.  Since then Dong and coauthours have studied also Stochastic Navier-Stokes Equations. For instance  Dong and  Xie in a recent paper
\cite{Dong_2009}    study the
stochastic NSEs driven by a  Poisson random measure   whose L\'evy
measure  is finite. This assumption implies that the jump times
form a discrete subset of the real half line $\mathbb{R}_+$ and it
is essentially  the same as the deterministic one. The proof of
the main result is based on the approximation of  the Poisson
random measure $N$ by a sequence of Poisson random measures $N_n$
whose L\'evy measures are finite.  Fernando and Sritharan in \cite{Fernando+Sritharan_2010} and \cite{BHZ_2012},
2-d stochastic Navier-Stokes Equations are studied by a local-monotonicity method of Barbu  \cite{Barbu_2011_private} which is restricted to $2$-dimensional NSEs and does not  require  compactness.  Also Debussche at all considered in \cite{Debussche+Imkeller_2011} a stochastic reaction diffusion equation driven by \textit{additive} L\'evy noise.

\par

The paper is organized as follows. In  Section \ref{sec_analytic} we introduce the
definitions necessary to formulate the main results. The main
results are presented in  sections \ref{sec_srde} and  \ref{sec_main-results}. In section \ref{sec_srde} we consider a
SPDEs of reaction diffusion type and list the exact conditions
under which a martingale solution exists. The main results for
this section is stated in Theorem \ref{Th:general}. To prove
Theorem \ref{Th:general} we first consider auxiliary SPDEs with
continuous and bounded coefficients. The main result for these
auxiliary SPDEs are formulated either in terms of a Poisson random measure, see Theorem \ref{Th:bound}, or
 in terms of a L\'evy process, see Theorem \ref{Th:bound-levy}.
Three examples  illustrating the applicability of our main results are
then presented in Sections \ref{reaction1}, \ref{reaction1-s-t} and \ref{det-set-B}:  %\ref{reaction1} we present
an SPDE of reaction diffusion type driven by a real-valued
$\alpha$-stable tempered L\'evy noise is solved in Section
\ref{reaction1} and
an SPDE of reaction diffusion with bounded drift (reps. polynomial nonlinearity) driven by a space time L\'evy
white noise in Section \ref{reaction1-s-t} (resp. Section \ref{det-set-B}). The remaining
sections are devoted to the proofs of our
results. In the Appendix we recall some definitions and well-known results in analysis and probability theory. We also prove new results, amongst them modified version of the Skorokhod embedding theorem, see Theorem \ref{thm-Skorokhod}, that are interesting
in themselves.

Let us finish this Introduction by pointing out that the approach presented in this paper
(or rather it's earlier arXiv version) has already been taken up and used for the proof of the existence of solutions to
Stochastic Navier-Stokes equations and second grade fluids driven by L{\'e}vy noise, see \cite{Motyl_2013} and \cite{EH+PR+MS}, respectively.

\subsection*{Acknowledgments}
E.~ Hausenblas and P.~A.~Razafimandimby are funded by the FWF-Austrian Science Fund through the projects P21622 and M1487, respectively.
The research on this paper was initiated during the visit of EH  to the University of York in October 2008.
She would like to thank the Mathematics Department at York for hospitality. An earlier version of this paper can be found on arXiv:1010.5933.
  The authors would like to thank Carl Chalk, Ela Motyl and Markus Riedle for a careful reading of the
  manuscript.

\begin{notation}\label{not-0}
By $\mathbb{N}$ we denote the set of natural numbers, i.e.
$\mathbb{N}=\{0,1,2,\cdots\}$ and by $\bar{\mathbb{N}}$ we denote
the set $\mathbb{N}\cup\{+\infty\}$. Whenever we speak about
$\mathbb{N}$ (or $\bar{\mathbb{N}}$)-valued measurable functions
we implicitly assume that the  set $\mathbb{N}$, resp. $\bar{\mathbb{N}}$, is equipped with the full
$\sigma$-field $2^\mathbb{N}$, resp.  $2^{\bar{\mathbb{N}}}$. By
$\Rb{+}$ we will denote the  interval $[0,\infty)$ and by
$\Rb{\ast}$ the set $\mathbb{R}\setminus\{0\}$. If $X$ is a
topological space, then by $\mathcal{B}(X)$ we will denote the
Borel $\sigma$-field on $X$. By $\Leb\deln{\lambda_d}$ we will denote the
Lebesgue measure on $(\mathbb{R}^d,\mathcal{B}(\mathbb{R}^d))$ or \deln{ by
$\lambda$  the Lebesgue measure on}
$(\mathbb{R},\mathcal{B}(\mathbb{R}))$.

If $(S,{\mathcal{S}})$ is a measurable space then by  $M(S)$ we denote the
set of all real valued  measures on $(S,{\mathcal{S}})$ and  by ${\mathcal{M}}(S)$
the $\sigma$-field
 on $M(S)$ generated by functions
$i_B:M(S) \ni\mu \mapsto \mu(B)\in {\mathbb{R}}$, $B\in {\mathcal{S}}$. By $M_+(S)$
we denote the set of all non-negative measures on $S$  and  by
${\mathcal{M}}_+(S)$ the $\sigma$-field
 on $M_+(S)$ generated by functions
$i_B:M_+(S) \ni\mu \mapsto \mu(B)\in {\mathbb{R}}_+$, $B\in {\mathcal{S}}$. Finally,
by $M_I( S)$ we denote the family of all
$\overline{\mathbb{N}}$-valued measures on $(S,{\mathcal{S}})$,   and  by
${\mathcal{M}}_I(S)$ the $\sigma$-field
 on $M_I(S)$ generated by functions
$i_B:M(S) \ni\mu \mapsto \mu(B)\in \bar{\mathbb{N}}$, $B\in {\mathcal{S}}$.
%
%
%If $(S,{\mathcal{S}})$ is a measurable space, then by  ${\mathcal{M}}^+(S)$ we denote the set of all positive measures on $S$,
If $(S,{\mathcal{S}})$ is a measurable space then  we will denote  by
${\mathcal{S}}\otimes \mathcal{B}({\mathbb{R}_+})$  the product
$\sigma$-field on $S\times \mathbb{R}_+$ and by $\nu\otimes
\leb$ the product measure of $\nu$ and the Lebesgue measure
$\Leb\deln{\lambda}$.

For an Banach space $Y$ we denote by $\mathbb{D}([0,T],Y)$ the space of all
c\`{a}dl\`{a}g functions $u:[0,T]\to Y$. We equip $\mathbb{D}([0,T],Y)$ with the Skorokhod topology.

The space of bounded linear operators from a Banach space $Y_1$ into another Banach space $Y_2$ will be denoted by $\mathcal{L}(Y_1,Y_2)$. The norm of $A\in\mathcal{L}(Y_1,Y_2)$
 is denoted by $\lVert A\rVert_{\mathcal{L}(Y_1,Y_2)}$.

For a bounded domain $\mathcal{O}\subset \mathbb{R}^d$ with smooth boundary we denote by $C_0(\bar{ \mathcal{O}})$ the space of continuous functions which vanish on the boundary of
$\mathcal{O}$.
\end{notation}

\section{Analytic Assumptions and Hypotheses}\label{sec_analytic}

Let us begin with a list of assumptions which will be frequently
used
 throughout this  and later sections. Whenever we use any of them this
will be specifically written. Throughout this section we fix a real number $p$ satisfying
\[ p\in
(1,2]\deln{,\text{ and }  \rho\in [0,\frac 1p)}.\] Let
$E$ be a Banach space and $A$ is a linear map on $E$. The norm of
$E$ is denoted by $\lve \cdot\rve$ and the norm of any other
Banach space $Y$ is denoted by $\lve \cdot \rve_Y$. For sake of simplicity the operator norm on $\mathcal{L}(E,E)$ will be denoted by $\lVert \cdot \rVert$.
The Banach
space $E$ and the linear map $A$ satisfy the following conditions.
\begin{assumption}\label{assum-1}
\begin{enumerate}[(i)]
\item \label{h1}$E$ is a separable,  type $p$ UMD Banach
space\footnote{Note if $E$ has the $UMD$ property and of type $p$ then it is a M-type $p$ Banach space (see, for instance, \cite{MR1313905}).}.
 \item \label{h2-a}
$-A$ is a positive operator\footnote{See Section I.14.1 in Triebel's monograph \cite{Triebel_1995}.} on $E$ with compact resolvent. In particular,
there
exists $M>0$ such that
\[ \Lve
(A+\lambda)^{-1}\Rve \le \frac{M}{1+\lambda}, \mbox{ for any }\lambda\ge 0;\]
 \item \label{h2-b} $-A$ is an
infinitesimal generator of an analytic semigroup of contraction
type in $E$.
\item \label{h3} $A$ has the BIP (bounded imaginary power) property, i.e. there exist some
 constants $K>0$ and $\vt\in [0,\frac\pi2)$
such that
\begin{equation}
\Vert A^{is} \Vert \le K e^{\vartheta |s|}, \; s \in \mathbb{R}.
\label{2.1}
\end{equation}

\end{enumerate}
\end{assumption}

 For
any $\gamma>0$, the completion
of $E$ with respect to the norm $\lVert A^{-\gamma}\cdot \rVert$ will be  denoted by $D(A^{-\gamma})$.

 Let us now formulate some consequences of the
last assumption. We begin with recalling a result from
\cite{Chalk-06}.
\begin{lemma}\label{lemma-carl}
Suppose that a linear operator $A$ in a Banach space $\E $
satisfies Assumption \ref{assum-1}-(\ref{h2-a},
\ref{h2-b},\ref{h3}). Then, if $\mu \geq 0$, the operator $\mu
I+A$ satisfies these conditions as well. Assumption
\ref{assum-1}.(\ref{h2-a}) is satisfied with the same constant
$M$. Also, there exists a constant $\tilde{K}$ such that for each
$\mu \geq 0$,
\begin{equation}
\Vert (\mu I+A)^{is} \Vert \leq \tilde{K}  e^{\vartheta |s|}, \; s
\in \mathbb{R}. \label{ineq-bip-2}
\end{equation}
\end{lemma}

\begin{theorem}\label{th-domains}
Suppose that a linear operator $A$ in a Banach space $\E $
satisfies Assumption \ref{assum-1}.(\ref{h2-a},
\ref{h2-b},\ref{h3}). Then there exists a constant $\tilde{K}$
such that for all $\mu \geq 0$ and $\alpha \in [0,1]$,
\begin{eqnarray}
\Vert (\mu I+A)^{\alpha}A^{-\alpha} \Vert &\leq& \tilde{K}^2 (\mu
M+1)^\alpha ,  \label{ineq-domains}
\\
\Vert A^{\alpha}(\mu I+A)^{-\alpha} \Vert &\leq& \tilde{K}^2
\left(1+\frac{\mu M}{1+\mu}\right)^\alpha \leq \tilde{K}^2
\left(1+M\right)^\alpha ,
\\
\Vert (\mu I+A)^{-\alpha} \Vert &\leq&
\left(\frac{M}{1+\mu}\right)^\alpha.
\end{eqnarray}
\end{theorem}
\begin{proof}
Let us fix  $\mu \geq 0$ and $\alpha \in [0,1]$. Consider an
$\mathcal{L}(\E,\E)$-valued  function $f(z)$, $0\leq \Re z \leq 1$, defined
by
$$f(z)=(\mu I+A)^{z}A^{-z}.$$ This function is continuous in the
closed strip $0\leq \Re z \leq 1$ and analytic in its interior.
From the assumptions we infer that
\begin{eqnarray} \Vert f(is)\Vert &\leq& \Vert (\mu I+A)^{is}\Vert \Vert A^{-is}
\Vert \leq \tilde{K}^2 e^{2\vartheta |s|}, \; s \in \mathbb{R};\\
\nonumber \Vert f(1+is)\Vert &\leq& \Vert (\mu I+A)^{is}\Vert \Vert
A^{-is} \Vert \Vert (\mu I+A)A^{-1}\Vert \\ &\leq&
 \tilde{K}^2 (\mu M+1) e^{2\vartheta |s|}, \; s \in \mathbb{R}.
 \end{eqnarray}
 Therefore, inequality (\ref{ineq-domains}) follows by applying the Hadamard Three
 Line Theorem, see e.g. Appendix to IX.4 in \cite{MR0493420}. We can prove the two other inequalities in an analogous way.
\end{proof}

\begin{lemma}\label{lem-increament}
Assume that $-A$ is an infinitesimal generator of an analytic
semigroup $\{e^{-tA}\}_{t\geq 0}$ on a Banach space $\E $. Assume
that\footnote{compare
with Assumption II.6.1 in \cite{516.47023}.} for an $\omega \in (0,\frac\pi2)$, $M>0$  and $r>0$,
\begin{eqnarray}\label{Pazy-6.1}
\Sigma^-&:=&\{z \in \mathbb{C}: |\Arg z| \leq \pi -\omega\}\cup
B_\mathbb{C}(0,r) \subset \rho(-A), \end{eqnarray}
 where $B_\mathbb{C}(0,r)=\{z\in\mathbb{C}: \vert z \vert < r\}$, $\rho(-A)$ is the
resolvent set of the operator $-A$, and
\begin{eqnarray}\label{Pazy-6.2}
\Vert(A+z)^{-1} \Vert &\leq& \frac{M}{ 1+|z| }, \;\; z\in
\Sigma^-.
\end{eqnarray}
Then, there exist a constant $\delta>0$ and constants
$M_\alpha>0$ for $\alpha \geq 0$, such that
\begin{eqnarray}\label{Pazy-6.6}
\Vert A^\alpha e^{-tA} \Vert &\leq& \frac{M_\alpha}{t^\alpha
}e^{-\delta t}, \;\; t>0; \\
\label{Pazy-6.60} \Vert A^{-\alpha}( e^{-tA} -I) \Vert &\leq&
{M_\alpha}{t^\alpha }e^{-\delta t}, \;\; t>0.
\end{eqnarray}
\end{lemma}
\begin{proof}
For the first part see Theorem 6.13(c) in \cite{516.47023}.
If $\alpha \leq 1$, for the second part see Theorem 6.13(c)
in \cite{516.47023}. The general case follows by induction on $[\alpha]$, the
integer part of $\alpha$.
\end{proof}

\section{Martingale Solutions of Stochastic Reaction-Diffusion Equations}\label{sec_srde}

This section is devoted to the statement of our main results. We will
introduce our problem, the concept of martingale solution, the
main assumptions and state our main result. Let us recall that we assume that the Banach space $\E$ and the linear operator
 $A$ satisfy Assumption
\ref{assum-1}. We will also need  two other
Banach spaces $X$ and $B$ such that the embeddings
\begin{equation}\label{eqn-spaces}
E\subset
 X\subset B
 \end{equation}
are continuous and dense and  the following condition:
\begin{assumption}\label{assum-main-1}
The linear map $A$ has a unique extension which is still denoted by $A$ and satisfies Assumption
\ref{assum-1}-(\ref{h2-a}, \ref{h2-b}) on $X$. Furthermore, we assume that the extension $-A$ generates a
$C_0$-semigroup on $X$.
\end{assumption}

Throughout we fix $T>0$ and we consider the following
stochastic evolution equations: \begin{equation} \label{eqn-3.1}
\left\{ \begin{array}{rcl}
du(t)&+&Au(t)\, dt=F(t,u(t))\,dt
\\
&&\hspace{1truecm}\lefteqn{
+\int_{Z} G(t,u(t);z)\,\tilde \eta(dz,dt),\quad t\in (0,T]\phantom{\Bigg|}
}\\
u(0)&=&u_0,\end{array}\right.  \end{equation}
where $\tilde{\eta}$ is a compensated Poisson random measure corresponding to the measure space $(Z, {\mathcal{Z}}, \nu)$ about which we make
following standing assumption
\begin{assumption}\label{ass-Z-standing}
We assume that  $\nu$ is a $\sigma$--finite {nonnegative} measure
on a measurable space   $(Z,{\mathcal{Z}})$, {i.e. $\nu \in M_+(Z)$.}
\end{assumption}

The assumption on the nonlinear map  $G$ is   as
follows.
\begin{assumption}\label{assum-main-2}
There exists  $\rho\in [0,\frac 1p)$ such that
$A^{\rho-\frac 1p}G: [0,T]\times X \to L^p(Z,\nu,E)$ is bounded and
separately continuous.
\end{assumption}

 Next we will present assumptions on the nonlinear part $F$ of the drift operator. For this purpose we first recall the notion of subdifferential of a norm
$\vp$, for more details see \cite{MR0500309}. Given $x,y \in X$
the map $\vp :\mathbb{R} \ni s \mapsto \vert x+sy\vert \in \mathbb{R}$ is
convex and therefore it is right and left differentiable. Denote
by $D_{\pm}\vert x\vert  y $ the right/left derivative of $\vp$ at
$0$. Then the subdifferential $\partial \vert x \vert $ of $\vert
x \vert$, $x \in X$, is defined by
\[
\partial \vert x \vert := \left\{ x^\ast \in X^\ast: D_{-}\vert x \vert y
\le \langle y, x^\ast \rangle \le D_{+}\vert x \vert y, \;\; y
\in X\right\},
\]
where  $X^\ast$ is the  dual space to $X$. One can show that not
only $\partial \vert x \vert$ is a nonempty, closed and convex
set, but also
\[
\partial \vert x \vert =\{ x^\ast \in X^\ast:
\langle x, x^\ast \rangle =\vert x \vert \; \;\mbox{\rm and }
\vert x^\ast \vert \le 1\}.
\]
In particular, $\partial |0|$ is the unit ball  in $X^\ast$.

\begin{assumption}\label{assum-main-3}

\begin{enumerate}[(i)]
\item The map $F:[0,T]\times X \to X$ is  separately continuous.
\item \label{ii}   There exist numbers
 $k_0>0$, $q>1$ and $k\ge 0$  such that with $$a(r)=k_0(1+r^q), r\ge 0,$$
 the following condition holds  for  $t\in [0,T]$ ,
\begin{equation*}
\langle -Ax+F(t,x+y), z\rangle \le a(\lve y\rve_X)-k\lve x\rve_X,\;\; x\in D(A), y\in X, z\in \partial \lve
x\rve.
\end{equation*}
 \item \label{iii} There exists a sequence $(F_n)_{n=1}^\infty$ of \textbf{bounded} nonlinear maps  from $[0,T]\times X$ to $X$ such that
 \begin{enumerate}
\item $F_n$ satisfies item \eqref{ii}
 uniformly in $n$,
 \item\label{assum-b} $F_n$ converges pointwise to $F$ in $X$.
 \end{enumerate}
\end{enumerate}
 \end{assumption}

 With all the notations and concepts given above we can define a martingale solution
to Problem  (\ref{eqn-3.1})

\begin{definition}\label{Def:mart-sol}
Let $p\in (1,2]$, $E$ and $B$ are separable, UMD and type $p$   Banach spaces. Let $X$ be another
Banach space such that the embeddings $E\subset X\subset B$ are continuous. Let $(Z,{\mathcal{Z}})$ be a measurable space and  $\nu \in M_+(Z)$.

An {\sl $X$--valued  martingale solution}   to the Problem
(\ref{eqn-3.1}) is a system

\begin{equation}
\left(\Omega ,{{\mathcal{F}}},\mathbb{P},{\mathbb{F}},
\eta, \,u\right)
\label{mart-system}
\end{equation}
such that
\begin{numlistn}
\item  $(\Omega ,{{\mathcal{F}}},{\mathbb{F}},\mathbb{P})$ is a complete filtered
probability space with filtration ${\mathbb{F}}=\{{{\mathcal{F}}}_t:t\in [0,T]\}$,
\item $\eta={\{\eta(t):t\in [0,T]\}}$  is a time homogeneous Poisson
random measure on $(Z,{\mathcal{B}}(Z))$ with intensity measure $\nu$ over
$(\Omega ,{{\mathcal{F}}},{\mathbb{F}},\mathbb{P})$,
\item $u:[0,T]\times \Omega \to X$  is an ${\mathbb{F}}$-progressively
measurable %and c\`{a}dl\`{a}g
 process such that for any $t\in [0,T]$, $\mathbb{P}$ a.s.
\begin{eqnarray}\label{finitness-2}
 \int_0^t \left| e^{-(t-r)A}F(r,u(r))\right|_{X}\, dr<\infty \\
 \int_0^t \int_Z \left|
e^{-(t-r)A} G(r,u(r);z)\right|_{{\BE}} ^ p\nu(dz)\, dr  <\infty, \end{eqnarray}
and for any $t\in [0,T]$, $\mathbb{P}$ a.s.,
\begin{eqnarray}
u(t) &=&
 e^{-tA}u_0 + \int^{t}_{0}e^{-(t-r)A}F(r,u(r))\, dr
\\ &&{} +
\int^{t}_{0}\int_Ze^{-(t-r)A}\;  G(r,u(r);z)\,\tilde  \eta(dz,dr).
\nonumber
\end{eqnarray}
\end{numlistn}
If in addition, $\mathbb{P}$-a.s. $u\in\mathbb{D}([0,T];B)$, then the system \eqref{mart-system} will be called a
 an $X$--valued mild solution with \cadlag paths in $B$ to problem \ref{eqn-3.1}.
We say  that the $X$--valued, $B$-valued \cadlag martingale  solution \ref{mart-system}
to  problem \ref{eqn-3.1} is unique iff given  another martingale solution
to \ref{eqn-3.1}
\begin{eqnarray*}
&&\left(\Omega^\prime ,{{\mathcal{F}}^\prime},\mathbb{P}^\prime,{{\mathbb{F}}^\prime} ,\eta^\prime,u^\prime\right) ,
\end{eqnarray*}
the laws of the processes $u$ and $u^\prime$ on the space $\mathbb{D}([0,T];B)$
are equal.
\end{definition}

Let us denote by (\textbf{P}) the following statement.
\begin{equation}
-A \textit{ generates a contraction type semigroup on } D(A^{\rho-\frac 1p}). \tag{\textbf{P}}
\end{equation}

\noindent Let $\delta= \frac1p - \rho$ and $q_{\max}$ be a number defined by
\begin{equation}\label{qu-max}
 q_{\max}=\begin{cases}
           \frac 1\delta \text{ if the statement (\textbf{P}) is true },\\
           \frac p{1-\rho} \text{ otherwise.}
          \end{cases}
\end{equation}
\begin{remark}
 Note that in both case $q_{\max}\ge p$ and hence $q_{\max}>1$.
 Moreover, in the 2nd case $$q_{\max}<\frac{p^2}{p-1}$$ while
 in the 1st case  $q_{\max}<\infty$.
\end{remark}
Now we formulate our main theorem.

\begin{theorem}\label{Th:general}
Let $E$ be a separable, UMD and type $p \in (1,2]$ Banach space and $\rho \in [0,\frac1p)$.
We assume that
$A$ is a positive operator on $E$, having the BIP property and with compact resolvent. We suppose that $-A$ generates a contraction type $C_0$-semigroup on $E$.\\
We assume that the Banach space $X$ satisfies $E\subset X \subset D(A^{\rho-\frac1p})$ and that $-A$ admits a unique extension, denoted with the same symbol, on $X$ and this
extension is the infinitesimal generator of a $C_0$-semigroup in $X$.\\
 Let $(Z,{\mathcal{Z}}, \nu)$ be a measure space with $\nu \in M_+(Z)$ and $G$ be a map defined on $[0,T]\times X$
 such that $A^{\rho - \frac1p} G$ is $L^p(Z, \nu, E)$-valued, bounded and separately continuous. \\
 Let  $F: [0,T]\times X \to X$
 be a separately continuous map such that there exist $k_0>0, k\ge 0$ and
 $q\in (1, \infty)$ such that
 \begin{equation}\label{Dissip-F}
\langle -Ax+F(t,x+y), z\rangle \le k_0 (1+\lve y\rve_X^q)-k\lve x\rve_X,
\end{equation}
for all $x\in D(A)$, $y\in X$ and  $z\in \partial \lve x\rve$. \\
 We also suppose that there exists a sequence $(F_n)_{n\in \mathbb{N}}$ of bounded and separately continuous maps from $[0,T]\times X$ into $X$, satisfying
 \eqref{Dissip-F} uniformly in $n$, and
 pointwise converging to $F$.

Furthermore, let the number $q_{\max}$ be defined as in \eqref{qu-max}. Assume that $q<q_{\max}$ and that there exists a UMD and separable space $Y$ such that :
 \begin{enumerate}[(1)]
  \item $X\subset Y$,
  \item $A$ has a unique extension $A_Y$ which satisfies on $Y$ the same properties as satisfied by $A$ on
  $E$,
    \item and $D(A_Y^\theta)\subset X$ for some $\theta \le 1-\frac{q}{q_{\max}}$.
 \end{enumerate}
Then for any $u_0\in X$ problem \eqref{eqn-3.1} has an $X$-valued, $B$-valued \cadlag
martingale solution
in the sense of Definition \ref{Def:mart-sol} with $B=D(A^{\rho-\frac1p})$. Moreover, for any $\tilde{q}\in (q, q_{\max})$ and $r\in (1,p)$ the
stochastic process $u$ satisfies
$${\mathbb{E}}\biggl(\int_0^T \lve u(t)\rve^{\tilde{q}}_X dt\biggr)^{\frac{r}{\tilde{q}^2}} <\infty.
$$

\end{theorem}

\begin{remark}\label{rem-M type-p}
 Since $E$ is a separable, UMD and M-type $p$ Banach space, we infer from
 \cite[Remark 4.2, also Theorem A.4]{MR1313905}  that $B_0=D(A^{\beta})$, $\beta \in \mathbb{R}_\ast$, is also a UMD, M-type $p$ Banach space.
 The linear map $A$ has an extension (or restriction depending whether $\beta$ is smaller or larger than $0$) $A_0$ (usually denoted by $A$) to $B_0$ which satisfies
Assumption \ref{assum-1}-(\ref{h2-a}, \ref{h2-b}, \ref{h3}).  The operator $-A_0$
 generates a contraction type $C_0$-semigroup denoted by $\big\{e^{-tA}\big\}_{t\geq 0}$ on $B_0$. Moreover, if $A$ has the BIP property then so is $A_0$.
\end{remark}

\begin{remark}{\rm
If  the maps $F(t,\cdot)$ and $A ^ {\rho-\frac 1p }G(t,\cdot)$ are Lipschitz continuous uniformly with
respect to $t\in[0,T]$, i.e. there exists $K>0$
such that for all  $t\in[0,T]$ and all  $ u_1,u_2\in X$,
\begin{eqnarray*}
|F(t,u_2)-F(t,u_1) |_{X}
 &\leq& K|u_2-u_1|_X,
\\
 \int_Z |A ^ {\rho-\frac1p}G(t,u_2;z)-A ^ {\rho-\frac1p}G(t,u_1;z) |_{E} ^ p\:
\nu(dz) &\leq& K|u_2-u_1|_X^p,
\end{eqnarray*}
then  the SPDEs \eqref{eqn-3.1} has  a unique strong solution. In
our work we are interested in the   case when  both these
conditions are relaxed. }\end{remark}

In order to prove Theorem \ref{Th:general} we will consider an
auxiliary problem for which we will prove an existence result that
is stated in Theorem
\ref{Th:bound}. This secondary existence theorem holds under more
restrictive conditions than the ones stated above, we mainly
assume boundedness of $F$. Thanks to Assumption
\ref{assum-main-3}-(\ref{iii}.b) we can recover our main result
from Theorem \ref{Th:bound}. We will also see in the next section
that Theorem \ref{Th:general} can be reformulated in term of
L\'evy processes (see Theorem \ref{Th:bound-levy} of Section
\ref{sec_main-results}.)

\section{An auxiliary existence result}\label{sec_main-results}
%%%%%%%%%%
%
In Section \ref{sec_srde} we considered problem \eqref{eqn-3.1}
with a nonlinearity $F$ which is possibly unbounded. We mainly assume that $F$ is of polynomial growth
and the drift $-A+F$ is dissipative, this is typical to
Reaction-Diffusion equations. In the present section we will
strengthen those conditions by assuming that $F$ is bounded.

As in the previous section let $(Z,{\mathcal{Z}}, \nu)$ be a measure space with $\nu \in M_+(Z)$. Throughout we fix $T>0$ and we consider the following
stochastic evolution equation: \begin{equation} \label{eqn-3.1-B}
\left\{ \begin{array}{rcl}
du(t)&+&Au(t)\, dt=F(t,u(t))\,dt
\\
&&\hspace{1truecm}\lefteqn{
+\int_{Z} G(t,u(t);z)\,\tilde \eta(dz,dt),\quad t\in (0,T],\phantom{\Bigg|}
}\\
u(0)&=&u_0,\end{array}\right.  \end{equation}
where $\tilde{\eta}$ is the compensated Poisson random measure (to be constructed) corresponding to $\nu$.
As in Theorem \ref{Th:general} we assume that we are given a separable, UMD and M-type $p$ (with $p\in (1,2]$) Banach space $E$.
We also suppose that $A$ satisfies the assumption of Theorem \ref{Th:general}.

For the nonlinear map $F$ we assume that it satisfies the following set of conditions.
\begin{assumption}
\label{assum-2}
Assume that there exists $\rho\in [0,\frac1p)$such that the function $F$ (resp. $G$) is a \textbf{bounded} nonlinear map
from  $[0,T]\times E$ into $D(A^{\rho-1})$ (resp. $D(A^{\rho-\frac 1p}$)) and separately
continuous w.r.t its variables.
\end{assumption}

\begin{assumption}\label{assum-3}
Assume that $u_0\in D(A^{\rho-\frac
1p})$.
\end{assumption}

We state the following theorem whose proof will be given in
Section \ref{proof-aux-result}. Even if it is our secondary result
it is still important as it is the main tool for the proof of
Theorem \ref{Th:general}.
\begin{theorem}
\label{Th:bound}
Let $p\in (1,2]$, $\rho\in [0,\frac1p)$ and $E$ be a separable, UMD and M-type $p$ Banach space. We assume that
$A$ is a positive operator on E, having the BIP property and with compact resolvent. We suppose that $-A$ generates a contraction type $C_0$-semigroup on $E$.
Let  $F: [0,T]\times E \to D(A^{\rho-1})$
 be a bounded and separately continuous map. Let $(Z,{\mathcal{Z}}, \nu)$ be a measure space with $\nu \in M_+(Z)$. We assume that $A^{\rho - \frac1p} G$ defined on $[0,T]\times E$
is $L^p(Z, \nu, E)$-valued, bounded and separately continuous.

 Then for every $u_0\in D(A^{\rho-\frac1p})$ problem \eqref{eqn-3.1-B} has an $E$-valued martingale solution which has \cadlag
paths in $B=D(A^{\rho-\frac1p})$. Moreover, the stochastic process
$u$ satisfies
\begin{eqnarray} \label{int-est-theo} \int_0 ^T {\mathbb{E}} |u(t)|_{E}^p\: dt<\infty.
\end{eqnarray}

%%%%%%%%%%%%%%%%%%%%%%%%%%%%%%%%%%%%%%%%%%%%%%%%
 \end{theorem}
 As pointed out in Subsection \ref{subsection-rel-l-prm},  a
Banach space valued L\'evy process  induces a time homogeneous Poisson
random measure on the same  Banach space and vice versa. Hence,
Theorem \ref{Th:bound} can be written in terms of a L\'evy
process. Let $Z$ be a Banach space, $L=\{L(t):t\ge 0\}$ a
$Z$--valued L\'evy process with finite $p$-variation. Here we keep all the assumptions on $A, B, E, F$  and assume that
the diffusion coefficient $\mathbb{G}$ is
such that $A^{\rho-\frac 1p}\mathbb{G}$ maps ${\mathbb{R}}_+\times E$ into
the space $\mathcal{L}(Z,E)$ of bounded linear operators from $Z$ to $E$.
In this sense  the following framework is less general than the Poisson Random Measure's setting
: \begin{eqnarray} \label{eqn-3.1-levy} \hspace{0.8truecm} \left\{
\begin{array}{rcl}
du(t)&+&Au(t)\, dt=F(t,u(t))\,dt+ \mathbb{G}(t,u(t))\,dL(t),\phantom{\Bigg|}\\
u(0)&=&u_0.\end{array}\right.  \end{eqnarray}

In order to formulate a version of Theorem \ref{Th:bound} in terms
of the L\'evy process $L$ we need to reformulate Assumption
\ref{assum-main-2}.

\begin{assumption}\label{ass-levy}
Let $Z$ be a Banach space and $\rho\in [0,\frac 1p)$ such that the
function $A ^ {\rho-\frac1p}\mathbb{G}:[0,T]\times E \to
\mathcal{L}(Z,E) $ can be seen as a map
$$A ^ {\rho-\frac 1p}\mathbb{G}:[0,T]\times E \to L^p (Z,\nu;E) $$
 satisfying Assumption \ref{assum-main-2}.
\end{assumption}

\begin{theorem}
\label{Th:bound-levy} Assume that $p\in (1,2]$, $\rho\in [0,\frac1p)$, $E$ and $A$ satisfy the assumptions of Theorem \ref{Th:bound}.
%$e^{-tA}$ acts on $B$.
Let $Z$ be a Banach space of type $p$ and  $\nu\in{\mathcal{M}}_+(Z)$ a L\'evy measure on $(Z,{\mathcal{Z}})$.

Assume that the functions $F$  and $\mathbb{G}$ satisfy
Assumptions \ref{assum-2} and \ref{ass-levy}, respectively.
Assume, that $u_0$ satisfy Assumption \ref{assum-3}.

%\\[0.2cm]
Then,  there exists a system
$$
\left(\Omega ,{{\mathcal{F}}},\mathbb{P},{\mathbb{F}},
\{L(t):t\in [0,T]\},\{u(t):t\in [0,T]\}\right) $$
% \]
such that
\begin{numlistn}
\item  $(\Omega ,{{\mathcal{F}}},{\mathbb{F}} ,
\mathbb{P})$ is a complete filtered probability space with
filtration ${\mathbb{F}} =\{{{\mathcal{F}}}_t\}_{t\in [0,T]}$; \item ${\{L(t):t\in
[0,T]\}}$ is a $Z$-valued L\'evy process with characteristic
measure $\nu$
 over $(\Omega
,{{\mathcal{F}}},{\mathbb{F}} ,\mathbb{P})$; \item $u=\{u(t): t\in [0,T]\}$  is a
$E$-valued and adapted  process, \cadlag  on $B$, such that
$(F(t,u(t)):t\in [0,T])$ and $(\mathbb{G}(t,u(t)): t\in [0,T])$ are
well defined $B$-valued, resp.\  $L^p (Z,\nu;B)$ valued and
progressively measurable processes, and for all $t\in [0,T]$,
$\mathbb{P}$-a.s.
\begin{eqnarray}
u(t) &=&
 e^{-tA}u_0 + \int^{t}_{0}e^{-(t-r)A}F(r,u(r))\, dr
 \nonumber \\&& \quad +
\int^{t}_{0}e^{-(t-r)A}\;  \mathbb{G}(r,u(r))\,dL(t) . \nonumber
\end{eqnarray}
\end{numlistn}
%
%\ref{eqn-3.1}. ,\xi\in\mathcal{O}
Moreover, \begin{eqnarray} \label{int-est-theo-1} \int_0 ^T{\mathbb{E}} |u(t)|_E ^
p\: dt<\infty. \end{eqnarray}
\end{theorem}

\begin{remark}\label{rem-ass-levy} A sufficient condition for Assumption \ref{ass-levy} to hold is that
the function
\[A ^ {\rho-\frac 1p}\mathbb{G}:\mathbb{R}_+\times E \to \mathcal{L}(Z,E) \]
 satisfies Assumption \ref{assum-main-2}.
\end{remark}
The proof of Theorem \ref{Th:bound} will be given in Section
\ref{proof-aux-result}. In the next three sections we give some applications of Theorem \ref{Th:general} and Theorem \ref{Th:bound}.

\section{Application I: The Reaction-Diffusion Equation with L\'evy Noise of Spectral Type}
\label{reaction1}
\nopagebreak

Throughout this section, ${\mathcal{O}}$ is a bounded open domain in
${{\mathbb{R}}}^d $, $d \ge 1$, with $C ^ \infty$ boundary. Let $T>0$, $\alpha\in
(0,2)$, $p\in(1,2]$ and $q\in (1,\infty)$ be fixed
real numbers. Let $L=\{L(t): t\in [0,\infty)\}$ be a real-valued
tempered $\alpha$-stable L\'evy process, i.e., a L\'evy process
with   the characteristic measure $\nu_\alpha$ given by
\begin{equation}
\label{eqn-nu_R}
\nu_\alpha(dz)=\lvert z\rvert^{-\alpha -1}e^{-\lvert z\rvert}\,dz,\quad z\in{\mathbb{R}}\setminus\{0\}.
\end{equation}
Our aim in this section is to specify the conditions under which  Theorem
\ref{Th:general} covers an equation of the following type, where $T>0$ is fixed,
\begin{eqnarray}\label{spec-type}
\\
\nonumber&& \left\{
\begin{array}{rcl} du(t,\xi) &=& \Delta u(t,\xi ) \: dt -|u(t,\xi)| ^{q-1}\sgn(u(t,\xi)) + u(t,\xi) \: dt
\\
&+&  \sqrt{|u(t,\xi)|}/(1+\sqrt{|u(t,\xi)|}) \,dL (t) ,\quad t\in (0,T],\, \xi\in\mathcal{O},\phantom{\Big|}
\\
u(t,\xi) &=& 0,\quad \xi\in\partial \mathcal{O},\,
t\in(0,T],\phantom{\Big|}
\\
u(0,\xi) &=& u_0(\xi),\quad \xi\in\mathcal{O},\phantom{\Big|}
\end{array}\right..
\end{eqnarray}
For this purpose we will reformulate problem \eqref{spec-type}
using  a more general setting and the language of the Poisson random measures.

Firstly, we denote by $X=C_0(\mathcal{O})$ the space of continuous functions $u: \bar{\mathcal{O}}\to\mathbb{R}$ which
vanish on the boundary $ \partial{\mathcal{O}}$. For $\gamma \in
\mathbb{R}$ and $r\in (1,\infty)$  the symbol $H^{\gamma,
r}(\mathcal{O})$ and $H^{\gamma,r}_0(\mathcal{O})$ denote the Sobolev spaces as defined in \cite[Definition 1, page 301]{Triebel_1995} and \cite[Definition 2, page 301]{Triebel_1995},
respectively. We just write $H^{\gamma, r}$ and $H_0^{\gamma,r}$
when there is no risk of ambiguity.  Let us briefly recall definitions of these spaces.
If $k$ is  a natural number  and $p\in [1,\infty)$ is  a real number,   we denote by $H^{k,p}(\mathcal{O})$,  see
\cite[section I.6]{Friedman_1969}, the space of all functions $u \in L^{p}(\mathcal{O})$ whose weak derivatives $D^\alpha u$ of degree $\vert \alpha\vert \leq k$ exist and belong to $L^{p}(\mathcal{O})$.  Endowed with a natural norm $\Vert \cdot \Vert_{k,p}$
\[
\Vert u \Vert_{k,p}^p :=\sum_{\vert \gamma\vert \leq k} \vert D^\gamma u\vert_{L^p}^p,\;\; u \in  H^{k,p}(\mathcal{O}),
\]
this space is a separable Banach space.  The closure of the space $C_0^\infty (\mathcal{O})$ in the space $H^{k,p}(\mathcal{O})$
is  denoted  by $H_0^{k,p}(\mathcal{O})$.
In the case $\beta \in \mathbb{R}^+\setminus \mathbb{N}$,  the fractional order Sobolev spaces $H^{\beta,p}(\mathcal{O})$
 can be defined by the complex interpolation method, i.e.
\begin{equation}\label{eqn-sobolev_spaces}
H^{\beta,p}(\mathcal{O})=[H^{k,p}(\mathcal{O}), H^{m,p}(\mathcal{O})]_{\vartheta},
\end{equation}
 where
$k, m \in \mathbb{N}, \vartheta \in (0,1)$, $k<m$, are chosen to
satisfy $ \beta=(1-\vartheta)k+\vartheta m$.   It is well known, see e.g. Theorem 11.1 in
\cite[chapter I]{Lions+Magenes_1972} and Theorem 1.4.3.2 on p.317 in \cite{Triebel_1995}   that
$H_0^{s,p}(\mathcal{O})=H^{s,p}(\mathcal{O})$ iff $s\leq \frac1{p}$.

We denote by  $A$ the negative of the Laplace operator in $\mathcal{O}$  with the Dirichlet
boundary conditions.

Secondly let us consider a separately continuous real valued
functions $f$ defined on $[0,T]\times \mathcal{O}\times \mathbb{R}$
satisfying the following condition. There exists a number $K>0$
such that for $t\ge0, x \in {\mathcal{O}}, u \in {{\mathbb{R}}}$
\begin{equation}
-K(1+\vert u\vert^q1_{[0,\infty)}(u))\le f(t,x ,u)\le K(1+\vert
u\vert^q1_{(-\infty, 0]}(u)).
\label{eq:f-diss}
\end{equation}
 It is not difficult to prove that  if $f$ satisfies
\ref{eq:f-diss} then  $$f(t,v+z)\sgn(v) \le K(1+|z|^q), \text{ for
all }v,z \in \mathbb{R}, t\in [0,T].$$ Therefore, by
\cite[Proposition 6.2]{Brz+Gat_1999} the Nemytskii map $F$ defined
by
\begin{equation}\label{mmmm-1}
F(t,u)(\xi) :=  f(t,\xi,u(\xi)), \; u \in X, \xi\in\mathcal{O}, t\in [0,T],
\end{equation}
satisfies Assumption \ref{assum-main-3}-(i, \ref{ii}) on $X $.\\
Approximating %the function
$f$ by a sequence $(f_n)_{n\in{\mathbb{N}}}$, where for any $t\in [0,T]$
and $\xi\in \mathcal{O}$
$$
f_n(t,\xi,u) := \bcase f(t,\xi  ,u), & \mbox{ if } u\in[-n,n]\\  f(t,\xi
,n), & \mbox{ if } u\ge n,\\  f(t,\xi ,-n), & \mbox{ if } u\le - n,\ecase
$$%
we obtain a  sequence $(F_n)_{n\in{\mathbb{N}}}$ defined by
$$
F_n: [0,T]\times X \ni (t,u) \mapsto
\{ \mathcal{O} \ni x \mapsto  f_n(t,\xi,u(\xi) ) \}\in X,
$$
which satisfies Assumption \ref{assum-main-3}\eqref{iii}.

Next we reformulate the noise appearing in the problem we want to study. %%%%%%%%%%%
%%%%%%%%%%%%%%%%%%%%%%%%%%%%%%%%%%
Let $\eta$ be a time homogeneous  Poisson random measure on ${\mathbb{R}}$
over a complete filtered probability space $(\Omega,{\mathcal{F}},{\mathbb{F}},\mathbb{P})$
with the intensity measure $ \nu_\alpha$ defined in
\eqref{eqn-nu_R} and let $\tilde \eta$ be the corresponding
compensated Poisson random measure. It is well known that the L\'evy processes $L=\{L(t): t\ge 0\}$ and $\tilde L=\{\tilde L(t):
t\ge 0\}$, where
\[
 \tilde L(t):=\int_0 ^ t \int_{\mathbb{R}} z \, \tilde  \eta(dz,ds), \;\; t\ge 0,
\]
have the same laws on $\mathbb{D}([0,\infty);\mathbb{R})$.
%%%%%%%%%%%%%%%%%%%%%%%

 Now let $g$ be a real-valued bounded and
separately continuous function defined on $[0,T]\times
\mathcal{O}\times\mathbb{R}$.
Let   $G$ be defined by
\begin{equation}\label{defi_G}
G(t,u;z)(\xi)= g(t,\xi,u(\xi))z,\;\; t\in [0,T], u\in L ^1 (\mathcal{O}),
z\in\mathbb{R}, \xi\in \mathcal{O}.
\end{equation}
With this notation problem \eqref{spec-type} can be rewritten
in the following form
\begin{equation}\label{spde011}
\begin{cases}
du(t) + A u(t) \: dt =F(u(t))\: dt+ {} \int_{{\mathbb{R}}}G(t,u,z)
\,\tilde{\eta}( dz,dt),\quad t\in(0,T],\phantom{\Big|}
\\
u(t) = u_0.\quad \phantom{\Big|}
\end{cases}
\end{equation} The following Theorem is a corollary of
Theorem \ref{Th:general}.
\begin{theorem}\label{ex01}
 Assume that $\alpha\in
(0,2)$, $p\in (1, 2]$, $p>\alpha$, $d\in \mathbb{N}$.
Let $\nu_\alpha$ be a L\'evy measure on $\mathbb{R}$ given by \eqref{eqn-nu_R}. Assume that $q>1$.

Then if $B=L^r(\mathcal{O})$ with $r$ satisfying
$$ r> \max\{qd, 2d \},$$
the following holds.\\
For any $u_0\in C_0(\mathcal{O})$ there exists an $C_0(\mathcal{O})$-valued, $L^r(\mathcal{O})$-valued c\`adl\`ag martingale solution $(\Omega, \mathcal{F}, \mathbb{F}, \mathbb{P}, \eta, u)$
 to problem \eqref{spde011} such that $\eta$ is a time homogeneous Poisson random measure $({\mathbb{R}}, \mathcal{B}({\mathbb{R}}))$ with intensity measure $\nu_\alpha$.

\end{theorem}
\begin{proof}
Let $\alpha$, $p$, $d$, $q$ as in the assumptions. We choose $r> \max\{qd, 2d\}$ and $\kappa \in (\frac dr, \frac 1q )$.
We put $\delta= \frac\kappa 2$. Then $\delta<\frac 1{2q}<\frac12$, and since $\frac12\le \frac1p$ we have $\delta<\frac1p$.
We also have that $q_{\max}=\frac1\delta>q$.
Next, we put $E=H^{\kappa,r}_0$ and we denote by $A=A_r$ the negative of the Laplace operator with the Dirichlet boundary conditions in the space $B$. We have $r\ge p$ because $d\ge 1$,
$r\ge 2d$ and $p \in (1,2]$. Thus, $E$ and $B$ are separable, UMD and type $p$ Banach spaces.
Now, it is well known that the assumptions of the first and second part of Theorem \ref{Th:general} are satisfied by $A_r$.

We put $Z=\mathbb{R}$ and $\nu=\nu_\alpha$. Then we immediately see that
\[
C_p:=\int_Z \vert z\vert^p \, \nu(dz)<\infty.
\]
Then we define a map $\tilde{G}$ by
\[
\tilde{G}(t,\cdot): X \ni u \mapsto \{ Z\ni z \mapsto z g\circ
u\}\in L^p(Z,\nu;E), t\in [0,T].
\]
The map $\tilde{G}$ may not be defined on the whole space $X$ but
the map $A_r^{-\delta}\tilde{G}$ is because $\delta=\frac
\kappa2$. Indeed, we have the following chain of
inequalities.
\begin{eqnarray*}
&&\hspace{-1truecm}\lefteqn{\int_Z \vert A_r^{-\delta}\tilde{G}(u)(z)\vert_E^p \, \nu(dz) = \int_Z \vert A_r^{-\delta}\left((g\circ u)z\right)\vert_E^p \, \nu(dz)}\\
&=& \int_Z \lvert A_r^{-\delta}(g\circ u)\rvert^p_E \lvert z\rvert^p
\nu(dz)
 \le  C_p \lvert A_r^{\frac{\kappa}{2}}A_r^{-\delta}
(g\circ u)\rvert^p_{L^r}
 \le  C_p\vert g\vert_{L^\infty}^p.
\end{eqnarray*}
Since the function $g$ is continuous one can easily check that the continuity condition in
Assumption \ref{assum-main-2} is satisfied. Observe that
$$\tilde{G}(t,u)(z)=G(t,u,z), t\in [0,T], u\in X, z\in Z,$$ where $G$ is defined in \eqref{defi_G}.

Let us choose $X=C_0(\mathcal{O})$. Since $\kappa>\frac dr $ we have $E\subset X $. Moreover, it is straightforward to check that the nonlinear map $F$ defined by \eqref{mmmm-1} satisfies
Assumption \ref{assum-main-3} on $X$.

Finally, let $Y=L^r(\mathcal{O})$ and $A_Y=A_r$. Since $1-\frac{q}{q_{\max}}>\frac12$ and $\frac 12 > \frac dr$, we can find $\kappa_1\in (\frac dr, \frac12 )$ such that
$$ D(A_r^{\frac{\kappa_1}{2}})\subset X \subset Y .$$

Thus all assumptions (with our choice of spaces and maps) of Theorem
\ref{Th:general} are satisfied and therefore the proof of the existence of a solution with requested properties follows.
\end{proof}

%%%%%%%%%%%%%%
%%%%%%%%%%%%%%%%%%%%%%%%%
%%%%%%%%%%%%%%%%%%%%%%%%%

\section{Application II: The Reaction-Diffusion Equation of an arbitrary Order with Space Time L\'evy Noise} \label{reaction1-s-t}

In this chapter we will apply Theorem \ref{Th:bound} to SPDEs of
reaction diffusion type driven by the so called space time L\'evy
noise or impulsive white noise. This kind of noise is a generalisation of the space time
white noise and is treated quite often in the literature, e.g.\ in
Peszat and Zabczyk \cite[Definition 7.24]{Peszat_Z_2007} or St
Lubert Bi\'e \cite{980.39765}.

First, we introduce the deterministic part of the equation, and later we will
present our result, i.e.\ the precise assumptions under which a
martingale solution of such a reaction diffusion type equation
with space time Poissonian noise exists.

\subsection{Deterministic part of the problem} \label{det-set-A}%A elliptic operator}
Let $d\ge 1$, $p\in (1,2]$, and $\mathcal{O}$ is a bounded open domain in
$\mathbb{R}^d$ with  boundary $\partial \mathcal{O}$ of ${\mathcal{C}}^\infty$
class. Let also $k$ be a positive integer. Borrowing the
presentation of \cite[Section 6.3]{Brz+Gat_1999} we introduce a
differential operator ${\mathcal{A}}$ of order $2k$ as follows.

\begin{enumerate}[(a)]
\item \label{item-a} The differential operator ${\mathcal{A}}$ defined by
\begin{equation}
{\mathcal{A}} u(x )=-\sum_{\lvert \alpha\rvert\le 2k }a_\alpha(x) D^\alpha
u(x), x\in \mathcal{O}, \label{A-high-1}
\end{equation}
is properly elliptic (see \cite[Section 4.9.1]{Triebel_1995}). The
coefficients $a_\alpha$ are $\mathcal{C}^\infty$ functions on the
closure $\bar{\mathcal{O}}$ of $\mathcal{O}$.

\item A system $\{\mathcal{B}_j: j=1,\cdots, k\}$ of differential
operators on $\partial \mathcal{O}$ is given,
\begin{equation}\label{A-high-2}
\mathcal{B}_j=\sum_{\lvert \alpha \rvert \le
m_j}b_{j,\alpha}D^\alpha,
\end{equation}
with the coefficients $b_{j,\alpha}$ being $\mathcal{C}^\infty$
functions on $\partial \mathcal{O}$. The orders $m_j$ of the operators
$\mathcal{B}_j$ are ordered in the following way:
$$0\le m_1<m_1<\ldots<m_k.$$
We assume that $m_k<2k$ and
\begin{equation}\label{A-high-3}
\sum_{\lvert \alpha\rvert=m_j}b_{j,\alpha}(\xi)n_\xi\neq 0, x\in
\mathcal{O}, j=1,2,\ldots,k,
\end{equation}
 where $n_\xi$ is the unit outer normal
vector to $\partial \mathcal{O}$ at $\xi\in \partial \mathcal{O}$.

\item For any $x\in \bar{\mathcal{O}}$ and $\xi \in \mathbb{R}^n\backslash
\{0\}$ let $a(x,\xi)=\sum_{\lvert
\alpha\rvert=2k}a_\alpha(x)\xi^\alpha$. We assume that
\begin{equation}\label{A-high-4}
(-1)^k\frac{a(x,\xi)}{\lvert a(x,\xi)\rvert}\neq -1, x\in
\bar{\mathcal{O}}, \xi \in \mathbb{R}^n\backslash \{0\}.
\end{equation}
\item \label{itewm-b} If $b_{x,\xi}=\sum_{\lvert
\alpha\rvert=2k}a_\alpha(x)\xi^\alpha$ then for all $x\in \partial
\mathcal{O}$, $\xi \in T_x(\partial \mathcal{O})$, $t\in (-\infty,0]$ the
polynomial
$$\{\tau\to b_j(x,\xi+\tau n_x)\}, j=1,\cdots,k$$
are linearly independent modulo polynomial $\{\tau\to \prod_{j=1}^k
(\tau-\tau^{+}(t)\}$. Here $T_x(\mathcal{O})$ is the set of all tangent
vectors to $\partial \mathcal{O}$ at $x\in\partial \mathcal{O}$ and
$\tau^{+}_j(t)$ are the roots with positive imaginary part of the
polynomial defined by $\mathbb{C}\ni \tau \to a(x,\xi+\tau
n_x)-t$.
\end{enumerate}
 The differential operator ${\mathcal{A}}$ induces a linear unbounded map
 $A_r$ on the Banach space $L^r(\mathcal{O})$, $r>1$, defined by
 \begin{equation}\label{A_high-5}
 \begin{cases}

D(A_r)=\{ u \in H^{2k,r}: \mathcal{B}_ju\big|_{\partial \mathcal{O}}=0
\text{ for } m_j< 2k-\frac1r\},\\
A_ru={\mathcal{A}} u, \quad u\in D(A_r).
 \end{cases}
 \end{equation}
\begin{assumption}\label{assum-Levy-stwn}
Assume that $\nu$
is a L\' evy measure on ${\mathbb{R}}$ such that
\begin{equation}\label{eqn-cond_pG}
\text{there exists } p\in (1,2]: C_p(\nu):=\int_{\mathbb{R}} |z|^p \nu(dz)
<\infty,
\end{equation}
and $L_{st}=\left\{ L_{st} (t):t\ge 0\right\}$ is a space time L\'evy
white noise  with intensity jump size measure $\nu$.

We fix $p$ as above for the remainder of this section.
We also put  $\hat \nu=\Leb\otimes \nu$,
where $\Leb$ is the Lebesgue measure on $\mathcal{O}$.
\end{assumption}
\begin{claim}\label{c7-1-high}
Assume that  $r\in (1,\infty)$ and $\theta \geq 0$.
Then we put $E=D(A_r^{\frac\theta{2k}})$, where $A_r$ is the linear  operator in the Banach space $L^{r}(\mathcal{O})$ defined in \eqref{A_high-5}. Then $A_r$ satisfies Assumption \ref{assum-1} on the space $E$.
\end{claim}
\begin{proof} It is enough to consider the case $\theta=0$.
The proof follows from Triebel's book \cite[Section 4.9.1]{Triebel_1995}, Seeley's paper \cite{MR0287376} and Lunardi's book \cite[Section
3.2]{Lunardi-book}. See also \cite[Section
6.3]{Brz+Gat_1999}.
\end{proof}
We introduce a nonlinear map which will play the role of the drift
for our stochastic equation.
\begin{assumption}\label{assum-F-high}
Assume that a function  $f;[0,T]\times \mathcal{O}\times \mathbb{R} \to \mathbb{R}$  is separately continuous and bounded. We denote by $F:
[0,T]\times L^r(\mathcal{O}) \to L^r(\mathcal{O})$ the Nemytskii map associated to $f$, i.e.
$F$ is defined by
\begin{eqnarray}\label{F-high} F(t,u)(x) & := &
f(t,x,u(x)) , \quad u \in L^r(\mathcal{O}), x\in\mathcal{O}, t\in [0,T].
\end{eqnarray}
The restriction of $F$ to $[0,T]\times D(A_r^{\frac\theta{2k}})$, $\theta\ge 0$, is also denoted by $F$.
\end{assumption}
\subsection{Coefficients of the noise}
 Now we introduce
the nonlinear coefficient of the noise.
\begin{assumption}\label{G-high}
 Let $g$ be a real bounded function defined on $[0,T]\times {\mathbb{R}}\times
\mathcal{O}$. We assume that $g$ is separately  continuous with
respect to the first and the second variables, and uniformly continuous with respect to
the third variable.
\end{assumption}

Let $\theta \geq 0$, $r\ge p$, and let us set $E=D(A_r^{\frac\theta{2k}})$.
 We define the nonlinear map $G_0:[0,T]\times L^
r(\mathcal{O})\to   L^ r(\mathcal{O})$ to be the Nemytskii operator associated to
$g$; that is,
\[
G_0(t,u) (x) := g(t,u(x),x), \quad
 \, u\in L^ r(\mathcal{O}),\; x\in\mathcal{O}.
\]
In view of the Assumption \ref{G-high}, by the Lebesgue dominated
convergence theorem, for every $t\geq 0$, $G_0(t, \cdot)$ is a
continuous map from $L^ r(\mathcal{O})$ into itself and for each
$u\in L^ r(\mathcal{O})$, the function $G_0(\cdot,u):[0,T]\to L^ r(\mathcal{O})$
is strongly measurable.

 By
Proposition \ref{prop-delta} along with Corollary \ref{cor-delta},
 we can define a bounded linear map
 \begin{equation}\label{eqn-Phi-elliptic case}
 \Phi: L^p(\mathcal{O})\to L^ p( \mathcal{O}\times {\mathbb{R}};
\hat \nu; B_{r,\infty}^ {-(d-\frac d{r})}(\mathcal{O}))
\end{equation}

 by
  \begin{equation}\label{eqn-Phi}
[\Phi v] (x, y) = (v(x)\delta_ x) y, \;\; (x,y)\in \mathcal{O}\times {\mathbb{R}}.
\end{equation}
Indeed, $\Phi$ is linear and by Corollary \ref{cor-delta} and
\eqref{eqn-cond_pG} we have the following chain of
equalities/inequalities
\begin{align*}
&\int_{\mathcal{O}\times \mathbb{R}} \lvert [\Phi
v](x,y)\rvert^p_{B^{-(d-\frac dr)}_{r,\infty}(\mathcal{O})}dx
\nu(dy)=\int_{\mathcal{O}\times \mathbb{R}} \lvert (v(x)\delta_ x)
y\rvert^p_{B^{-(d-\frac dr)}_{r,\infty}(\mathcal{O})} \hat
\nu(dx,dy)\\
=&\int_\mathcal{O}\lvert v \delta_x\rvert^p_{B^{-(d-\frac
dr)}_{r,\infty}(\mathcal{O})} dx \times \int_\mathbb{R} \lvert y\rvert^p
\nu(dy) \le C C_p(\nu)  \lvert v\rvert^p_{L^p(\mathcal{O})}.
\end{align*}
Finally, by the choise of $\theta$, $r$, and $p$ above, the embeddings $D(A_r^{\frac\theta{2k}})\subset L^r(\mathcal{O})\subset L^p(\mathcal{O})$ are continuous, so we can define  a nonlinear map $G$ as follows:
 \begin{eqnarray}
 G&:=& \Phi \circ G_0: [0,T]\times C_0(\mathcal{O}) \to
L^ p( \mathcal{O}\times {\mathbb{R}},\hat{\nu};  B_{r,\infty}^ {-(d-\frac
d{r})}(\mathcal{O})) .\label{DEF-G-high}
\end{eqnarray}
For the definition of the space $B_{r,\infty}^ {s}$ we refer the reader to Appendix \ref{besov}.

In what follows we still denote by $G$ the restriction
 of $G$ to $[0,T]\times D(A_r^{\frac\theta{2k}}) $ with $r\in (p,\infty)$ and $\theta \geq 0$.
It follows from the corresponding properties of the map $G_0$ that for every $t\geq 0$,
$G(t, \cdot)  $  is  continuous. Moreover, for each $u\in D(A_r^{\frac\theta{2k}})$ the function
$G(\cdot,u)$ is strongly measurable.
\begin{claim}\label{c7-3-high}
Assume that  $p\in (1,2]$, $d\in \mathbb{N}$, $r\ge p$, $k\in \mathbb{N}$ and $\theta\ge 0$ such that
\begin{equation}\label{neq-theta} \theta +d-\frac{d}{r}< \frac{2k}p.\end{equation}
Put $E= D(A_r^{\frac\theta{2k}})$. Then there exists $\delta<\frac1p$ such that  the map $A_r^{-\delta}G $ defined on $[0,T] \times E$ is $L^p(Z,\nu,E)$-valued, bounded and
 continuous w.r.t. $E$ and measurable with respect to time.
\end{claim}
\begin{proof}
Let us fix $k$, $r$, $d$, $\theta$ and $p$ as in the assumptions. Let us choose $\gamma>d -\frac{d}{r}$ such that   $\delta:=\frac{\theta+\gamma}{2k}<\frac1p$.
Then, since
$A_r^{-\delta}$ maps $H^{-\gamma,r}(\mathcal{O})$ into
$H^{\theta,r}(\mathcal{O})= E$ and, by \cite[Theorem
4.6.1-(a,b)]{Triebel_1995},  the Banach space $B ^ {-(d-\frac
d{r})}_{r,\infty}( \mathcal{O})$ is continuously embedded in $H^{-\gamma,r}(\mathcal{O})$, we infer
 that the  map
$A^{-\delta}G$ is $L^ p( \mathcal{O}\times {\mathbb{R}},\hat \nu;E)$-valued  continuous. Therefore, by  the continuity and boundedness
of the function $g$,  the function $A^{-\delta}G$ from
$[0,T]\times E$ into $L^ p( \mathcal{O}\times {\mathbb{R}},\hat \nu;E)$ is
separately continuous and bounded. Since $\delta<\frac1p$ we
deduce that  $G$ satisfies Assumption \ref{assum-main-2} with $\rho=\frac1p -\delta$.
\end{proof}
\begin{remark}\label{REM-c7-3-high}
 If $d< \frac{2k}p$, then we can find $\theta>\frac{d}r$ such that condition \eqref{neq-theta} is satisfied and the space $E$ is continuously embedded in $C_0(\mathcal{O})$.
\end{remark}

With the functional setting as in  Claim \ref{c7-1-high} and the
 mappings $A_r$, $F$ and $G$ defined above, the SPDEs that we
are interested in is
  \begin{equation} \label{eqn-HIGH}
 \begin{cases}
&du(t)+A_r u(t)\, dt= \int_{\mathcal{O}\times {\mathbb{R}}} G(u(t))[\xi,\zeta]\,\tilde \eta(d\xi\times d\zeta\times dt)\\
& \hspace{4truecm}+F(t,u(t))\,dt,\,\, t\in(0,T],\\
&u(0)=u_0.
\end{cases}
\end{equation}

\begin{remark}
A particular example of problem \eqref{eqn-HIGH} is the following SPDE
\begin{equation}
\begin{cases}
 \frac{\partial }{
\partial t} u(t,\xi)&+{\mathcal{A}} u(t)\, dt=f( u
(t,\xi))\phantom{\Big|}+  g(u(t ,\xi))[\dot{L}(\xi,t)] \,\,
\xi\in\mathcal{O},\, t\in(0,T],
\\
u(0,\xi) =& u_0(\xi),\quad\xi\in\mathcal{O},\phantom{\Big|}
\\
 u(t,\xi)=&
0, \hbox{ for }\xi \in \partial \mathcal{O}, \,  t\in (0,T].
\end{cases}
\end{equation}
where ${\mathcal{A}}$ is a second order differential operator, both $f$ and $g$ are continuous and bounded real
functions defined on $\mathbb{R}$ and, roughly speaking,  $\dot{L}$
denotes  the Radon-Nikodym derivative of the space time L\'evy
white noise $L$, i.e.\
$$
 \dot{L}(\xi,t) := \frac{ \partial  L (\xi,t)}{ \partial t\, \partial \xi}
.$$
\end{remark}

 We define the solution to problem  \eqref{eqn-HIGH} in
the following sense.

\begin{definition}\label{Def:mart-sol-s-t} Let $p\in (1,2]$ and $\nu$ a L\' evy measure on ${\mathbb{R}}$ satisfying condition
\eqref{eqn-cond_pG}. For $\theta \ge 0$ and $r\ge p$ we put $E=D(A_r^{\frac\theta{2k}})$, where $A_r$ is the linear  operator in the Banach space $L^{r}(\mathcal{O})$ defined
in \eqref{A_high-5}.  \delzb{be a separable Banach space of martingale type $p$.} An $E$-valued {\sl mild martingale
solution} to equation \eqref{eqn-HIGH} is a system
\begin{eqnarray} & \left(\Omega ,{{\mathcal{F}}},\mathbb{P},{\mathbb{F}}, \eta, u
\right)& \label{mart-system-s-t}
\end{eqnarray} %\end{equation}
where
\begin{numlistn}
\item  $(\Omega ,{{\mathcal{F}}},{\mathbb{F}},\mathbb{P})$ is a complete filtered
probability space with filtration ${\mathbb{F}}=\{{{\mathcal{F}}}_t\}_{t\ge 0}$,
\item
$ \eta$
is a space time Poisson white noise on $\mathcal{O}$ with jump size
intensity  $\hat \nu=\Leb\otimes \nu$.
\item  $u$ is a $E$-valued ${\mathbb{F}}$--adapted stochastic process such that  \begin{eqnarray} \label{sp-cont-f}
 {\mathbb{E}} \int_0^ T |u(s)|_{E}^ p\, ds
 <\infty,
 \end{eqnarray}
 \item for every $t\in [0,T]$,
    $u$ satisfies the following equation $\mathbb{P}$--a.s.
    \begin{eqnarray}
\lefteqn{ u(t) =
 e^{-tA}u_0 + \int^{t}_{0}e^{-(t-r)A}F(r,u(r))\, dr
}
\\
&&{}  +
\int^{t}_{0} \int_{\mathcal{O}\times {\mathbb{R}}} e^{-(t-r)A}\;   G(u(r))[\xi,\zeta]\,\tilde \eta(d\xi\times d\zeta\times dr)  . %\tilde  \eta(dz,dr).
\nonumber
\end{eqnarray}
\end{numlistn}
\end{definition}

The following result will be shown by applying Theorem
\ref{Th:bound}.
\begin{theorem} \label{stpn} %\label{ex01}
%Fix $r\ge p$.
Let $p\in (1,2]$ and $\nu$ be a L\' evy measure on ${\mathbb{R}}$ satisfying condition
\eqref{eqn-cond_pG}, and $\mathcal{A}$ a differential operator satisfying the properties
\eqref{item-a}-\eqref{itewm-b} above. For $\theta \ge 0$ and $r\ge p$ we put $E=D(A_r^{\frac\theta{2k}})$, where $A_r$
is the linear  operator in the Banach space $L^{r}(\mathcal{O})$ defined
in \eqref{A_high-5}.
Let $F$ and $G$ be two maps defined in \eqref{F-high} and \eqref{DEF-G-high}, respectively.

In addition to Assumptions
\ref{assum-F-high} and  \ref{G-high} assume that the numbers $p, r, \theta, d$ and $k$ satisfy \eqref{neq-theta}.  Then for every $u_0 \in D(A_r^{\frac{\theta}{2k}})$ there
exists a $D(A_r^{\frac{\theta}{2k}})$-valued martingale solution $u$ to \eqref{eqn-HIGH}. Moreover, $u$ is $L^r(\mathcal{O})$-valued c{\`a}dl{\`a}g.
\end{theorem}

The above Theorem can be reformulated in terms of  space time
L\'evy  noise. But, since such a result  would not be
significantly different from the last one, we omit it and leave as an exercise to an interested reader.
\begin{proof}[Proof of Theorem \ref{stpn}]
Let us fix the numbers $d$, $k$, $p$, and $r$, the space $E$ and the operator $A_r$ as in  statement of the
theorem.
Also, let  $F$ (resp. $G$) be defined by equality
\eqref{F-high} (resp.  \eqref{DEF-G-high}).

Since $r\ge p$ then the separable Banach spaces $E$ and $B$ are UMD and M-type $p$.
Thanks to Claim
\ref{c7-1-high},  $A_r$ has the BIP property on $E$, is a positive operator with compact resolvent and its negative $-A_r$ generates a contraction type
$C_0$-semigroup on $E$. Owing to Claim \ref{c7-3-high} we can find $\rho\in [0,\frac1p)$ such that  the map $A_r^{\rho-\frac1p}G $
defined on $[0,T] \times E$ is $L^p(Z,\nu,E)$-valued, bounded and
 continuous w.r.t. $E$ and measurable with respect to time. All assumptions but
Assumption \ref{assum-2} of Theorem \ref{Th:bound} are satisfied
by problem \eqref{eqn-HIGH}. However, it follows from Assumption
\ref{assum-F-high} that the Nemytskii map $F$ defined by
\begin{eqnarray}\nonumber F(t,u)(x) & := & f(t,x,u(x)) , \quad u \in E,
x\in\mathcal{O}, t\in [0,T],\end{eqnarray}
satisfies Assumption \ref{assum-2}. Hence, from the applicability of Theorem
\ref{Th:bound} we easily conclude that problem \eqref{eqn-HIGH} has a $E$-valued martingale solution $u$.
Since $A_r$ is the infinitesimal generator of a contraction $C_0$-semigroup on $L^r(\mathcal{O})$, the paths of the
martingale solution $u$ is \cadlag on $L^r(\mathcal{O})$.
\end{proof}
%%%%%%%%
%%%%%%%%%%%%%%
%%%%%%%%%%%%%
\begin{remark}
Let $n$ be a positive integer and for each $j=1,\cdots,n$ let $\{L_{j,st}(t); t\ge 0\}$ be a L\'evy process with intensity jump size measure $\nu_j$. We assume that $\nu_j$
are L\' evy measures on ${\mathbb{R}}$ such that
\begin{equation}\label{eqn-cond_pG-high}
\text{ there exists } p\in (1, 2] :  \int_{\mathbb{R}} |z|^p \nu_j(dz)
<\infty.
\end{equation}
For any $T\in (0,\infty)$ we consider the following system
  \begin{equation} \label{eqn-SYS-HIGH}
\begin{cases}
&du_i(t)+{\mathcal{A}}_i u_i(t)\, dt=f_i(t,x,u_1(t,x),\dots,u_n(t,x)) dt\\& \quad +
\sum_{j=1}^n g_{ij}(t,x,u_1(t,x),\ldots,u_n(t,x))dL_{j,st}(t),\,\, t\in (0,T], x\in \mathcal{O}, \\
&u_i(0)= u_{i,0}, \,\, x\in \mathcal{O},\\
&\mathcal{B}_{j,i}u_i(t,x)=0, \,\, t\in [0,T], x\in \partial \mathcal{O}.
\end{cases}
\end{equation}
Here $\mathcal{O}$ is a bounded open set of $\mathbb{R}^d$, with $d\ge 1$. For each $i=1,\cdots, n$, ${\mathcal{A}}_i$
are differential operators of order $2k$ satisfying \eqref{item-a}-\eqref{itewm-b} above.
Furthermore, we assume that
\begin{equation*}
 f=[f_i]:[0,T]\times \mathcal{O}\times \mathbb{R}^n\to \mathbb{R}^n, \quad g=[g_{i,j}]\times \mathcal{O}\times \mathbb{R}^n\to \mathbb{R}^n,
\end{equation*}
are separately continuous and bounded. In addition, we assume that $g(t,x,\cdot)$ is uniformly continuous.
By making use of the same approach as above we can show that if $p\in (1,2]$, $r\ge p$, $\theta+d-\frac dr<\frac{2k}p $,\,\, then for any
$u_{i,0}\in H^{1,r}_0(\mathcal{O}, \mathbb{R}^n)\cap  H^{\theta,r}(\mathcal{O}, \mathbb{R}^n)$ there exists a $H^{1,r}_0(\mathcal{O}, \mathbb{R}^n)\cap  H^{\theta,r}(\mathcal{O}, \mathbb{R}^n)$-valued stochastic process $u$ and a space-time L\'evy noise
 $\{L_{j,st}(t); t\ge 0\}, \,\, j=1,\cdots,n,$ defined on some
filtered probability space $(\Omega, \mathcal{F}, \mathbb{F}, \mathbb{F})$ such that
$\mathbb{P}$--a.s.
    \begin{eqnarray*}
\lefteqn{ u_i(t) =
 e^{-tA_i}u_{i,0} + \int^{t}_{0}e^{-(t-s)A_i}f_i(s,x,u_1(s,x),\ldots,u_n(s,x))\, ds
}
\\
&&{}  +\sum_{j=1}^n
\int^{t}_{0} \int_{\mathcal{O}\times {\mathbb{R}}} e^{-(t-s)A_i}\;   g_{ij}(s,x,u_1(s,x),\ldots,u_n(s,x))dL_{j,st}(s). %\tilde  \eta(dz,dr).
\end{eqnarray*}
Here $\{e^{-A_it}: t\ge 0\}, i=1,\cdots,n,$ are the semigroup generated by $-A_i$ which are the linear
unbounded map induced by $-{\mathcal{A}}_i$ on $L^r(\mathcal{O},\mathbb{R}^n)$.

To check this we will apply the above theorem on the Banach space $E=H^{\theta,r}(\mathcal{O},\mathbb{R}^n$). For this purpose we consider the diagonal matrix
\begin{equation*}
 A=\begin{pmatrix}
A_1 &  0  & \ldots & 0\\
0  &  A_{2} & \ldots & 0\\
\vdots & \vdots & \ddots & \vdots\\
0  &   0       &\ldots & A_{n}
\end{pmatrix},
\end{equation*}
and denote by $G_0:[0,T]\times L^r(\mathcal{O},\mathbb{R}^n) \to L^r(\mathcal{O},\mathbb{R}^n)$ (resp. $F$) the Nemytskii operator associated to the matrix $g$ (resp. the vector $f$).
We also set $Z=\mathbb{R}^n$ and define the L\'evy measure $\hat{\nu}$ on $\mathcal{O}\times Z$ by
$\hat{\nu}=\Leb \otimes(\nu_1\otimes \cdots\otimes \nu_n)$.  As above we can define a bounded linear map
$\Phi: E\to L^p(\mathcal{O} \times Z; B^{-(d-d/r)}_{r,\infty} (\mathcal{O}, \mathbb{R}^n) )$ by \eqref{eqn-Phi}, i.e.
  \begin{equation*}
[\Phi v] (x, y) = (v(x)\delta_ x) y, \;\; (x,y)\in \mathcal{O}\times {\mathbb{R}}^n,
\end{equation*}
and  $G=\Phi \circ G_0$. The restriction of $F$ and $G$ to
$[0,T]\times E$ are still denoted by $F$ and $G$, respectively. We
denote by $\eta$ the Poisson random measure with intensity measure
$dt\times \hat{\nu}(dx,dz) $ on $[0,T]\times \mathcal{O}\times Z$. Then we
can rewrite problem \eqref{eqn-SYS-HIGH} in the following form
\begin{equation}
\begin{cases}
 du&+Au dt =F(t,u)dt+ \int_{\mathcal{O}\times Z} G(t,u)[x,z]\tilde{\eta}(dx\times dz\times dt),\\
 u(0)&=u_0.
 \end{cases}
 \end{equation}
 The claim is now a straightforward corollary of the the above theorem.
\end{remark}
%%%%%%%%%%
%%%%%%%%%
%%%%%%%%%%%
%
%%%%%%%%%%%%%%%%%%%%%%%%%%%%%
%%%%%%%%%%%%%%%%%%%%%%%%%%%%%%%%%%%%%%%
\section{Application III: Stochastic evolution equations with fractional generator and polynomial nonlinearities} \label{det-set-B}%A elliptic operator}
%%%%%%%%%%%%%%%%%%%%%%%%%%%%%%%%%%%%%%%%
%%%%%%%%%%%%%%%%%%%%%%%%%%%%%%%%%%%%%%%%%
%%%%%%%%%%%%%%%%%%%%%%%%%%%%%%%%%%%%%%%%%
%%%%%%%%%%%%%%%%%%%%%%%%%%%%%%%%%%%%%%%%%
 We will deal with a similar problem as in the previous section, but we will assume that the nonlinear term $F$ is of polynomial type.
Let $k=1$, $\mathcal{A}$
and $\mathcal{B}$ be a differential operators satisfying the items \eqref{item-a}-\eqref{itewm-b} described in
page \pageref{det-set-A}. Let $\gamma\in (0,1]$ and $A_r$ be the map induced by $(-\mathcal{A})^\gamma $ on the Banach space $L^r(\mathcal{O})$, $r>1$.
We also consider a space-time levy noise $L_{st}=\{L_{st}(t); t \ge 0\}$ with intensity measure $\nu$ satisfying \eqref{eqn-cond_pG} for some $p\in (1,2]$.
Next, let $g: [0.T]\times \mathbb{R}\times \mathcal{O} \to \mathbb{R}$ be a separately continuous wrt to the first and second variables, uniformly continuous wrt the third variable.
Furthermore we assume that
$g$ is bounded. For any $r\ge p$ we define on $[0,T]\times L^r(\mathcal{O})$ a map $G$ by
\begin{equation}\label{G-poly}
 [G(t,u)](x,y)= [(g(t,u(x), x) \delta_x]y, \,\, u\in L^r(\mathcal{O}),\,\ (x,y)\in \mathcal{O}\times \mathbb{R}.
 \end{equation}
We also consider a separately continuous
real valued functions $f$ defined on $[0,T]\times \mathcal{O}\times
\mathbb{R}$ such that there exists a constant $K>0$ and
\begin{equation}
-K(1+\vert u\vert^q1_{\{u\ge 0\}})\le f(t,x ,u)\le K(1+\vert
u\vert^q1_{\{u\le 0\}}),
\label{F-diss-Poly}
\end{equation}
for all $t\ge0, x \in {\mathcal{O}}, u \in {{\mathbb{R}}}$. It is not difficult to prove that  if $f$ satisfies
\ref{F-diss-Poly} then  $$f(t,v+z)\sgn v \le K(1+|z|^q),$$ for all
$v,z \in \mathbb{R}$ and $t\in [0,T]$. We denote by $F$ the Nemitskii operator associated to $f$ which is defined by
\begin{eqnarray}\label{Nem-Poly}
F(t,u)(x) & := &  f(t,x,u(x)) , \quad u
\in C_0(\mathcal{O}), x\in\mathcal{O}, t\in [0,T].
\end{eqnarray}
We approximate %the function
$f$ by a sequence $(f_n)_{n\in{\mathbb{N}}}$ of functions defined by
$$
f_n(t,x,u) := \bcase f(t,x  ,u), & \mbox{ if } u\in[-n,n]\\  f(t,x ,n), &
\mbox{ if } u\ge n,\\  f(t,x ,-n), & \mbox{ if } u\le - n,\ecase
$$%
for any $t\in [0,T]$, $x\in \mathcal{O}$ and $n\in \mathbb{N}$. By setting $
F_n(t,u)(\xi)= f_n(t,\xi,u(\xi) )$ for $(t,u, \xi )\in [0,T]\times C_0(\mathcal{O})\times \mathcal{O}$ we obtain a  sequence $(F_n)_{n\in{\mathbb{N}}}$ of
functions defined on $[0,T]\times C_0(\mathcal{O})$ into $C_0(\mathcal{O})$ which are bounded, separately
 continuous maps  satisfying
 \eqref{Dissip-F} uniformly in $n$, and
 pointwise converging to $F$ in $C_0(\mathcal{O})$.  Hence we
have the following result.
\begin{claim}\label{c7-Poly}
The nonlinear map $F$ defined by \eqref{Nem-Poly} satisfies
Assumption \ref{assum-main-3} with $X=C_0(\mathcal{O})$.
\end{claim}
\begin{remark}
An example of a real valued function $f$ satisfying the above
conditions is
 \begin{eqnarray}\label{EX-POLY} f:[0,\infty)\times \mathcal{O}
\times {\mathbb{R}} \ni (t,\xi,u) \mapsto -|u|^ q \sgn(u).
\end{eqnarray}
\end{remark}
With the various mappings we have introduced above we consider the following SPDEs
 \begin{equation} \label{eqn-Poly}
  \begin{cases}
&du(t)=-(-{\mathcal{A}} )^\gamma u(t)\, dt+f(t,x,u(t,x)) dt\\& \quad \quad \quad +
g(t,x,u(t,x))dL_{st}(t),\,\, t\in (0,T], x\in \mathcal{O}, \\
&u(0)= u_{0}, \,\, x\in \mathcal{O},\\
&\mathcal{B}u(t,x)=0, \,\, t\in [0,T], x\in \partial \mathcal{O}.
\end{cases}
\end{equation}
\begin{theorem} \label{Poly} %\label{ex01}
%Fix $r\ge p$.
Let $\gamma\in (0,1]$, $p\in (1,2]$, $\nu$ be a L\' evy measure on ${\mathbb{R}}$ satisfying condition
\eqref{eqn-cond_pG}, and $\mathcal{A}$ a differential operator satisfying the properties
\eqref{item-a}-\eqref{itewm-b} from \pageref{det-set-A}.
Let $F$ and $G$ be two maps defined in \eqref{c7-Poly} and \eqref{G-poly}, respectively.
In addition to the assumptions on $G$ above we also assume that $pd<2\gamma$ and F satisfies Assumption \ref{F-diss-Poly} with
$q\in (1,p)$. If $B=L^r(\mathcal{O})$ with r satisfying $r>\max\{p, \frac{p d}{p-q}\}$, then for any $u_0 \in C_0(\mathcal{O})$ there
exists a $C_0(\mathcal{O})$-valued, $L^r(\mathcal{O})$-valued \cadlag martingale solution to \eqref{eqn-Poly}.
\end{theorem}
Before we prove this result let us make the following remark.
\begin{remark}
Let us assume that the space time white noise has a jump size intensity measure $\nu$ which is finite\footnote{Note that  the times between the jumps are exponential distributed.}.
Then, as in the proof of Theorem IV.9.1 in \cite{Ikeda+Watanabe_1981}  the solution can be written as a concatenation
of  solutions to the deterministic reaction diffusion equation over random intervals with the initial data being a measure valued random variable.

To be more precise, let $\lambda:= \Leb(\mathcal{O}) \times \nu({\mathbb{R}})$, $\{\tau_i:i\in {\mathbb{N}}\}$ be a family of independent, exponential distributed with parameter $\lambda$ random variables and
$$N(t)=\sum_{n=1} ^\infty 1_{[T_n,\infty)}(t), \quad t\ge 0,
$$
 where  $T_n=\sum_{i=1} ^n\tau_i$, $n\in{\mathbb{N}}$. Let also $\{Y_i:i\in{\mathbb{N}}\}$ be a family of independent $\nu/\nu(\mathbb{R})$ % \coma{I guess this should be $\nu/\nu(\mathbb{R})$ }
 distributed random variables
and  $\{ x_i:i\in{\mathbb{N}}\}$ be a sequence  of independent and uniformly distributed random variables in $\mathcal{O}$.
Then, the space time white noise $\eta $ can be written as follows: for any $A\in{\mathcal{B}}(\mathcal{O})$, $B\in {\mathcal{B}}({\mathbb{R}})$ and $ I\in{\mathcal{B}}([0,\infty))$
$$
\eta(A\times B\times I)  = \bcase 0 & \mbox{ if } N(t)=0,\\\sum_{i=1}^{N(t)} \delta_{x_i,Y_i,T_i}(A\times B\times I)& \mbox{ if } N(t)>0 .
\ecase
$$
Using this representation the above SPDEs can be described by a deterministic PDE with inital condition
being a measure
in the time intervals $[T_n,T_{n+1})$, i.e. $u$ solves the deterministic PDE
 \begin{eqnarray}\left\{
\begin{array}{rcl} \frac{\partial }{
\partial t} u(t,\xi)&+&{\mathcal{A}} u(t)\, dt=f( u
(t,\xi))\phantom{\Big|} \,\,
\xi\in\mathcal{O},\, t\in (T_n,T_{n+1}),
\\
u(T_n ^+,\xi) &=& u(T_n ^-)+ Y_n\delta_{x_n},\quad\xi\in\mathcal{O},\phantom{\Big|}
\\
 u(t,\xi) &=&
0, \hbox{ for }\xi \in \partial \mathcal{O}, \,  t\in (T_n,T_{n+1}).\end{array}\right.
\end{eqnarray}
It follows that our conditions have to be stronger than the conditions in \cite{brezis}, which is indeed the case.
In fact, for  $\gamma=1$
 we assume that $d<\frac{2}{q}$ which is stronger than $d\le\frac2{q-1}$ imposed by Brezis and Friedman in \cite{brezis}.
\end{remark}
\begin{proof}[Proof of Theorem \ref{Poly}]
We just give a sketch of the proof because it is very similar to the proofs of Theorem \ref{ex01} and Theorem \ref{stpn}.
Let us fix the numbers $d$, $\gamma$, $p$, $q$, and $r$ as in the statement of the theorem. We denote by $A_r$ the operator induced by $(-\mathcal{A})^\gamma$ on
the Banach space $B$. Also, let $F$ be defined by equality
\eqref{Nem-Poly}.

 By Remark \ref{REM-c7-3-high} we can find $\theta>\frac{d}r$ such that $ \theta +d-\frac{d}{r}< \frac{2\gamma}p$ and the Banach space
 $E=D(A_r^{\frac{\theta}{2\gamma}})$ is continuously embedded in $X:=C_0(\mathcal{O})$. Thus, owing to the assumption on $g$ (resp. $f$) we can argue as in the proof of Claim
 \ref{c7-3-high} (resp. Claim \ref{Nem-Poly}) to prove the map $G$ (resp. $F$) satisfies Assumption \ref{assum-main-2}
(resp. Assumption \ref{assum-main-3})  with $\rho\in (0, \frac1p - \frac d {2\gamma})$. Note that since $r\ge p$, $E$ and $B$ are separable, UMD and type $p$ Banach spaces.
Finally, let $Y=L^r(\mathcal{O})$ and $A_Y=A_r$. Since $1-\frac{q}{q_{\max}}>1-\frac qp$ and $r> \frac{ pd}{p-q}$, we can find $\kappa_1\in (\frac dr, 1-\frac qp)$
such that
$$ D(A_r^{\frac{\kappa_1}{2\gamma}})\subset X \subset Y .$$ Note also that $-A_r$ is the infinitesimal generator of a contraction type
$C_0$-semigroup on $Y=B=L^r(\mathcal{O})$. Thus, all the assumptions of Theorem \ref{Th:general} are satisfied
by problem \eqref{eqn-HIGH}. Hence, we easily conclude the
proof of Theorem \ref{Poly} from the applicability of Theorem
\ref{Th:general}.
\end{proof}
%%%%%%%%%%%%%%%%%%%%%%%%%%%%%%%%%%%%%%%%%%%%%%%%%%%%%
%%%%%%%%%%%%%%%%%%%%%%%%%%%%%%%%%%%%%%%%%%%%%
\section{Some Preliminary Results}\label{preliminaries}

The purpose of this section is twofold. First we will summarize
some results concerning the deterministic convolution process.
Here we will use results already shown in \cite{634.60053} and
\cite{MR1488138}. Secondly, we will state several results
concerning the stochastic convolution process. We begin with
introducing some notation.

\begin{remark}
Without loss of generality we can assume that $-A$ generates a
contraction semigroup on a Banach space $E$. \delb{This is no
restriction, since}Indeed,  if $A$ satisfies Assumption
\ref{assum-1}-(\ref{h2-b}), then there exists   $\mu>0$
such that $-A-\mu I$ generates a contraction semigroup. \delb{In this
case $\lambda$ appearing in the estimate \eqref{int-est-theo} has
to be chosen at least larger than $\nu$.}
\end{remark}

\begin{definition}\label{def-notation-1}
For any separable Banach space $Y$ and a real number $q\in[1,\infty)$, we denote
by $L ^ q (0,T; Y)$, see e.g.
\cite{Diestel-Uhl_1977},
 the Lebesgue   space consisting of (equivalence classes of)  measurable functions $u:[0,T]\to Y$ such that
\delb{$\int_0 ^T  \: |u(s)| ^ q _Y \: ds <\infty$.
 Equipped with the norm}
$$
\left\| u\right\|_{L ^ q(0,T; Y)} := \left(\int_0 ^T  \: |u(s)| ^ q _Y
 \: ds\right) ^ \frac 1q <\infty.  \delb{,\quad u\in L ^ q(0,T; Y).}
$$
The Besov-Slobodetski space $W^{\alpha,q}(0,T;Y)$, where $\alpha \in (0,1)$,  consists of
 all $u\in L ^ q(0,T;Y)$ such that
 \begin{equation}\label{eqn-Besov-seminorm}
\Lve\; u\Rve_{W^{\alpha,q}(0,T;Y)}^q:=\int_0 ^T \int_0 ^T \frac{|u(t)-u(s)|_Y ^ q}{|t-s| ^ {1+\alpha q}}\,ds\,dt<\infty.
\end{equation}
Both the spaces $L ^ q (0,T; Y)$ and  $W^ {\alpha,q}(0,T;Y)$ are equipped with the natural norms
\delb{$$
\left\| u\right\|_{ W^ {\alpha,q}(0,T;Y)}:=\left( \int_0 ^T \int_ 0
^{T}\,{|u(t)-u(s)|_Y
 ^ q\over |t-s| ^ {1+\alpha q}}\,ds\,dt\right) ^ \frac 1q,
$$}
are  separable Banach spaces.

We will denote by ${\mathcal{C}}([0,T]; Y)$  the space of all $Y$-valued
continuous functions.  We will denote by
${{\mathcal{C}}}^\beta([0,T];Y)$, where $\beta\in (0,1]$,  the space  of  all functions
 $u \in{\mathcal{C}}([0,T];Y)$ such that
\begin{eqnarray}\label{notation-3}
\Vert u \Vert_{{{\mathcal{C}}}^\beta([0,T];Y)} &:=& \sup_{0\le t\le T }
|u(t)|_Y + \sup_{0 \le s < t \le T} \; \frac{|u(t)-u(s)|_Y }{
|t-s|^\beta} <\infty. \label{3.holder}
\end{eqnarray}
The spaces ${\mathcal{C}}([0,T]; Y)$   and ${{\mathcal{C}}}^\beta([0,T];Y)$ are Banach spaces. The latter one is not separable.

By $H^ {1,q}(0,T,Y)$ we will denote  the Banach spaces of
(classes of) functions $u\in L^ {q}(0,T,Y)$ whose weak
derivative $u^\prime$ belongs to $L ^ q(0,T;Y)$ as well. \delb{By $H_{0} ^
{1,q}(0,T;Y)$} We will also use the following notation  \delb{the closure in $H ^ {1,q}(0,T,Y)$
of the space $\{ u\in{\mathcal{C}} ^ \infty(0,T;Y): u(0)=0\}$. It is known
that} $$ H_{0} ^ {1,q}(0,T;Y)=\{u\in H ^ {1,q}(0,T,Y):
u(0)=0\}.$$
It is well known
that $ H_{0} ^ {1,q}(0,T;Y)$ is a closed subspace of $ H^ {1,q}(0,T;Y)$.
Moreover, if $X$ is another separable Banach space which is compactly embedded in $Y$, then, provided that $\alpha-\frac1p >-\frac 1q$ and $\beta>0$, the embeddings $$
W^ {\alpha,p}(0,T;X) \embed L^q(0,T;Y),$$
$${{\mathcal{C}}}^\beta([0,T];X) \embed {\mathcal{C}}([0,T]; Y),$$  and $$ H^ {1,q}(0,T;X) \embed  L^ {q}(0,T;Y)$$  are  compact.

\end{definition}

\begin{remark}\label{skorohodd} A known fact is that the Skorokhod topology is
weaker than the uniform topology (see e.g.\ \cite[Lemma 6.8, p.\
248]{para}), i.e.\ the embedding $ C([0,T];Y)\hookrightarrow
\mathbb{D}([0,T];Y) $ is continuous.
\end{remark}

In the part dealing with the stochastic convolution process we
will need also the following notation.
\begin{definition}\label{def-notation-2}
Assume that $Y$ is separable Banach space and   $q\in[1,\infty)$.

We will denote by
 ${\mathbf{L}}^q(\Omega_T;Y)$  the equivalence classes of all $\mathcal{B}([0,T])\times\mathcal{F}/\mathcal{B}(Y)$-measurable
functions $\xi :[0,T]\times\Omega\to Y$ such that \begin{eqnarray*} &&
\mathbb{E}\int_0^T \vert\xi (t)\vert_Y^q\,dt<\infty.\
\label{def:Mq-B} \end{eqnarray*}
 By
${\mathbf{L}}^p\left(\Omega; Y\right)$ we will denote the Banach space of equivalence classes of all $\mathcal{F}/\mathcal{B}(Y)$-measurable functions $\xi: (\Omega,{\mathcal{F}})
\to Y$ such that ${\mathbb{E}} |\xi|^p_{Y }<\infty$, see
\cite{Diestel-Uhl_1977}.\\
We will denote by $\mathcal{N}(0,T;Y)$ we denote the space of
(equivalence classes) of $\mathbb{F}$-progressively-measurable processes $\xi
:[0,T]\times\Omega\to Y$.

By ${\mathcal{N}}^q (0,T;Y)$, resp. $\mathcal{M}^q(0,T;Y)$,  we will denote the
space, resp. Banach space,   consisting of all processes
$\xi\in \mathcal{N}(0,T;Y)$ satisfying, respectively,
 \begin{eqnarray*} && \int_0^T \vert\xi
(t)\vert_Y^q\,dt<\infty, \mbox{  $\mathbb{P}$-a.s.} \label{def:Nq}\\
 && \mathbb{E}\int_0^T
\vert\xi (t)\vert_Y^q\,dt<\infty.\ \label{def:Mq} \end{eqnarray*}

If $X$ is a metric space,  %
by ${\mathbf{L}} ^ 0\left(\Omega; X\right)$ we denote the set of measurable
functions from $(\Omega,{\mathcal{F}})$ to $X$.
\end{definition}
\subsection{The operator $\Lambda$}
Let $E$ be a separable Banach space and let us fix two real numbers $q\in (1,\infty)$ and $\mu\geq 0$.

Similar to Brze{\'z}niak and G{\c{a}}tarek \cite{Brz+Gat_1999} we
define  a linear operator ${{\mathcal{A}}}$ by the formula
\begin{eqnarray}\label{operator_def}
D({{\mathcal{A}}})&=&\left\{ u\in L ^q(0,T ;\BB): \; Au \in L ^q(0,T
;E)\right\},
\nonumber\\
 {{\mathcal{A}}}u&:=&\left\{(0,T) \ni t \mapsto A(u(t)) \in \BB \right\},\quad u\in D({\mathcal{A}}).
\label{2.8}
\end{eqnarray}
It  is  known, see  \cite{MR910825}, that  if $A+\mu I$ satisfies
Assumption \ref{assum-1}-(\ref{h2-a}, \ref{h2-b}, \ref{h3}) then
${{\mathcal{A}}} + \mu I$ satisfies them as well.
With $q$ and $\BB$ as above we define two operators ${\mathcal{B}}$ and $\Lambda$, see \cite{Brz+Gat_1999}, by
 \begin{eqnarray} {\mathcal{B}} u&=&u^\prime,\quad u\in
D({\mathcal{B}}):= H_{0} ^ {1,q}(0,T;\BB).
\\
\Lambda&:=& {\mathcal{B}} +{{\mathcal{A}}}, \quad\quad D(\Lambda):=D({\mathcal{B}})\cap
D({{\mathcal{A}}}). \label{2.10}
\end{eqnarray}

If $E$ is a UMD Banach space and $A+\mu I$, for some $\mu \ge 0$,
satisfies Assumption \ref{assum-1}-(\ref{h2-a}, \ref{h2-b},
\ref{h3}) then, by
\cite{MR910825} and \cite{Giga+Sohr_1991}, since $\Lambda={\mathcal{B}}-\mu I +{{\mathcal{A}}} +\mu I$,  $\Lambda $ is a positive
operator. In particular, $\Lambda$ has a bounded inverse. The
domain $D(\Lambda)$ of $\Lambda$ endowed with the `graph' norm
\begin{equation}
\Vert u\Vert = \left\{ \int_0^T |u^\prime(s)|^q \, ds + \int_0^T
|Au(s)|^q \, ds \right\}^{\frac1q} \label{2.11a}
\end{equation}
is a Banach space.

Before continuing, we present two results on the fractional powers
of the  operator $\Lambda$, see  \cite{MR1488138} for the proof.

\begin{proposition}\label{Prop:2.0} Assume that Assumption \ref{assum-1}-(\ref{h1}, \ref{h2-a}, \ref{h2-b}) are satisfied.
Assume also that  for some
$\mu \ge 0$,   $A+\mu I$ satisfies Assumption \ref{assum-1}-(\ref{h3}). \\
Then,  for any $ \alpha \in(0, 1]$,  the operator
$\Lambda^{-\alpha}$ is a bounded linear operator in $L^q
(0,T;E)$, and for $f \in L^q(0,T;E)$,
\begin{eqnarray}
\left(\Lambda^{-\alpha}f\right) (t)&=&\frac1{\Gamma (\alpha)}  \int_0^t (t-s)^{\alpha -1}
e^{-(t-s)A}f(s) \, ds, \;\; t \in [0,T].
\label{2.13a}
\end{eqnarray}
\end{proposition}

\begin{lemma}\label{L:reg}
Let Assumption \ref{assum-1}-(\ref{h2-a}, \ref{h2-b}) holds.
Suppose that the positive numbers $\alpha, \beta,\delta$ and $q>1$
satisfy
\begin{equation}
\alpha  -\frac1q   +  \gamma>\beta +\delta.
\label{cond:1}
\end{equation}
    If $T  \in (0,\infty)$, then  the operator
\begin{equation}
\Lambda ^ {-\alpha} : L^q (0,T;D(A^\gamma)) \to
{{\mathcal{C}}}^\beta ([0,T];D(A^\delta)) \; \label{2.13b}
\end{equation}
is bounded.
\delb{
If $T=\infty$  and   the semigroup $\{e^{-tA}\}_{t\ge 0} $ is
exponentially bounded on $E$, i.e., for some  $ a>0$, $C>0$
\begin{eqnarray}
 |e^{-tA}|_{L(E,E)} &\le& Ce^{-at}, \; t \ge 0,
\label{exp-bound}
\end{eqnarray}
then  for any $f \in L^q_0(0,T;D(A^\gamma))$ the function
$u =\Lambda ^ {-\alpha} f$ belongs to
${{\mathcal{C}}}^\beta([0,T];D(A^\delta))$.
 Moreover, the operator $\Lambda ^ {-\alpha}$ is a bounded map in the
above spaces.
\begin{equation}
\Lambda ^ {-\alpha} : L^q(0,T;D(A^\gamma)) \to
{{\mathcal{C}}}^\beta([0,T];D(A^\delta)), \label{2.13c}
\end{equation}
is bounded.}
\end{lemma}

Finally, we present   slight modifications of the Proposition 2.2
and Theorem 2.6 from  \cite{MR1488138}.
\begin{theorem}
\label{Th:compact} Assume that $\BB$ is an UMD Banach space and
an operator $A$ satisfies   Assumptions \ref{assum-1}-(\ref{h2-a},
\ref{h2-b},\ref{h3}) \delb{is such that $A+\mu I$, for a $\mu \ge 0$, satisfies
condition Assumption and \ref{assum-1}-(\ref{h3})}. Assume that
$ \alpha \in (0,1]$ and $\delta,\gamma\geq  0$ are such that
\begin{equation}
\alpha-\frac1q +\gamma-\delta>-\frac1r.
\label{cond:2}
\end{equation}
 Then the   operator
 \begin{equation}
 A ^{\delta} \Lambda^{-\alpha}A ^ {-\gamma} : L^q (0,T;\BB)
\to  L^r(0,T;\BB)
\label{cond:3}
\end{equation}
is bounded. Moreover, if the operator $A^{-1}:\BB\to \BB$ is   compact, then the
operator in \eqref{cond:3}   is compact.
\end{theorem}
\begin{remark}\label{Prop:2.1-remark-F}
In view of Theorem \ref{Th:compact}
$$\Lambda ^ {-1}:L ^ p(0,T;E)\to L ^ p(0,T;E)
$$
is a well defined bounded linear operator for $p\in [1,\infty)$.
\end{remark}
\begin{corollary}\label{C:comp} Assume the first set of assumptions of Theorem
\ref{Th:compact} are satisfied. Assume  that  three  non-negative
numbers $\alpha,  \beta, \delta$ satisfy  the  following condition
\begin{equation}
\alpha - \frac1q>\beta + \delta.
\label{cond:1a}
\end{equation}
 Then the operator
$\Lambda^{-\alpha}: L^q(0,T;\BB ) \to
{{\mathcal{C}}}^\beta([0,T];D(A^\delta))$ is bounded. Moreover, if the
operator $A^{-1}:\BB\to \BB$ is  compact, then  the
operator $\Lambda^{-\alpha}: L^q(0,T;\BB ) \to {{\mathcal{C}}}^\beta([0,T];D(A^\delta))$  is also compact.\\
{ In particular, if $\alpha > \frac1 q$ and the operator $A^{-1}:\BB\to \BB$ is compact, the map
$\Lambda^{-\alpha}:L^q(0,T;\BB ) \to {{\mathcal{C}}}([0,T];\BB )$ is compact.
\delb{- Compare with the beginning of Section 8.}}
\end{corollary}
\subsection{The stochastic convolution}
Let $E$ be a separable, UMD and type $p$ Banach space. Let $X$ be a (arbitrary) separable Banach space and $G$ a nonlinear map defined on
$[0,T]\times X$ satisfying Assumption \ref{assum-main-2} and for
any $u\in \mathcal{N}(0,T;X)$ we will denote $G(s,u(s),z)$ by
$\xi(s,z)$.

We assume that $\rho\in [0,\frac 1p)$ and we fix  $M>0$. Throughout this section we
denote by $\mathcal{B}_M(X)$ the set of all $\mathbb{F}$-progressively measurable processes $\xi$
satisfying
\begin{equation}\label{Bound-xi}
\int_Z\lvert A^{\rho-\frac1p} \xi(s,z)\rvert^p_E\nu(dz)\leq M^p, \mbox{ for a.a. } s\in
[0,T].
\end{equation}
 We also set
\begin{equation}\label{eqn-stoch_ter}
\mathfrak{G}(\xi)=\Big\{ [0,T]\ni t\mapsto \int_0 ^ t\int_Z e ^ {-(t-s)A}
\xi(s,z)\;\tilde \eta(dz,ds), \xi \in \mathcal{B}_M(X) \Big\}.
\end{equation}

We will prove some results concerning the family $\{\mathfrak{G}(\xi): \xi \in
\mathcal{B}_M(X)\}$.
\begin{lemma}\label{est-stoc-conv-1} If  $\rho^\prime \in (0,\rho)$, then
there exists a constant $C_1>0$ such that
\begin{equation*}
{\mathbb{E}} \Lve A^{\rho^\prime} \mathfrak{G}(\xi)\Rve^p_{L^p(0,T:E)}\le C_1 T M^p,\;\; \xi \in
\mathcal{B}_M(X).
\end{equation*}
\end{lemma}
\begin{proof} Set $u=\mathfrak{G}(\xi)$. Let us fix $\rho^\prime \in (0,\rho)$. Since  $\rho^\prime+\frac1p-\rho<\frac1p$ and
\[
u(t)=A^{\frac1p-\rho} \int_0 ^ t\int_Z e ^ {-(t-s)A} \big[A^{\rho-\frac1p} \xi(s,z)\big]\;\tilde \eta(dz,ds),\;\; t\in [0,T],
\]
by \cite[Theorem 2.1]{maxreg}, see also \cite[Lemma 4.6]{uniqueconv}, we infer that
\begin{equation}\label{ineq-stoc-conv-1b}
{\mathbb{E}} \Lve A^{\rho^\prime}  u \Rve^p_{L^p(0,T:E)}
 \leq C_1\int_0^T\int_Z
\lve A^{\rho-\frac1p}\xi(s,z)\rve^p \nu(dz)dt \leq C_1T M^p.
\end{equation}
This completes the proof of the lemma.
\end{proof}
 Let us note that since, by
assumption, $A^{-1}$  exists and is bounded  the fractional powers
$A^{-\gamma}$, $\gamma \ge 0,$ are bounded too.
\begin{lemma}\label{lem-est-stoc-conv-2}
Assume that $\rho^\prime \in (0,\rho)$ and put $B=D(A^{\rho^\prime-1})$.
If $\xi \in \mathcal{B}_M(X)$, then
\begin{enumerate}[(i)]
 \item
\label{lem-est-stoc-conv-2-i} \delb{For any $1\le \bar{p}\le p$}  there exists
a constant $C_2=C_2(T)>0$ such that
\begin{equation*}
{\mathbb{E}} \sup_{0\le t\le T}\lve
A^{\rho^\prime-1}\mathfrak{G}(\xi)(t)\rve^{p}\le C_2 M^p,
\end{equation*}
\item\label{lem-est-stoc-conv-2-ii}  and the  process $u=\mathfrak{G}(\xi)$ admits a
$B$-valued \cadlag\, modification (which will be still denoted by
$\mathfrak{G}(\xi)$).
\item\label{lem-est-stoc-conv-2-dissipative}
If in addition the operator $A$ generates a contraction type semigroup on the space $D(A^{\rho-\frac1p})$, then the parts (i)-(ii) are true for $B=D(A^{\rho-\frac1p})$.
\end{enumerate}
\end{lemma}

\begin{proof} Part \eqref{lem-est-stoc-conv-2-dissipative} is a consequence of \cite{Brz+Haus+Zhu_2010}. \\
Let us fix $\xi \in
\mathcal{B}_M(X)$ and set $u=\mathfrak{G}(\xi)$ .
 In order to prove part \eqref{lem-est-stoc-conv-2-i} we recall that according to  \cite[Lemma
4.7]{uniqueconv} we have, with $\beta:=1-\frac1p+ \rho -\rho^\prime>0$, for $t\in [0,T]$,
\begin{equation}\label{est-conv-aux}
A^{\rho^\prime-1}u(t)=\int_0^t A^{\rho^\prime}u(s)
ds+A^{-\beta}\int_0^t\int_Z
A^{\rho-\frac1p}\xi(s,z)\tilde{\eta}(dz,ds).
\end{equation}
 Using \cite[Corollary C.2]{maxreg},
 we obtain
\begin{eqnarray*}
{\mathbb{E}}\sup_{0\le t\le T}\biggl \lve \int_0^t \int_Z
A^{\rho-1}\xi(s,z)\tilde{\eta}(dz,ds)\biggr\rve^{p}
&\leq& C {\mathbb{E}}\int_0^T \int_Z \lve A^{\rho-\frac1p}
\xi(s,z)\rve^{\bar{p}} \nu(dz)ds\\ &\leq&  C T M^p.
\end{eqnarray*}
Next applying H\"older's inequality twice and invoking inequality \eqref{ineq-stoc-conv-1b} we get
\begin{equation*}
{\mathbb{E}}\sup_{0\le t\le T}\biggl\lve \int_0^t A^{\rho^\prime}
u(s) ds\biggr\rve^p\leq
 T^{p-1}{\mathbb{E}}\int_0^T \lve
A^{\rho^\prime} u(t) \rve^p \,dt
\le  C_1 T^p M^p.
\end{equation*}
This completes  the proof of the 1st part  with $C_2(T)=CC_{21}T+C_1T^p$.\\
 Part \eqref{lem-est-stoc-conv-2-ii} is a
consequence of \cite[Lemma 4.7]{uniqueconv}. It can also be deduced from identity \eqref{est-conv-aux}.
\end{proof}

The next lemma is about  estimates of $\mathfrak{G}(\xi), \,\xi\in
\mathcal{B}_M(X)$ in the Besov-Slobodetski spaces $W^{\alpha, p}(0,T; E)$, see Definition \ref{def-notation-1}.
\begin{lemma}\label{est-stoc-conv-3}
Assume that  $\alpha\in (0,\rho)$. Then there exists a number $C_3>0$
such that
\begin{equation*}
{\mathbb{E}} \Lve
\mathfrak{G}(\xi)\Rve^p_{W^{\alpha,p}(0,T;E)}\le C_3 M^p,\;\;\ \xi \in
\mathcal{B}_M(X).
\end{equation*}
\end{lemma}
\begin{proof} Let us fix  $\xi \in
\mathcal{B}_M(X)$ and put $u=\mathfrak{G}(\xi)$. Let us fix $\alpha \in
(0,\rho)$ and let us choose an auxiliary $\rho^\prime\in
(\alpha,\rho)$. We will estimate the seminorm
$\Lve\; \cdot \Rve_{W^{\alpha,p}(0,T;Y)}$ of $u$. For this aim, without loss of generality,  we can take $s<t\in [0,T]$. As in the Gaussian case we have
\[u(t) - u(s) = \tilde{S}_1+\tilde{S}_2,\]
where
\begin{align*}
\tilde{S}_1=&\int_s ^{t}\int_Z e ^ {-(t-r)A}
   \xi(r;z)\; \tilde \eta(dz;dr)
  & \tilde{S}_2= & \left( e ^ {-(t-s)A}  -I\right) u(s).
\end{align*}
 In view of the definition \eqref{eqn-Besov-seminorm}  \delb{ $$ {\mathbb{E}}\Lve
u\Rve^p_{W^{\alpha,p}(0,T:E)} =2\int_0^T\int_0^t
\frac{\lve u(t)-u(s)\rve^p \,\,ds\,\,dt}{\lve
t-s\rve^{1+p\alpha}},$$} it is sufficient to prove that there exist two positive numbers $C_{31}, C_{32}$ such that
\begin{align*} S_1:=&{\mathbb{E}} \int_0^T \int_0^t
\frac{dsdt}{\lve t-s\rve^{\alpha p +1}} \biggl\lve \int_s^
t\int_Z e^{-(t-s)A} \xi(r,z)\tilde{\eta}(dz,dr)\biggr\lve^p\le C_{31} M^p,\\
S_2:=&{\mathbb{E}} \int_0^T \int_0^t \frac{\lve \left( e ^ {-(t-s)A}  -I\right)
u(s)\rve^p \,\,ds \,\, dt}{\lve t-s\rve^{1+\alpha p}}\le C_{32} M^p.
\end{align*}
 Let us start with $S_1$. By
using the Fubini Theorem,  \cite[Corollary C.2]{maxreg}, the estimate $\Lve A^{\frac1p-\rho}e^{-(t-r)A}\Rve^p_{\mathcal{L}(E)}\le
C (t-r)^{-p(\frac1p-\rho)}$ and the definition  \eqref{Bound-xi}  of the class $\mathcal{B}_M(X)$, we infer that

\begin{eqnarray*}
 S_1 &\leq& C \int_0^T \int_0^t \frac{dsdt}{\lve t-s\rve^{1+\alpha
p}} \times \\ && \quad  \quad {\mathbb{E}}  \int_s^t \int_Z \Lve
A^{\frac1p-\rho}e^{-(t-r)A}\Rve^p_{\mathcal{L}(E)}\lve
A^{\rho-\frac1p}\xi(r,z)\rve^p\nu(dz)dr
\\
&\leq & C M^p \int_0^T \int_0^t \frac{dsdt}{\lve t-s\rve^{1+\alpha
p}}  \int_0^{t-s} r^{p\rho-1} dr\\
&\le& CM^p \int_0^T \int_0^t \frac{dsdt}{\lve t-s\rve^{1+(\alpha-\rho
)p}} \leq CM^p  T^{1+p(\rho-\alpha)}.
\end{eqnarray*}
In order to study the term $S_2$ let us  recall, see \cite[Theorem II.6.13]{516.47023}, that there exists a $C>0$ such that
\begin{equation}\label{est-auxi}
 \left|A^{-\gamma}\left(e ^ {-hA}  -I\right)\right|  _{\mathcal{L}(E)}\le C\: h
^ {\gamma},\quad \; h>0.
\end{equation}
Therefore, by applying  the Young inequality for convolutions we infer that
\begin{align*}
S_2\le & C {\mathbb{E}} \int_0 ^ T \int_ 0^t \frac{\Lve
A^{-\rho^\prime}(e^{-(t-s)A}-I)\Rve^p_{\mathcal{L}(E)} \lve
A^{\rho^\prime}
u(s)\rve^p}{\lve t-s\rve^{1+p\alpha}} ds dt\\
\le & C {\mathbb{E}} \int_0^T \lve A^{\rho^\prime} u(s)\rve^p
\biggl(\int_0^{T-s}\frac{\Lve A^{-\rho^\prime}(e^{-\tau
A}-I)\Rve^p_{\mathcal{L}(E)} }{\lve \tau
\rve^{1+p\alpha}}d\tau\biggr)ds\\
\le & C {\mathbb{E}}\Lve
A^{\rho^\prime}u\Rve^p_{L^p(0,T;E)}\int_0^T
\tau^{-1+p(\rho^\prime-\alpha)}d\tau\\
\le & C {\mathbb{E}}\Lve A^{\rho^\prime}u\Rve^p_{L^p(0,T;E)}
T^{p(\rho^\prime-\alpha)}.
\end{align*}
Invoking Lemma \ref{est-stoc-conv-1} \delb{we derive from the last inequality that
\begin{equation*}
\sup_{\xi\in\mathcal{B}_M(X) } S_2\le C T^{p(\rho^\prime-\alpha)},
\end{equation*}
for any $\alpha\in (0,\rho)$. This fact}  and the estimate for $S_1$
concludes  the proof of the lemma.
\end{proof}

\begin{remark}
Since $\alpha\in (0,\frac1p)$ we cannot (and rightly!) infer from the above lemma
that $\mathfrak{G}(\xi), \xi\in\mathcal{B}_M(X)$ has almost surely a
\cadlag  modification.
\end{remark}

The next three lemmata are about tightness of the family of laws of $\{\mathfrak{G}(\xi): \xi\in\mathcal{B}_M(X)\}$
\begin{lemma}\label{tight-sto-conv-1}
The family of laws of $\{\mathfrak{G}(\xi): \xi\in\mathcal{B}_M(X)\}$ is tight on $L^p(0,T;E)$.
\end{lemma}
\begin{proof}
As in Lemma \ref{lem-est-stoc-conv-2} we
choose an auxiliary  $\rho^\prime \in (0,\rho)$ and put $B=D(A^{\rho^\prime-1})$.
Let us also put $Y=D(A^{\rho^\prime})$. Since, by the assumptions \ref{assum-1}, $A$ has
compact resolvent,  it
follows from the combination of \cite[Proposition 5.8]{engel},
\cite[Theorem 1.15.3, pp 103]{Triebel_1995} and \cite[Theorem
1.16.4-2, pp 117]{Triebel_1995} that the embeddings $Y \embed E$ and $E\embed B$ are compact.
 Thanks to Lemma \ref{est-stoc-conv-1} and Lemma
\ref{est-stoc-conv-3} $\{\mathfrak{G}(\xi): \xi\in\mathcal{B}_M(X)\}$ is
uniformly bounded on
$\mathcal{M}^p(0,T;Y)\cap L^p(\Omega;W^{\alpha,p}(0,T;E))$.
Hence, since the embedding $$ W^{\alpha,p}(0,T;E)\cap L^p(0,T:D(A^{\rho^\prime}))\hookrightarrow L^p(0,T;E)$$  is
compact, see \cite[Step 1 of Proof of Theorem 2.1]{MR84d:34066},
it follows from the Chebyshev inequality and \cite[Theorem
2.1]{MR84d:34066} that the laws of $\{\mathfrak{G}(\xi):
\xi\in\mathcal{B}_M(X)\}$ are tight on $L^p(0,T;E)$.
\end{proof}

\begin{lemma}\label{tight-sto-conv-2}
\begin{enumerate}[]
 \item Let the assumption of Lemma \ref{lem-est-stoc-conv-2} hold. Then the family of laws of $\{\mathfrak{G}(\xi): \xi\in\mathcal{B}_M(X)\}$ is
tight on $\mathbb{D}([0,T];D(A^{\rho^\prime-1}))$ for any $\rho^\prime\in (0,\rho)$.
\end{enumerate}

\end{lemma}

For the proof of this lemma we need the following general result.
\begin{lemma}\label{GEN-TIGHT}
Assume that $p\in (1,2]$, $T>0$. Assume that $X$ and $Y$ are two
M-type $p$ Banach spaces such that the embedding $X\hookrightarrow Y$ is
compact.
For every $\xi\in\mathcal{B}_M(X)$ let a process   $v=\Lambda (\xi)$ be defined by
\begin{equation*}
v(t)=\int_0^t \int_Z \xi(s,z) \tilde{\eta}(dz,ds),\;\; t\in [0,T].
\end{equation*}
Then the family of laws of $\{v=\Lambda(\xi): \xi \in \mathcal{B}_M(X)\}$ on $\mathbb{D}([0,T];Y)$ is tight.
\end{lemma}
\begin{proof}
We need to check item (a) and (b) of Corollary \ref{comp-2}.

By the Burkholder inequality from \cite{maxreg} there exists $C>0$
such that for any $\xi\in \mathcal{B}_M(X)$ we have
\begin{equation*}
{\mathbb{E}}\sup_{s\in [0,T]}\lvert v(s)\rvert^p_X\le C M^p.
\end{equation*}
Let $\eps>0 $ and $K_\eps=\{ y \in X: \lvert y\rvert\le
(C\eps^{-1})^\frac1p M\}$. It follows easily from the Chebyshev inequality that
\begin{align*}
\mathbb{P}(v(t)\notin K_\eps, t\in [0,T])\le &[(C\eps^{-1})^\frac1p
M]^{-p} {\mathbb{E}}\sup_{t\in [0,T]}\lvert v(t)\rvert^p_X \le \eps.
\end{align*}
Corollary \ref{comp-2}-(a) follows from this inequality and the
compactness embedding $X\subset Y$.

Next, let us fix $0\le \sigma \le \tau \le T$. Then by \cite[Corollary C.2]{maxreg} and the Jensen inequality
\begin{align*}
{\mathbb{E}}\sup_{t\in [\sigma, \tau]}\lvert
v(t)-v(\sigma)\lvert_Y=&{\mathbb{E}}\sup_{t\in [\sigma, \tau]}\biggl\lvert
\int_\sigma^t\int_Z\xi(s,z)\tilde{\eta}(dz,ds)\biggr\rvert_Y
\\
\le & C{\mathbb{E}}\biggl( \int_\sigma^\tau\int_Z\rvert \xi(s,z)\rvert_Y^p
\nu(dz)ds\biggr)^\frac1p
\\
\leq & C
\biggl( {\mathbb{E}} \int_\sigma^\tau\int_Z\rvert \xi(s,z)\rvert_Y^p
\nu(dz)ds\biggr)^\frac1p
\le & C M (\tau-\sigma)^\frac1p.
\end{align*}
Thus Corollary \ref{comp-2}-(b) follows.
\end{proof}
\begin{proof}[Proof of Lemma \ref{tight-sto-conv-2}] Let us fix $\xi \in \mathcal{B}_M(X)$ and put $u=\mathfrak{G}(\xi)$.  Let us fix auxiliary numbers $\rho^\prime\in (0,\rho)$ and $\gamma \in (\rho^\prime,
\rho)$. Let us put   $\beta=\rho-\rho^\prime+1-\frac1p>0$ and   let us  rewrite identity  \eqref{est-conv-aux} as follows
\begin{align*}
A^{\rho^\prime-1}u(t)&=A^{\rho^\prime -\gamma}\int_0^t
A^{\gamma}u(s) ds+A^{-\beta}\int_0^t\int_Z
A^{\rho-\frac1p}\xi(s,z)\tilde{\eta}(dz,ds)\\
&=v_1(t)+v_2(t),\;\; t\in [0,T].
\end{align*}
It follows from Lemma \ref{GEN-TIGHT} that the family of laws of $\{v_2: \in \mathcal{B}_M(X)\}$ is tight on
$\mathbb{D}([0,T];E)$.\\
On the other hand, since $\gamma<\rho$, by Lemma \ref{est-stoc-conv-1} there exists $C>0$ such that
\begin{equation*}
{\mathbb{E}} \int_0^t \rvert A^\gamma u(s)\rvert^p ds \le C M^p.
\end{equation*}
Since $\rho^\prime-\gamma<0$, the map $A^{\rho^\prime-\gamma}:E \to E$ is compact. Therefore, since
  $v_1(t)=A^{\rho^\prime-\gamma}\int_0^t A^{\gamma}
u(s) \,ds $  we infer  that the family of laws of $\{v_1(\xi): \xi  \in \mathcal{B}_M(X) \} $ are tight on $C([0,T];E)$.
 Now we easily conclude the proof of (i) since by \cite[Theorem 4.1]{Whitt} the function $\phi: C([0,T];E)\times \mathbb{D}([0,T];E) \ni (x,y)  \mapsto x+y \in \mathbb{D}([0,T];E)$ is continuous.
\end{proof}
 %%%%%%
We also need the following auxiliary result.

\begin{lemma}\label{lem-tight-st0-conv-3}
Assume that the assumptions of part \eqref{lem-est-stoc-conv-2-dissipative} of Lemma \ref{lem-est-stoc-conv-2} are satisfied. Then for every $q\in (p,\frac{1}{\frac1p-\rho})=(p,\frac{p}{1-p\rho})$ and every  $r\in (1,p)$  there exists  $C>0$  such that
\begin{equation}\label{eqn-est-stoc-conv-4}
{\mathbb{E}} \lvert \mathfrak{G}(\xi)
\rvert^r_{L^{q} (0,T; E)}\le CM^r, \;\; \xi\in\mathcal{B}_M(X).
 \end{equation}
Moreover, the family of laws of $\{\mathfrak{G}(\xi): \xi \in\mathcal{B}_M(X)\}$ on $L^{q} (0,T; E)$ is tight.
\end{lemma}
\begin{proof}[Proof of Lemma \ref{lem-tight-st0-conv-3}]
Let us fix $q\in (p,\frac{1}{\frac1p-\rho})$. Then $q(\frac1p-\rho)<1$ and $q>p$. Thus there exists $\rho^\prime \in (0,\rho)$ such that $q(\frac1p-\rho)=1-p(\rho-\rho^\prime)$. Define
\[\theta=\frac{\frac1p-\rho}{\frac1p-\rho+\rho^\prime}=\frac{1-p\rho}{1-p\rho+p\rho^\prime}=\frac{p}{q}\in (0,1)
\]
and set, as in proof of Lemma \ref{tight-sto-conv-1}, $Y=D(A^{\rho^\prime})$. We also put $B=D(A^{\rho-\frac1p})$.
Then by the reiteration property of the complex interpolation we have,
\[
E=[B,Y]_\theta.
\]
Therefore,
\[
\vert y \vert_E \leq  \vert y \vert_B^{1-\theta} \,\vert y \vert_Y^{\theta}, \;\;\; y\in Y.
\]
Take any $r\in (1,p)$ and put $s=\frac{q}{r}$, $\frac1{s}+\frac{1}{s^\ast}=1$. Note that then $sr\theta=p$ and, since $r<p$, $rs^\ast(1-\theta)<p$. Let us choose an auxiliary $\delta>1$ such that $\delta rs^\ast(1-\theta)=p$.

Let us fix $\xi \in\mathcal{B}_M(X)$ and put $u=\mathfrak{G}(\xi)$. Then by the H\"older and Jensen inequalities
\begin{eqnarray*}
\mathbb{E} \vert u \vert_{L^q(0,T;E)}^r &\leq &
\mathbb{E}\Big[ \vert u \vert_{L^{q\theta}(0,T;Y)}^{r\theta}   \vert u \vert_{L^\infty(0,T;B)}^{r(1-\theta)}\Big]
=\mathbb{E}\Big[ \vert u \vert_{L^{p}(0,T;Y)}^{r\theta}   \vert u \vert_{L^\infty(0,T;B)}^{r(1-\theta)}\Big]
\\
&\leq & \mathbb{E}\Big[ \vert u \vert_{L^{p}(0,T;Y)}^{sr\theta}    \Big]^{\frac1s} \mathbb{E}\Big[ \vert u \vert_{L^\infty(0,T;B)}^{\delta rs^\ast(1-\theta)}\Big]^{\frac1{\delta s^\ast}}
\\&=&\mathbb{E}\Big[ \vert u \vert_{L^{p}(0,T;Y)}^{p}
\Big]^{\frac{r}{q}} \mathbb{E}\Big[ \vert u \vert_{L^\infty(0,T;B)}^{p}\Big]^{\frac{1}{\delta}-\frac{r}{\delta q}}\leq C_4M^r,
\end{eqnarray*}
where Lemma \ref{est-stoc-conv-1} and Lemma \ref{lem-est-stoc-conv-2}-\eqref{lem-est-stoc-conv-2-dissipative} were used to obtain the last inequality.
The proof of the  first part is complete. \\
To prove the second part we observe that by the same argument as above, given $q$ and $r$, we can find $\eps>0$ and $C>0$  such that
\begin{equation}\label{eqn-est-stoc-conv-4'}
{\mathbb{E}} \lvert \mathfrak{G}(\xi)
\rvert^r_{L^{q} (0,T; D(A^\eps))}\le CM^r, \;\; \xi\in\mathcal{B}_M(X).
 \end{equation}
Moreover, by Lemma \ref{est-stoc-conv-3}, for
any fixed  $\alpha\in (0,\rho)$,  we can find  $C_3>0$
such that
\begin{equation*}
{\mathbb{E}} \Lve
\mathfrak{G}(\xi)\Rve^p_{W^{\alpha,p}(0,T;E)}\le C_3 M^p,\;\;\ \xi \in
\mathcal{B}_M(X).
\end{equation*}
Now, since the embedding $D(A^\eps) \embed E$ is compact, the embedding
\[
L^{q} (0,T; D(A^\eps)) \cap W^{\alpha,p}(0,T;E) \embed L^{q} (0,T; E)
\]
is compact. Hence the second part of the Lemma follows.

\end{proof}
Next we will formulate an analogous result in the case when only the assumptions of parts
 \eqref{lem-est-stoc-conv-2-i} and \eqref{lem-est-stoc-conv-2-ii} of Lemma \ref{lem-est-stoc-conv-2} are satisfied. The proof will be similar with the difference that instead of taking $B=D(A^{\rho-\frac1p})$ we will need to take $B=D(A^{\rho^\prime-1})$.

\begin{lemma}\label{lem-tight-st0-conv-3b}
Let the assumptions of parts \eqref{lem-est-stoc-conv-2-i} and \eqref{lem-est-stoc-conv-2-ii} of Lemma \ref{lem-est-stoc-conv-2} be satisfied with
$\rho^\prime\in (0,\rho)$.
Then for every $q\in (p,\frac{p}{1-\rho^\prime})$ and every  $r\in (1,p)$  there exists $C>0$  such that
\begin{equation}\label{eqn-est-stoc-conv-4-a}
{\mathbb{E}} \lvert \mathfrak{G}(\xi)
\rvert^r_{L^q (0,T; E)}\le CM^r, \;\; \xi\in\mathcal{B}_M(X).
 \end{equation}
 Moreover, the family of laws of $\{\mathfrak{G}(\xi): \xi \in\mathcal{B}_M(X)\}$ on $L^{q} (0,T; E)$ is tight.
\end{lemma}
\begin{proof}[Proof of Lemma \ref{lem-tight-st0-conv-3b}]
Let $q\in (p,\frac{p}{1-\rho^\prime})$, $\rho^\prime\in (0,\rho)$ and define $\theta=\frac{p}{q}$.
As in proof of Lemma \ref{lem-tight-st0-conv-3} we let $Y=D(A^{\rho^\prime})$ and
$B=D(A^{\rho^\prime-1})$. Then by the reiteration property of the complex interpolation we have the continuous embedding,
$$[B,Y]_\theta=D(A^{\theta+\rho^\prime-1})\subset E. $$
Owing to Lemma \ref{est-stoc-conv-1} and parts \eqref{lem-est-stoc-conv-2-i} and \eqref{lem-est-stoc-conv-2-ii} of Lemma \ref{lem-est-stoc-conv-2} we can argue exactly as in the proof of Lemma \ref{lem-tight-st0-conv-3} and show that for any $r\in (1,p)$ we have
\begin{equation}\label{eqn-est-stoc-conv-4-b}
{\mathbb{E}} \lvert \mathfrak{G}(\xi)
\rvert^r_{L^q (0,T; D(A^{\theta+\rho^\prime-1}))}\le CM^r, \;\; \xi\in\mathcal{B}_M(X),
 \end{equation}
 which implies inequality \eqref{eqn-est-stoc-conv-4-a}. Thanks to \eqref{eqn-est-stoc-conv-4-b} we can again use the same argument as in proof of Lemma \ref{lem-tight-st0-conv-3}
  to deduce that the family of laws of $\{\mathfrak{G}(\xi): \xi \in\mathcal{B}_M(X)\}$ on $L^{q} (0,T; E)$ is tight.
\end{proof}
\section{Proof of Theorem \ref{Th:bound}}\label{proof-aux-result}
%
%%%%%%%%%%%%%%%%%%%%%%%%%%%%%%%%%%%%%%%%%%%%%%
%%%%%%%%%%%%%%%%%%%%%%%%%%%%%%%%%%%%%%%%%%%
%%%%%%%%%%%%%%%%%%%%%%%%%%%%%%%%%%%%%%%%%
Now, we will begin the proof of Theorem \ref{Th:bound} by first
defining a sequence of approximating processes. Let us fix for the whole section a number $T>0$. Consider a  sequence
$\{x_n\}\subset E$ such that $ x_n\to  u_0$ strongly in $D(A^{\rho-\frac1p})$
as $n\to\infty$. Define a function
$\phi_n:[0,\infty)\to[0,\infty)$ by $\phi_n(s)=\frac{k}{ 2^n}$, if
$k\in\mathbb{N}$ and $\frac{k}{ 2^n}\le s<\frac{k+1}{2^n}$, i.e.\
$\phi_n(s)= 2 ^ {-n}[2 ^ ns]$, $s\ge 0$, where $[t]$ is the
integer part of $t\in \mathbb{R}$.
 Let us define  a sequence $\{  u _n\}$ of adapted
$E$--valued processes by

\begin{eqnarray} \nonumber {u}_n(t) &=&e^{-tA}x_n+\int^t_0e^{-(t-s)A}F(s, \hat
{u}_n(s))\,ds
\\
&&{} +\int^t_0 \int_Z e^{-(t-s)A}G(s,\hat {u}_n(s);z)\,\tilde
\eta(dz;ds), \label{a}\quad t\in [0,T], \end{eqnarray} where $\hat u_n$ is
defined by \begin{eqnarray}\label{hatdefined} \hat u_n(s):= \bcase x_n, &
\mbox{ if } s\in [0,2^{-n}),
\\
 \dashint_{\phi_n(s)-2^{-n}}^{\phi_n(s)} u_n(r)\: dr, &\mbox{ if
} s\geq 2^{-n},\ecase \end{eqnarray}
and where we have used the following shortcut notation

\[
\dashint_A f(t) dt:= \frac{1}{\Leb(A)}\int_A f(t)\, dt, \quad A\in{\mathcal{B}}([0,T]).\]

({Here $\Leb$ denotes the Lebesgue measure.})
Note, that $\hat u$ is a
progressively measurable, piecewise  constant,  $E$-valued
process.
Between the grid points, equation \eqref{a} is linear, therefore,
$u_n$ is well defined for all $n\in{\mathbb{N}}$.

{Secondly, we prove some uniform estimates (w.r.t $n$) for the
solution $u_n$ of \eqref{a}. Recall that $F$ is a bounded
nonlinear map defined on $[0,T]\times E$ and taking values in $D(A^{\rho-1})$.
Furthermore, it is continuous w.r.t to the first and
second variables.}{
\begin{proposition}\label{unif-est-un}
For any $\alpha\in (0,\rho)$ and $\rho^\prime\in (0,\rho)$, there exists a constant $C$ such
that the following inequalities hold
\begin{align}
\sup_{n\in \mathbb{N}}{\mathbb{E}}\Lve A^{\rho^\prime}
u_n\Rve^p_{L^p(0,T;E)}\le
C,\label{est-un-1}\\
\sup_{n\in \mathbb{N}}{\mathbb{E}}\Lve A^{\rho^\prime}
\hat{u}_n\Rve^p_{L^p(0,T;E)}\le
C,\label{est-un-2}\\
\sup_{n\in \mathbb{N}}{\mathbb{E}}\Lve u_n\Rve^p_{W{\alpha,p}(0,T;E)}\le
C.\label{est-un-3}
\end{align}
\end{proposition}
\begin{proof}
Without loss of generality we take $T=1$. For each $n\in
\mathbb{N}$ we divide the interval [0,1] into small intervals of
length $2^{-n}$ each by setting: $I_k=[\frac{k}{2^{n}}, \frac{k+1}{2^{n}})$, $k=0,...,2^n-1$. We also put $J_0=I_0$ and $J_k=\bigcup\limits_{\ell=1}^k I_\ell$,
$k=1,...,2^n-1$. Define the sequences of processes $\{u_n^k: k=0,...,2^n -1\}$
and $\{\hat{u}^k_n: k=0,...,2^n-1 \}$ inductively by
\begin{equation*}
\begin{cases}
%\begin{split}
u_n^0(t)=e^{-t A}x_n+\int_0^t e^{-(t-s)A}F(s,x_n)ds\\
\quad \quad \quad +\int_0^t\int_Z e^{-(t-s)A}G(s,x_n) \tilde{\eta}(dz,ds),\;\; t\in I_0,
%\end{split}
\\
\hat{u}^0_n(t)=x_n ,\;\; t\in I_0;
\end{cases}
\end{equation*}
\begin{equation*}
\begin{cases}
%\begin{split}
u_n^k(t)=e^{-t A}x_n+\int_0^t e^{-(t-s)A}F(s,\hat{u}^k_n(s))ds\\
\quad  \quad\quad +\int_0^t\int_Z e^{-(t-s)A}G(s,\hat{u}^k_n(s))
\tilde{\eta}(dz,ds),\\
\quad \quad =:g_n(t)+y^k_n(t)+z^k_n(t), \text{ if } t\in  J_{k},
%\end{split}
\\
 \text{}
\\
 \hat{u}^k_n(t)= \begin{cases}
\hat{u}^{k-1}_n(t) \text{ if } t\in  J_{k-1},\\
\dashint_{I_{k-1}}u_n^{k-1}(s)ds \text{ if } t\in I_k,\\
\end{cases} \\
 k=1,\cdots,\,\,2^n-1.
\end{cases}
\end{equation*}
\delb{for any $t\in \bigcup\limits_{j=1}^k I_j$.} Note that, by
definition, $u_n^k$ is equal to the restriction of $u_n$ to $J_k$
and $ u_n=u_n^{2^n-1}$. Hence to prove our proposition it is
sufficient to check that the estimates
\eqref{est-un-1}-\eqref{est-un-3} are true and uniform w.r.t $k$
on $J_k$ for $u^k_n$ and
$\hat{u}^k_n$ with $k=0,...,2^n-1$. \\ On the interval $J_0$ we have
\begin{equation}
\label{eqn-10.1}
\begin{split}
u^0_n(t)=e^{-t A}x_n+\int_0^t e^{-(t-s)A}F(s,x_n)ds
\\ \quad \quad +\int_0^t\int_Z e^{-(t-s)A}G(s,x_n) \tilde{\eta}(dz,ds)\\
=g_n(t)+y^0_n(t)+z^0_n(t).
\end{split}
\end{equation}
First, it follows from \cite[Theorems 2 and 7]{Di-Blasio-84}  that  there exists a constant $C>0$ such that for any $n\in \mathbb{N}$
\begin{equation*}
%\begin{split}
{\mathbb{E}}\biggl[ \Lve A^{\rho} g_n\Rve^p_{L^p(J_0;E)}+\Lve
A^{\rho} y^0_n\Rve^p_{L^p(J_0;E)}\biggr]\le C \lve
x_n\rve^p+ C {\mathbb{E}}\Lve A^{\rho-1} F\Rve^p_{L^p(J_0;E)}
\le C.
%\end{split}
\end{equation*}
Secondly, we derive from \cite[Theorem 7 and 19]{Di-Blasio-84}  that for any  $\alpha\in
(0,\rho)$ there exists a constant $C>0$ such that for any $n\in
\mathbb{N}$
\begin{equation*}
%\begin{split}
{\mathbb{E}} \biggl[\Lve g_n\Rve^p_{W^{\alpha,p}(J_0;E)}+\Lve
y^0_n\Rve^p_{W^{\alpha,p}(J_0;E)}\biggr]\le C \lve x_n\rve^p+ C
{\mathbb{E}}\Lve
A^{\rho-1}F\Rve^p_{L^p(J_0;E)}
\le C.
%\end{split}
\end{equation*}
Hence combining these two remarks with Lemma
\ref{est-stoc-conv-1}, Lemma \ref{est-stoc-conv-3} and Proposition
\ref{convergen} we infer that Eqs.
\eqref{est-un-1}-\eqref{est-un-2} are true on $J_0$ for $u^0_n$
and $\hat{u}^0_n$. Using the same approach we can prove by
induction that for each $\alpha$ and $\rho$ as above there exists a constant $C>0$ such that for any $n\in \mathbb{N}$ and $k\in
\{0,...,2^n-1\}$ we have
\begin{equation*}
\begin{split}
{\mathbb{E}}\biggl[ \Lve A^{\rho} g_n\Rve^p_{L^p(J_k;E)}+ \Lve
A^{\rho} y^k_n\Rve^p_{L^p(J_k;E)}\biggr]\le C \lve
x_n\rve^p+ C {\mathbb{E}}\Lve A^{\rho-1}
F(.,\hat{u}^k _n)\Rve^p_{L^p(J_k;E)}\\
\le C,
\end{split}
\end{equation*}
 and
\begin{equation*}
\begin{split}
{\mathbb{E}}\biggl[ \Lve g_n\Rve^p_{W^{\alpha,p}(J_k;E)}+\Lve
y^k_n\Rve^p_{W^{\alpha,p}(J_k;E)}\biggr]\le C \lve x_n\rve^p+ C
{\mathbb{E}}\Lve A^{\rho-1}
F(.,\hat{u}^k_n)\Rve^p_{L^p(J_k;E)}\\
\le C.
\end{split}
\end{equation*}
With the same argument as above we check that Eqs.
\eqref{est-un-1}-\eqref{est-un-3} are correct and uniform w.r.t
$k$ on each $J_k$ with $k=0,...,2^n-1$. With this fact and the
identity $u_n=u^{2^n-1}_n$, we conclude the proof of our
proposition.
\end{proof}
We will also need the following result.
\begin{proposition}\label{prop-initial}
Suppose that $(x_n)_n\subset E$ is a sequence such that $$\lve
A^{\rho-\frac1p}[x_n-u_0]\rve\rightarrow 0 ,$$ as $n\rightarrow
\infty$ .
Then the sequence $\{ g_n: n\in{\mathbb{N}}\}$ defined by \eqref{eqn-10.1}  is convergent (and hence is
precompact) in the following space  $C([0,T]; D(A^{\rho-\frac1p})) \cap L ^ p(0,T;E)$.
\end{proposition}%
\begin{proof}[Proof of Proposition \ref{prop-initial}] The convergence in $C([0,T]; D(A^{\rho-\frac1p}))$ is obvious.
From \cite[Theorem 2]{Di-Blasio-84} we infer that for any
$\theta\in[\frac1p-\rho,\frac1p]$ there exists some $C>0$ such
that
\begin{equation}
 \sup_{n\in \mathbb{N}}\Lve A^{\theta}[g_n-e^{-tA}u_0]\Rve^p_{L^p(0,T;E)}\le C \lve A^{\rho-\frac1p}[x_n-x]\rve^p,
\end{equation}
from which we derive the convergence in $L^p(0,T;E)$.
\end{proof}
}

 After these preliminary claims we are now ready for the proof of Theorem \ref{Th:bound} which will be divided
into several steps. But before we go further let us define a
sequence of Poisson random measures $\{\eta_n \}_{n\in{\mathbb{N}}}$ by
putting $ \eta_n=\eta$ for all $n\in{\mathbb{N}}$.
\usecounter{lil}
\begin{steps}\label{step1}
\usecounter{lil}\setcounter{lil}{0}
\item The family of laws of  $\{  (u_n, \eta_n):n\in{\mathbb{N}}\}$ is
tight on \\ $\big[ L ^ p(0,T;E) \cap \mathbb{D}([0,T];D(A^{\rho^\prime-1}))\big] \times M_{\mathbb{N}}(Z\times [0,T])$, for any $\rho^\prime\in (0,\rho)$.
\end{steps}
\begin{proof} To simplify notation we set $B_0=D(A^{\rho^\prime-1})$ for any $\rho^\prime\in (0,\rho)$. Define three functions $f_n$, $g_n$ and $v_n$ by
\begin{eqnarray}\label{def-f-n}
 {f}_n(t)=F(t,\hat {u}_n(t)),\;\; t\in [0,T]
\end{eqnarray} \begin{eqnarray}\label{def-g-n}
 {g}_n(t;z)=G(s,\hat {u}_n(t);z), \quad t\in [0,T],\quad z\in Z
\end{eqnarray} and \begin{eqnarray}\label{def-v-n} v_n(t)=\int^t_0\int_Z e^{-(t
-s)A}G(s,\hat {u}_n(s); z)\,\tilde \eta(dz;ds). \end{eqnarray}
We argue exactly as in \cite{Brz+Gat_1999}. We recall that space ${\mathcal{M}}^ p (0,T;E)$, $\Lambda$ and $\mathcal{A}$ are defined on page
 \pageref{def-notation-2} and \pageref{operator_def}, respectively. Since, by estimates
\eqref{est-un-2} and Assumption \ref{assum-2}, the family $\{ A ^
{\rho-1} f_n : n\in{\mathbb{N}}\}$ is bounded in ${\mathcal{M}}^ p (0,T;E)$
and $\mathcal{A}\Lambda^{-1}$ is bounded on $L^p(0,T; E)$ it
follows from \cite[Theorem 2.6]{Brz+Gat_1999} and Corollary
\ref{C:comp} that $\Lambda^{-1} f_n=\Lambda^{-\rho}
(\mathcal{A}\Lambda^{-1})^{1-\rho} A^{\rho-1} f_n$
is tight on $L ^ p (0,T;E)\cap C([0,T];E)$. This fact, the compact
embedding $E\subset B_0$ and the continuity of the embedding
$$C([0,T];B_0)\subset \mathbb{D}(0,T;B_0)$$ imply that $\Lambda^{-1}f_n$
is tight on $L ^ p (0,T;E)\cap \mathbb{D}([0,T];B_0)$. Next by estimates
\eqref{est-un-2}, Lemma \ref{tight-sto-conv-1} and Lemma
\ref{tight-sto-conv-2}, we infer that  the laws of the family $\{ v_n: n\in{\mathbb{N}}\}$
are tight on  $L ^ p(0,T;E)\cap \mathbb{D}([0,T];B_0)$. Finally, from
Proposition \ref{prop-initial} it follows that  the family of
functions $\{ e ^ {-\cdot  A}x_n:n\in{\mathbb{N}}\}$ is precompact in $ L
^ p (0,T;E) \cap \mathbb{D}([0,T];B_0)$. Since
$$   u_n = v_n+\Lambda^{-1}{   f}_n+ e ^ {-\cdot A}x_n ,\quad n\in{\mathbb{N}},
$$
we easily conclude that the laws of the family  $\{     u _n :
n\in{\mathbb{N}}\}$ are tight on $  L ^ p(0,T;E)\cap\mathbb{D}([0,T];B_0)$. Since
$M_{\mathbb{N}}(Z\times [0,T])$ is a separable metric space, by
\cite[Theorem 3.2]{para} the laws of the family
$\{\eta_n:n\in{\mathbb{N}}\}$ are  tight on $M_{\mathbb{N}}(Z\times [0,T])$.
Consequently The laws of the family  $\{  (u_n, \eta_n):n\in{\mathbb{N}}\}$
are tight on $L ^ p(0,T;E)\cap\mathbb{D}([0,T];B_0)\times M_{\mathbb{N}}(Z\times
[0,T])$.
\end{proof}
From \textbf{Step (I)}  and Prokhorov Theorem (see, for instance, \cite[Theorem 2.3]{761.60052}) we deduce that there
exist   a subsequence of $\{ ( u_n, \eta_n) , n\in{\mathbb{N}}\}$, still
denoted by $\{ ( u_n, \eta_n):n\in{\mathbb{N}}\}$, and a Borel probability
measure  $\mu_\ast$ on $ \left[\mathbb{D}([0,T];D(A^{\rho^\prime-1}))\cap L ^
p(0,T;E)\right]\times M_{\mathbb{N}}(Z\times [0,T]) $ such that
 ${\mathcal{L}}( u_n, \eta_n )\to \mu_\ast $
weakly.
By Theorem \ref{thm-Skorokhod} there exists a probability space
$(\bar \Omega,\bar {\mathcal{F}},\bar{\mathbb{P}})$ and \\ $  L ^ p(0,T;E)\cap
\mathbb{D}([0,T];D(A^{\rho^\prime-1}))\times M_{\mathbb{N}}(Z\times [0,T])$-valued random variables
$(\bar u_1,\bar \eta_ 1)$, $(\bar u_2,\bar \eta_2)$, $\ldots $,
having the same law as the random variables $( u_1, \eta_ 1)$, $(
u_2,\eta_2)$, $\ldots $, and a $  L ^ p(0,T;E)\cap
\mathbb{D}([0,T];D(A^{\rho^\prime-1}))\times M_{\mathbb{N}}(Z\times [0,T])$-valued random variable
$(u_\ast,\eta_\ast)$ on $(\bar \Omega,\bar {\mathcal{F}},\bar{\mathbb{P}})$ with
${\mathcal{L}}((u_\ast,\eta_\ast))=\mu_\ast$ such that $\mathbb{P}$-a.s.
\begin{eqnarray}\label{limita-s}
\lefteqn{ %\lim_{n\to \infty }
\lim_{n\to\infty} \left( \bar u _n,\bar \eta_ n \right) =  \left( u
_\ast,\eta _\ast\right) } &&
\\&&\nonumber
 \mbox{ in }   L ^ p(0,T;E)\cap
\mathbb{D}([0,T];D(A^{\rho^\prime-1}))\times M_{\mathbb{N}}(Z\times [0,T]),
\end{eqnarray} and $\bar \eta_n=\eta_\ast$ for all $n\in{\mathbb{N}}$. The sequence
$\{\bar{u}_n; n\in \mathbb{N}\}$ satisfies the same properties as the original sequence $(u_n)$, some of these are stated in part (i) of the next step.
\begin{steps}
\addtocounter{lil}{1}
\item The following holds
\begin{steps-2}
\item $ \sup_{n\in{\mathbb{N}}} \left\|\bar u_n
\right\|_{{\mathbf{L}}^p(\bar{\Omega}\times [0,T];E)}<\infty $ and
\item  for any $r\in (1,p)$ we have
$$\lim_{n\to \infty}\bar{{\mathbb{E}}}\left\|\bar
u_n- u_\ast \right\|^r_{{\mathbf{L}}^p(0,T;E)} = 0. $$
\end{steps-2}
 \end{steps}
 \begin{proof}
Let us recall that by our construction which used the Skorokhod
embedding Theorem, the laws of $u_n$ and $\bar u_n$ on $L ^
p(0,T;E)$ are identical
 for any $n\in{\mathbb{N}}$.
Hence, $\|u_n\|_{{\mathbf{L}}^p(\bar{\Omega}\times [0,T];E)}= \|\bar
u_n\|_{{\mathbf{L}}^p(\bar{\Omega}\times [0,T];E)}$ and part (i) easily follows from
estimates \eqref{est-un-1}. It follows from part (i) that $\Lve \bar u_n \Rve_{{\mathbf{L}}^r(0,T;E) }$ is uniformly integrable with respect to  the probability measure $\bar{\mathbb{P}}$. Since $ \bar{\mathbb{P}}$-a.s.\ $\bar u_n\to
u_\ast$ and $\Lve u_n \Rve_{{\mathbf{L}}^r(0,T;E) }$ is uniformly integrable wrt the probability measure $\bar{ \mathbb{P}}$, we obtain part (ii) owing to the applicability of the Vitali Convergence Theorem.
 \end{proof}
Before we continue we should note that the random variables
$\bar{u}_n, u_\ast: \bar{\Omega} \to L^p(0,T;E)$ induce two
$E$-valued stochastic processes still denoted with the same
symbols, see for example \cite[Proposition B.4]{ZB+MO} for a proof
for the space $L^\infty_{loc}(R_+; L^2_{loc}(\mathbb{R}^d))$. Now, let
$\bar {\mathbb{F}}=(\bar {\mathcal{F}}_t)_{t\ge 0}$
be the
filtration defined by
\begin{equation} \label{eq:filtrat}
 \bar {\mathcal{F}}_t= \sigma( {\bar\eta_n}(s) ,\{u_m(s), m\in{\mathbb{N}}\},\, u_\ast(s); 0\le s\le t),\quad t\in[0,T].
\end{equation}
Since $\bar \eta_n=\eta_ \ast$, it is easy to show that the filtration obtained by
replacing $\eta_n$ with $\eta_\ast$ in
\eqref{eq:filtrat} is the equal to $\bar {\mathbb{F}}$.

The next two
steps imply that the following two $E$-valued integrals over the filtered
probability
space $(\bar \Omega,\Bar {\mathcal{F}},\bar{\mathbb{F}}, %(\bar {\mathcal{F}}_t)_{t\ge 0},
\bar{\mathbb{P}})$
$$
\int_0 ^ t  \int_Z e ^ {-(t-s)A}\: G(s,\bar u_n(s),z))\,
\tilde{\bar \eta}_n(dz,ds),\quad t\ge 0,
$$
and
$$
\int_0 ^ t  \int_Z e ^ {-(t-s)A}\: G( s,u_\ast(s),z))\, \tilde{
\eta}_\ast(dz,ds),\quad t\ge 0,
$$
do exist.
\pagebreak[3]
\begin{steps}
\addtocounter{lil}{1}
 \item {\samepage The following holds
\begin{steps-2}
\item for every $n\in{\mathbb{N}}$, $\bar \eta_n$ is a time homogeneous
Poisson random measure on ${\mathcal{B}}(Z)\times {\mathcal{B}}([0,T])$ over $(\bar
\Omega,\Bar {\mathcal{F}},\bar {\mathbb{F}} ,\bar{\mathbb{P}})$ with intensity measure $\nu$;
\item $\eta_\ast$ is a time homogeneous Poisson random measure on
${\mathcal{B}}(Z)\times {\mathcal{B}}([0,T])$ over $(\bar \Omega,\Bar {\mathcal{F}},\bar
{\mathbb{F}},\bar{\mathbb{P}})$ with intensity measure $\nu$;
\end{steps-2}}
\end{steps}
%%%%%%%%%%%%%
\begin{proof}
 Before embarking on the proof, let us first recall that the
modified version of the Skorokhod embedding Theorem, i.e.\ Theorem
\ref{thm-Skorokhod}, implies that
$\bar\eta_n(\bar \omega)=\eta^\ast(\bar \omega)$ for all $\bar \omega\in \bar \Omega$
and $n\in{\mathbb{N}}$.

For a random measure $\mu$ on $S\times [0,T]$ and for any $A\in {\mathcal{S}} $ %( \times {\mathcal{B}}(0,T)$
let us define  an  $\bar{\mathbb{N}}$-valued process
$(N_\mu(t,A))_{t\ge 0}$ by
$N_\mu(t,A):= \mu(A \times (0,t]), \;\; t\ge 0.
$
In addition, we denote by  $(N_\mu(t))_{t\ge 0}$ the measure
valued process defined by $N_\mu(t) =\{ {\mathcal{S}} \ni A \mapsto
N_\mu(t,A)\in\bar{\mathbb{N}}\}$, $t\in[0,T]$.

{\sl Proof of \textbf{Step (III)}-(i):} Since $\bar\eta_n$ and $\eta_n$ have the same law and $\eta_n$ is a time homogeneous Poisson
 random measure, it follows from Proposition \ref{PROP-EQ-PRM} and Remark \ref{REM-EQ-PRM} that
 $\bar \eta_n$ satisfies Definition \ref{def-Prm} \eqref{prm-i}-\eqref{prm-3}. Therefore, in order to prove
 part (i) of \textbf{Step (III)} we only need to prove that $\bar \eta_n$ satisfies Definition \ref{def-Prm} \eqref{prm-4}
  with the filtration defined in \eqref{eq:filtrat}.
  For this purpose let us fix $m\in \mathbb{N}$, $t_0\in [0,T]$ and $r\ge s\ge t_0$.
  It follows from the definition of $\bar{{\mathbb{F}}}$ that $\bar \eta_n$ is $\bar{{\mathbb{F}}}$-adapted and it remains to prove that
  $\bar X_m= N_{\bar \eta_n}(r)-N_{\bar \eta_n}(s)$ is independent of
  $\bar{\mathcal{F}}_{t_0}$.
  By Definition \ref{def-Prm} \eqref{prm-ii} the random variable $\bar X_m= N_{\bar \eta_n}(r)-N_{\bar \eta_n}(s)$ is independent of
  $N_{\bar \eta_n}(t_0)$, so we only need to show that $\bar X_m$ is independent of $\bar u_m(\sigma)$ and $u_\ast(\sigma)$ for any $\sigma\le t_0$.
  In what follows we also fix $\sigma \in [0,t_0]$. Since $\mathcal{L}(\bar{u}_m, \bar{\eta}_m)=\mathcal{L}(u_m, \eta_m)$, it follows that
  \begin{equation}\label{EQ-LAW}
\mathcal{L}(\bar{u}_m|_{[0,\sigma]}, \bar{X}_m)=\mathcal{L}(u_m|_{[0,\sigma]}, X_m),
  \end{equation}
 where $X_m=N_{\eta_m}(r)-N_{\eta_m}(s)$.
  Recall that $\eta_m=\eta_\ast$ and $u_m$ is the unique solution to the linear stochastic evolution equation \eqref{a}, hence it is adapted to
  the $\sigma$-algebra generated by
  $\eta_m$. Consequently, $u_m|[0,\sigma]$ is independent of $X_m$ and we infer from this last remark and \eqref{EQ-LAW} that $\bar{u}_m|_{[0,\sigma]}$ is independent of
  $\bar{X}_m$. It remains to prove that $\bar{X}_m$ is independent of $u_\ast|[0,\sigma]$, but, this is a subject
of  the next Lemma.
\begin{lemma}\label{independenc}
Let $B$ be a Banach space, $z$ and $y_\ast$ be two B-valued random
variables over $(\Omega,{\mathcal{F}},\mathbb{P})$. Let $\{y_n:n\in{\mathbb{N}}\}$  be a
family of $B$-valued random variables over a probability space
$(\Omega,{\mathcal{F}},\mathbb{P})$ such that $y_n\to y_\ast$ weakly, i.e.\ for all
$\phi\in B ^ \ast$, $ \mathbb{E} e ^ {i\la \phi,y_n\ra }\to \mathbb{E}   e ^ {i\la
\phi,y\ra }$.
If for all $n\ge 1$ the two random variables $y_n$ and $z$ are
independent, then $z$ is also independent of $y_\ast$.
\end{lemma}
\begin{proof}[Proof of Lemma \ref{independenc}]
The random variables $y_\ast$ and $z$ are independent iff \begin{eqnarray*}
{\mathbb{E}} e ^{i(\theta_1 z + \theta_2 y_\ast )} = {\mathbb{E}} e ^{i\,\theta_1
z}\,{\mathbb{E}} e^{ i\theta_2 y_\ast }, \quad \theta_1,\theta_2\in B  ^
\ast. \end{eqnarray*}
The weak convergence and the independence of $z$ and $y_n$ for all $n\in{\mathbb{N}}$ justify the following chain of equalities.
\begin{eqnarray*} {\mathbb{E}} e ^{i(\theta_1 z + \theta_2 y_\ast )} =
\lim_{n\to\infty}{\mathbb{E}} e ^{i(\theta_1 z + \theta_2 y_n )} =
\lim_{n\to\infty}{\mathbb{E}} e ^{i\theta_1 z}\,{\mathbb{E}} e^{ \theta_2 y_n}={\mathbb{E}} e
^{i\theta_1 z}\,{\mathbb{E}} e^{ \theta_2 y_\ast}. \end{eqnarray*}
\end{proof}
Since $\bar u_m|[0,\sigma]$ is independent
from $ \bar{X}_m$, Lemma \ref{independenc} implies that $u_\ast|[0,\sigma]$ is independent from $\bar{X}_m$.

{\sl Proof of \textbf{Step (III)}-(ii):} We have to show that
$\eta_\ast\in{\mathcal{M}}_I(S\times {\mathbb{R}} ^ +)$ is a time homogeneous Poisson random measure
with intensity $\nu$. But this will follow from \textbf{Step
(III)}-(i), since $\eta_\ast(\omega)=\bar \eta_m(\omega)$ for all
$\omega\in\Omega$ and $m\in \mathbb{N}$.
\end{proof}
\begin{steps}
\addtocounter{lil}{1} \item {\samepage The following  holds
 \begin{steps-2} \item  for all $n\in{\mathbb{N}}$, $\bar u_n$ is a  $\bar{\mathbb{F}}$-progressively measurable process;
\item the process $ u_\ast$ is a
$\bar{\mathbb{F}} $
-progressively measurable process.
\end{steps-2}}  \end{steps}
\pagebreak[3] %%%%%%%%%%%%%%%
\begin{proof}As we noted earlier, one can argue as in \cite[Proposition B.4]{ZB+MO} and prove that
the random variables $\bar{u}_n, u_\ast: \bar{\Omega} \to
{\mathbf{L}}^p(0,T;E)$ induce two $E$-valued stochastic processes still
denoted with the same symbols. Here, we have to show that for each
$n\in{\mathbb{N}}$, $\bar u_n$ and
$u_\ast$ are
$\bar {\mathbb{F}}$-progressively measurable. By definition of $\bar {\mathbb{F}}$,
for fixed $n\in{\mathbb{N}}$ the process $\bar u_n$ is adapted to $\bar
{\mathbb{F}}$ by the definition of $\bar {\mathbb{F}} $.
Let us fix $r\in (1,p)$. By Step (II) the process  $\bar u_n$ is bounded in
${\mathbf{L}}^r(\bar{\Omega}_T;E)$, hence,  there exists a sequence of
simple functions $\{\bar u_n^m,m\in{\mathbb{N}}\}$ such that $\bar u_n^m\to
\bar u_n$ as $m\to\infty$ in
 ${\mathbf{L}}^r(\bar{\Omega}_T;E)$. In particularly, by using the shifted Haar projections
used in  \cite[Appendix B]{uniqueconv} we can choose $\{\bar
u_n^m,m\in{\mathbb{N}}\}$ to be progressively measurable. It follows that
$\bar u_n$ is progressively measurable as a
${\mathbf{L}}^r(\bar{\Omega}_T;E)$-limit of a sequence of progressively
processes. Since $\bar u_n^m\to \bar u_n$ as $m\to\infty$ in
${\mathbf{L}}^r(\bar{\Omega}_T;E)$ and $\bar u_n\to u_\ast$ as $n\to\infty$
also in ${\mathbf{L}}^r(\bar{\Omega}_T;E)$, it follows that $u_\ast$ is a
limit  in ${\mathbf{L}}^r(\bar{\Omega}_T;E)$ of some progressively
measurable step functions. In particular, $u_\ast$ is also
progressively measurable, i.e.\ (ii) holds.
\end{proof}
%%%%%%%%%%%%%%%%%%%%%%%%%%%%%%%%%%%%%%%%%
  Let $\mu$ be a time homogeneous Poisson random measure
over $(\bar \Omega,\bar {\mathcal{F}},\bar {\mathbb{F}},\bar{\mathbb{P}})$ with intensity
measure $\nu$, $v$ be an $E$-valued progressively measurable process, $u_0 \in
D(A^{\rho-\frac1p})$ and ${\mathcal{K}}$ be a nonlinear map defined by
\begin{eqnarray}\label{definition-ck}\nonumber
 \lefteqn{ {\mathcal{K}} (x,v,\mu) (t) := e ^ {-tA} u_0 + \int_0 ^ t e ^
{-(t-s)A}\: F(s,v  (s))\: ds  {}}\\&&{}+
  \int_0 ^ t  \int_Z e ^ {-(t-s)A}\: G(s,v(s);z))\tilde \mu(dz,ds)
,\,\quad  t\in [0,T].
\end{eqnarray}
Here, as usual, $\tilde \mu$ denotes the compensated Poisson
random measure of $\mu$.
%\
%\
\begin{steps}
\addtocounter{lil}{1} \item For all $t\in [0,T]$ and $n\in
\mathbb{N}$ we have $\bar{\mathbb{P}}$-almost surely  $$\bar
u_n(t)-{\mathcal{K}} (x_n,\hat{\bar{u}}_n,\bar \eta_n )(t)=0,$$ where $\hat
{\bar{u}}_n$ is defined by \begin{eqnarray}\label{hatdefined-123} \hat
{\bar{u}}_n:= \bcase x_n, & \mbox{ if } s\in [0,2^{-n}),
\\
2 ^ n \int_{\phi_n(s)-2^{-n}}^{\phi_n(s)} \bar{u}_n(r)\: dr,
&\mbox{ if } s\geq 2^{-n}.\ecase \end{eqnarray}
\end{steps}
\begin{proof}
First, let
$$\mathfrak{S}_1=L^p(0,T;E)\cap L^\infty(\mathbb{R}_+; D(A^{\rho^\prime-1}))$$ and
$$\mathfrak{S}_2= M_{\mathbb{N}}(Z\times [0,T]).$$ Again for simplicity we set $B_0=D(A^{\rho^\prime-1})$ for any $\rho^\prime \in (0,\rho)$.
It is
proved in \cite{unique} that
the map  %have the same law.
 ${\mathcal{G}}:\mathfrak{S}_1\to \mathfrak{S}_1$ defined by %eqref{hatdefined}
\begin{equation*} {\mathcal{G}}(u)(s) :=  \bcase x_n, & \mbox{ if }
s\in [0,2^{-n}),
\\
2 ^ n \int_{\phi_n(s)-2^{-n}}^{\phi_n(s)} u(r)\: dr, &\mbox{ if }
s\geq 2^{-n}\ecase \end{equation*} is linear and bounded from
$\mathfrak{S}_1$ into itself. Therefore, for any $n\in{\mathbb{N}}$, the
two triplets of random variables $\{(u_n,\eta_n,\hat u_n)\}$ and
$\{(\bar u_n,\bar \eta_n,\hat{\bar u}_n)\}$, where $ \hat u_n=
{\mathcal{G}}(u_n))$ and $\hat{\bar u}_n= {\mathcal{G}}(\bar u_n)$, have equal laws on
$\mathfrak{S}\times \mathfrak{S}_1$. Second, let us define
processes $\tilde z_n$ and $ \breve z_n $ by
 \begin{eqnarray}
\label{abb}\nonumber
 \lefteqn{ \tilde z_n (t) := e ^ {-tA} x_n + \int_0 ^ t e ^
{-(t-s)A}\: F( u_n  (s))\: ds {}}
 && \\&&{} +
  \int_0 ^ t  \int_Z e ^ {-(t-s)A}\: G(u_n (s,z))\tilde \eta(dz,ds)
,\,t\in [0,T],
\\
\nonumber
 \lefteqn{ \breve z_n (t) := e ^ {-tA} x_n + \int_0 ^ t e ^
{-(t-s)A}\: F( \hat u_n  (s))\: ds {}}
 &&\\  \label{abb2} &&{} +
  \int_0 ^ t  \int_Z e ^ {-(t-s)A}\: G(\hat u_n (s,z))\tilde \eta(dz,ds)
,\, t\in [0,T].
\end{eqnarray}
 Let us also define processes $\tilde {\bar z}_n$ and
$\breve {\bar z}_n$ by replacing $(u_n, \eta)$ and $(\hat u_n,
\eta)$
 by $(\bar{u_n}, \bar{\eta}_n)$ and $(\hat {\bar u}_n, \bar{\eta}_n)$ in
formula \eqref{abb} and \eqref{abb2}, respectively.
Since $F\circ {\mathcal{G}}$ and
$G\circ {\mathcal{G}}$ are continuous, it follows from \cite[Theorem
1]{uniqueconv} that the quintuples of random variables
$\{(u_n,\eta_n,\hat u_n,\tilde z_n,\breve z_n)\}$ and $\{(\bar
u_n,\bar \eta_n,\hat{\bar u}_n,\tilde {\bar z}_n,\breve {\bar z}_n
)\}$ have equal laws on $\mathfrak{S}\times \mathfrak{S}_1\times
\mathfrak{S}_1\times \mathfrak{S}_1$. Consequently $\Psi(u_n,
\breve{z}_n)$ and $\Psi(\bar{u}_n, \breve{\bar{z}}_n)$ have equal
laws on $\mathfrak{S}_1\times\mathfrak{S}_1$, where the continuous
functional $\Psi$ is defined by
\begin{equation*}
 \Psi(v,w)=\int_0^T \lvert
v(t)-w(t)\rvert_{B_0} \,\, dt , \text{ for } v\in \mathfrak{S} \text{ and
} w\in \mathfrak{S}.
\end{equation*}
The above fact implies that
for any any real-valued $\varphi\in C(\mathbb{R}_+)$ we have
\begin{equation}\label{law-eq-funct}
\En[\varphi(\Psi(\bar{u}_n,
\breve{{z}}_n))]=\bar{{\mathbb{E}}}[\varphi(\Psi(\bar{u}_n,
\breve{\bar{z}}_n))].
\end{equation}
 Now let $\eps$ be an
arbitrary positive number and $\phi_\eps\in C(\mathbb{R}_+)$
defined by
\begin{equation*}
\phi_\eps(y)=\begin{cases} \frac{y}{\eps} \text{  if  } y \in
[0,\eps),\\
\mathbf{1}_{[\eps,\infty)}(y) \text{   otherwise.}
\end{cases}
\end{equation*}
%%%%%%%%%%%%%
It is easy to check that
\begin{equation*}
\begin{split}
\bar{\mathbb{P}}\left(\Psi(\bar{u}_n, \breve{\bar{z}}_n)\ge
\eps\right) &\le \int_{\bar{\Omega}}
\mathbf{1}_{[\eps,\infty)}(\Psi(\bar{u}_n, \breve{\bar{z}}_n))
d\bar{\mathbb{P}} \\ & \quad +
\int_{\bar{\Omega}}\mathbf{1}_{[0,\eps)}(\Psi(\bar{u}_n,
\breve{\bar{z}}_n)) \frac{\Psi(\bar{u}_n
\breve{\bar{z}}_n)}{\eps}d\bar{\mathbb{P}}\\
& =\bar{\mathbb{E}}\phi_\eps(\Psi(\bar{u}_n, \breve{\bar{z}}_n)).
\end{split}
\end{equation*}
The last inequality altogether with \eqref{law-eq-funct} implies that
\begin{equation}\label{eq-gal}
\bar{\mathbb{P}}\left(\Psi(\bar{u}_n, \breve{\bar{z}}_n)\ge
\eps\right)\le \mathbb{E}_n\phi_\eps(\Psi(u_n,\breve{z}_n)).
\end{equation}
%%%%%%%%%%
 Since for any $t$ and $\mathbb{P}$-almost surely $u_n(t)-\breve{z}_n(t)=0$
 we obtain that $\bar{\mathbb{P}}$ almost surely $\Psi(u_n,\breve{z}_n)=0$ ,
which along with \eqref{eq-gal} yield that for any $\eps>0$
\begin{equation*}\label{eq-gal2}
\bar{\mathbb{P}}(\Psi(\bar{u}_n, \breve{\bar{z}}_n)\ge \eps)=0.
\end{equation*}
Since $\eps>0$ is arbitrary,  we infer from the last equation that
\text{ $\bar{\mathbb{P}}$-a.s. ,}
\begin{equation*}
\Psi(\bar{u}_n, \breve{\bar{z}}_n)=0.
\end{equation*}
This implies that for $\bar{\mathbb{P}}$ almost all $t\in [0,T]$ and almost surely
$\bar{u_n}(t)=\breve{\bar{z}}_n(t)$. Since two \cadlag \,functions
which are equal almost all $t\in [0,T]$ must be equal for all
$t\in [0,T]$, we derive that almost
surely
$$\bar{u}_n=\mathcal{K}(x_n, \hat{\bar{u}}_n,
\tilde{\bar{\eta}}_n),$$
for all $t\in [0,T]$ .
\end{proof}

\begin{steps}
\addtocounter{lil}{1} \item We have
$$\left\|{\mathcal{K}} (x_n,\hat{\bar{u_n}},\bar \eta_n)-{\mathcal{K}}(u_0 ,u_\ast,\eta_\ast )\right\|_{{\mathbf{L}}^p(\bar{\Omega}\times[0,T];E)}
\to 0,\mbox{ as } n\to\infty.
$$
\end{steps}
%%%%%%%%%%%%%%%%%%%%%%
\begin{proof}
First, notice that since $\bar{\eta}_n=\eta_\ast$ for any $n\in
\mathbb{N}$, the convergence in \textbf{Step (VI)} is equivalent
to
$$
\bar {\mathbb{E}} \left\| {\mathcal{K}}( x_n,\hat{\bar{u _n}},\eta _\ast)-{\mathcal{K}}( u_0,u _\ast,
\eta _\ast)\right\|^{p} _{{\mathbf{L}}^p(0,T;E)}\to 0 \mbox{ as   }
n\to\infty.
$$

Observe, that  for any any $n$
%%%%%%%%%%%%%%%%%
\begin{eqnarray*}\nonumber \lefteqn{\bar{{\mathbb{E}}} \left\| {\mathcal{K}}(
x_n,\hat{\bar{u}}_n,\eta_\ast)-{\mathcal{K}}( u_0 ,u _\ast,\eta_
\ast)\right\|_{{\mathbf{L}}^p(0,T;E)}^ {p  }
 \le C\;  \Big\{\bar{\mathbb{E}} \int   _0 ^ T \left| e ^ {-tA}\left[ x_n-u_0\right]\right| ^ {p  }_E
 \: dt{} } &&
\\ &+&
\bar{\mathbb{E}}  \int _0  ^ T \left| \int  _0  ^ t
 e ^ {-(t-s) A} \left[  F(s, \hat{\bar{u}}_n(s))- F(s,u _\ast (s))\right] \: ds \right|_E
^ {p  }\: dt {}
\\
\lefteqn{{}+
 \bar{\mathbb{E}} \int_0 ^ T \Big| \int_0 ^ t\int_Z
 e ^ {-(t-s) A}  G(s,\hat{\bar{u}}_n(s);z)  \tilde{\eta}_\ast(dz,ds) {}  }&&
 \\ && \qquad\qquad\qquad{}   -\int_0 ^ t\int_Z
 e ^ {-(t-s) A}  G(s,u_{\ast}(s);z) \tilde{ \eta}_\ast(dz,ds)
 \Big|_E
^ {p  }\: dt
\Big\}
\\ &=:& C\left(\, S_0 ^ n+S_1 ^ n +S_2 ^ {n}\, \right).\phantom{\Big|}
\end{eqnarray*}
Since $ A ^ {\rho-\frac 1p} x_n\to A ^ {\rho-\frac1p}u_0$ in $E$, the Lebesgue DCT implies
\begin{eqnarray*}
 S_0 ^ n&\le&
 C_1 \int_0 ^ T \Vert A^{-(\rho-\frac1p)}e ^ {-tA}\Vert_{L(E,E)} |(A^{\rho-\frac1p}u_0-A^{\rho-\frac1p}x_n)|_E ^ p\: dt
 \longrightarrow 0, \end{eqnarray*}
as
 $n\to\infty$.\\
Since $\bar{u}_n=\mathcal{K}(x_n,\hat{\bar{u}}_n,
\tilde{\eta}_\ast)$, arguing as in Proposition \ref{unif-est-un}
we can show that $\bar{u}_n\in W^{\alpha,p}(0,T,E)$ for any
$\alpha\in (0,\frac 1p)$. Hence it follows from inequality
\eqref{aninequality} that
\begin{equation*}
\Lve \hat{\bar{u}}_n-u_\ast\Rve^p_{L^p(0,T;E)}\le C \Lve
\bar{u}_n-u_\ast\Rve^p_{L^p(0,T;E)}+C
2^{-np\alpha}\Lve\bar{u}_n\Rve^p_{W^{\alpha,p}(0,T;E)},\,\bar{\mathbb{P}}
\mbox{ a.s.\  }
\end{equation*}
from which and equality \eqref{limita-s} we infer that
\begin{eqnarray*}
 \lim_{n\to \infty }  \hat{\bar{u _n}}= u _\ast  \mbox{ in }   L ^ p(0,T;E), \quad \bar{\mathbb{P}} \mbox{ a.s.\  } .
 \end{eqnarray*}
Next,
by the  Young inequality we infer that
\begin{eqnarray*}
&&\hspace{-2truecm} \lefteqn{\int_0^T \vert \int_0 ^ t
  e ^ {-(t-s) A} \big[ F(s, \hat{\bar{u}}_n(s)) -F(s,   u_\ast(s))\big] ds  \vert^{p} dt} \\
 &\leq &  C
\int_0^T \vert A^{\rho-1}F(s, \hat{\bar{u}}_n(s))
-A^{\rho-1} F(s, u_\ast(s)) \vert^{p}\, ds,
\end{eqnarray*}
where $C= \left(\int_0^T \vert A^{1-\rho}e^{-sA}\vert\,
ds\right)$. By Assumption \eqref{assum-2} and since $\bar
u_n,u_\ast\in L ^ p(0,T;E)$ $\bar{\mathbb{P}}$ a.s.\ , there exists a constant $C>0$ such that
$$
 \int_0 ^ T \vert A^{\rho-1}F(s, \hat{\bar{u}}_n(s))\vert^{p}\, ds
 \le C,
$$
and
$$
 \int_0 ^ T \vert A^{\rho-1} F(s,u _\ast (s)) \vert^{p}\, ds
 \le C.
$$
The Lebesgue's dominated convergence theorem and the continuity of
$F$ (see Assumption \eqref{assum-2}) give
$$
S ^n_1 \to 0\; \mbox{ as }\;
n\to\infty.
$$
%%%%%%%%%%%%%%%%%%%%%%%%%%%%%%%%%%%%%%%%%%%%%%%%%%%%%%%%%%%%%%%%%%%%%%%%%%%%%%%%%%%%%%%%%%%%%%%%%%%%%%%%%%%%%%%%%%%%%
To  prove the convergence of   $S_2 ^ {n}$, we proceed in a similar
way as we have done to show convergence  for $ S_1 ^ n$.
In particular,  applying the Burkholder inequality and the Fubini
Theorem as well as \eqref{limita-s}, give
\begin{eqnarray}\label{ggob}\nonumber
\lefteqn{ S_2 ^ {n}  \le \bar {\mathbb{E}} \int_0^T\int_0 ^ t
 \int_Z }&&\\
\nonumber && \vert  e ^ {-(t-s) A} \big[ G(s,
\hat{\bar{u}}_n(s);z) -G(s,   u_\ast(s);z)\big] \vert^{p}
\nu(dz)\;ds  \; dt
 \\\nonumber
&\leq& C
 {\mathbb{E}} \int_0^T \int_Z
\\
&& \vert A^{\rho-\frac1p} G(s, \hat{\bar{u}}_n(s);z)
-A^{\rho-\frac1p} G(s,  u_\ast(s);z) \vert^{p}\,\nu(dz)\;
 ds,
\end{eqnarray}
where $C= \left(\int_0^T\vert A^{\frac1p-\rho}e^{-sA}\vert^p\,
ds\right)$.
By  Assumption \ref{assum-main-2} there exists a $C>0$ such that
\begin{eqnarray}\label{ggob2}
 &\sup_{n}&\int_Z \vert A^{\rho-\frac1p}G(s, \hat{\bar{u}}_n(s);z)\vert^{p}\,\nu(dz)
 \leq  C,
\end{eqnarray}
hence by Lebesgue Dominated Convergence Theorem we infer that
$$
S ^n_2\to 0\; \mbox{ as }\;
n\to\infty.
$$
\end{proof}
 To establish Theorem \ref{Th:bound} we need to check the following claim.
\begin{steps}
\addtocounter{lil}{1}
\item We have that for all $t\in [0,T]$ $\mathbb{P}$ a.s.\  $u_\ast (t) =
{\mathcal{K}}(u_0,u_\ast,\eta_\ast)(t) $.
\end{steps}
\begin{proof}
Let us fix $r\in (1,p)$. From Steps (II) to (VI) we infer that $\bar u_n\to u_\ast$ in
${\mathbf{L}}^r(\bar{\Omega}_T;E)$, % and
$$\bar u_n={\mathcal{K}} (x_n,\hat{\bar{u_n}},\bar \eta_n )\text{ in }{\mathbf{L}}^r(\bar{\Omega}_T;E)$$ and $${\mathcal{K}} (x_n,\hat{\bar{u_n}},\bar
\eta_n)-{\mathcal{K}}(u_0,u_\ast,\eta_\ast )\to 0
 \text{ in } {\mathbf{L}}^r(\bar{\Omega}_T;E).$$
 By the uniqueness of the limit, we infer that $u_\ast= {\mathcal{K}}(u_0,u_\ast,\eta_\ast)$ in
 ${\mathbf{L}}^r(\bar{\Omega}_T;E) $,
which implies that
  $\bar{\mathbb{P}}$-a.s.\
$u_\ast(t) = {\mathcal{K}}(u_0,u_\ast,\eta_\ast)(t)$ a.e. $t\in [0,T]$.
By equality \eqref{limita-s} we infer that  $\bar{\mathbb{P}}-a.s.$,  $u_\ast\in \mathbb{D}([0,T];D(A^{\rho^\prime-1}))$.
Hence  by combination of Lemmata \ref{est-stoc-conv-1} and
\ref{lem-est-stoc-conv-2},  and \cite[Theorem 2.8]{Di-Blasio-84}
 we deduce that $\bar{\mathbb{P}}-a.s.$, ${\mathcal{K}}(u_0,u_\ast,\eta_\ast)(\cdot)\in \mathbb{D}([0,T];D(A^{\rho^\prime-1})))$.
 Hence $\bar{\mathbb{P}}$-a.s. $u_\ast(t)={\mathcal{K}}(u_0,u_\ast,\eta_\ast)(t)$ for all $t\in [0,T]$.
\end{proof}
Since the solution to \eqref{eqn-3.1} is a fix point of the
operator ${\mathcal{K}}$, it follows from Step (VII) that $u_\ast$ is a
martingale solution to Equation  \eqref{eqn-3.1}.
The paths of the stochastic process $u_\ast$ are \cadlag in $D(A^{\rho^\prime-1})$ for any $\rho^\prime\in (0,\rho)$.
The assertion
\eqref{int-est-theo} follows from Step (II) and this ends the
proof of the first part of Theorem \ref{Th:bound}.
However, since $-A$ is the infinitesimal generator of a contraction type semigroup in $D(A^{\rho-\frac1p})$ and $u_\ast\in L^p(0,T;E)$ almost surely,
then thanks to the boundedness of $F$ and $G$ we easily infer from
Lemma \ref{lem-est-stoc-conv-2}-\eqref{lem-est-stoc-conv-2-dissipative} that the paths of $u_\ast$ are \cadlag on $D(A^{\rho-\frac1p})$
and this completes the proof of Theorem \ref{Th:bound}.
%%%%%%%%%%%%%%%%%%%%%%%%%%%%%%%%%%%%%%%%%%%%%%%%
\section{Proof of Theorem \ref{Th:general}}
In this section we replace  the boundedness
assumption
 on $F$ by the dissipativity of the drift $-A+F$.
The space $E$, $X$ are as before and we recall that
$E\subset
 X\subset D(A^{\rho-1})$.
Before we proceed let us  state the following
important consequence of Assumption \ref{assum-main-3}.
\begin{lemma}\label{Lem:a'priori}(See Da Prato \cite{MR0500309})
Assume that  $X$ is a Banach  space, $-A$ a generator  of a
strongly continuous semigroup of bounded linear  operators on  $X$
and a   mapping $F:[0,T]\times X\to X$ satisfies Assumption
\ref{assum-main-3}. Assume  that for $\tau\in[0,\infty]$ two
continuous functions $z,v:[0,\tau )\to X$ satisfy
\begin{eqnarray*}
z(t)&=&\int^t_0e^{-(t-s)A}F(s,z(s)+v(s))\,ds,\;t<\tau .
\end{eqnarray*}
Then
\begin{eqnarray}
|z(t)|_X&\le&\int^t_0e^{-k(t-s)}a(|v(s)|_X)\,ds,\;\;0\le t< \tau .
\label{eq:a'priori}
\end{eqnarray}
\end{lemma}
In this section we show the existence of a martingale solution
to equation \ref{eqn-3.1}. First, let us notice   that Assumption
\ref{assum-main-3} implies that
\begin{equation}
|F(t,y)|_X\le a(|y|_X), \;\; \;t\ge 0, \; y \in X. \label{dar-2.3}
\end{equation}
\begin{proof}[Proof of Theorem \ref{Th:general}]  Without loss of generality we
assume that $k=0$.
Let  $(F_n)_{n\in{\mathbb{N}}}$ be %, :[0,\infty)\times X\to X$ ba
a sequence of functions from $[0,\infty)\times X$ to $X$ given by
Assumption \ref{assum-main-3}-\eqref{iii}. In particular, there
exists a sequence $(R_F^n)_{n\in{\mathbb{N}}}$ of positive numbers, such
that $|F_n(s,y)|_X \le R_F^n$ for all $(s,x)\in[0,\infty)\times
X$, $n\in{\mathbb{N}}$, and $|F_n(s,x)-F(s,x)|_X \to 0$ as $n\to\infty$
for all $(s,x)\in[0,\infty)\times X$. % converges to $F(s,x)$.
For $n\in{\mathbb{N}}$ let $F_n ^ E$ be the restriction of $F_n$ to
$[0,\infty)\times E$. Since $E\hookrightarrow X$ continuously,
$F_n ^ E:[0,\infty)\times E\to X$ is also separately continuous
with respect to both variables and bounded by $R_F^n$.  Finally, by the conitnuity of the embeddings $X\hookrightarrow
D(A^{\rho-\frac1p})\subset D(A^{\rho-1})$ the function $F ^ E_n$
satisfies all the assumptions of Theorem \ref{Th:bound}. It is clear that the restriction of $G$ to $[0,T]\times E$ also satisfies the assumptions of Theorem \ref{Th:bound}.
In what follows we still denote by $G$ this restiction.
In  view of that theorem
there exists an $E$-valued martingale solution to the following Problem
\begin{equation}\label{dar-2.66}
 \begin{cases}
  d {u}_n(t)=&[-A {u}_n(t)+F_n(s,
{u}_n(t))]\,dt+\int_Z G(s, {u}_
n(t);z)\, \tilde \eta_n (dz;dt),\\\
 {u}_n(0)=&u_0.
 \end{cases}
\end{equation}
%\end{eqnarray}
%
Let us denote this martingale solution by
\begin{eqnarray*}&&
\left(\Omega_n,{{\mathcal{F}}}_n,\mathbb{P}_n,{\mathbb{F}}_n,
\{\eta_n(t,z)\}_{t\ge 0,z\in Z},\{ {u}_n(t)\}_{t\ge 0}\right).
\end{eqnarray*}
We donote by $\En$ the mathematical expectation on
$\left(\Omega_n,{{\mathcal{F}}}_n,\mathbb{P}_n\right)$.

In view of Theorem \ref{Th:bound}, for each $n\in{\mathbb{N}}$, $u_n$ is a mild
solution \cadlag on $D(A^{\rho-\frac 1p})$ over the probability space
$(\Omega_n,{{\mathcal{F}}}_n,{\mathbb{F}}_n , \mathbb{P}_n) $.
In particular $$ {u}_n(t) = e^{-tA}u_0+z_n(t)+v_n(t), \, \, t\in [0,T],$$
where

\begin{eqnarray}
v_n(t)&=& \int^t_0\int_Z e^{-(t
-s)A}G(s, {u}_n(s);z)\,\tilde \eta_n (dz;ds),\\
z_n(t)&=&\int_0^te^{-(t-s)A}F_n(s, {u}_n(s)))\,ds.
\label{dar-2.61}
\end{eqnarray}
Notice, that $z_n(t)= {u}_n(t)-v_n(t)-e^{-tA}u_0$, $t\in [0,T]$. %%%%%%%%%%%%%%%%%%%
%%%%%%%%%%%%%%%%%%%%%%%%%%%
%%%%%%%%%%%%%%%%%%%%%%%%%%%%%%%%%%%%%%%%%%%%%%%%%%%%%%%
%
%
Similarly to the proof of Theorem \ref{Th:bound} the proof of
Theorem \ref{Th:general} will be divided into several steps. The
first two steps are the following.
\usecounter{lil}
\begin{steps}
\setcounter{lil}{1}
\item
\addtocounter{lil}{1} Let $q_{\max}$ be defined by \eqref{qu-max}. Then for any $\tilde{q}\in (q, q_{\max})$  and $r\in (1,p)$, we have
\begin{equation}\label{EST-VN-STEP1}
\sup_{n\in{\mathbb{N}}}\En \left\| v_n \right\|^r_{L ^ {
\tilde{q}}(0,T;E)}<\infty.
\end{equation}
\end{steps}
\begin{proof}
\textbf{Step (I)} follows from Lemma \ref{lem-tight-st0-conv-3} and Lemma \ref{lem-tight-st0-conv-3b}.
\end{proof}
\begin{steps}
\item \addtocounter{lil}{1} For any $\tilde{q}\in (q, q_{\max})$ and $r\in (1,p)$ defined in \textbf{Step(I)}, we have
 \begin{equation}\label{EST-ZN-STEP2}
  \sup_{n\in{\mathbb{N}}} \En \sup_{0\le t\le T} |  z_n(t)|_X^{\frac{r}{
 \tilde{q}}}<\infty.
 \end{equation}

  Moreover, the laws of the family  $\{  z_n:n\in{\mathbb{N}}\}$ are tight on ${\mathcal{C}}([0,T];X)$.
\end{steps}
\begin{proof}
\addtocounter{zaehler}{1}
 By Lemma \ref{Lem:a'priori}
we infer that for any $T\ge 0$ %\in0,T$
\begin{eqnarray} \label{dar-2.66-z} \nonumber \sup_{0\le t\le T}
|z_n(t)|_X\le \int_0 ^ T e ^ {-k(t-s)} a\left( |v_n(s)|_X+ |e ^
{-sA} x|_X\right) \: ds
\\
\le C\: \int_0 ^ T \left( 1+ |v_n(s)|_X ^ q + |e ^ {-sA} x|_X ^
q\right)\, ds. \end{eqnarray}
Since the embedding $E\subset X$ is continuous,
 it follows from Step (I) that
there exists $\tilde{q}>q$ such that for any $r\in (1,p)$
\begin{eqnarray}\label{inte-est} %$$
\sup_n \En |v_n| ^ r _{L ^ {\tilde{q}}(0,T;X)} <\infty. \end{eqnarray}
Therefore,
$$\sup_n \En  \sup_{0\le t\le T} |z_n(t)|_X^ \frac{ r}{\tilde{q}}<\infty.$$
Hence, we proved the first part of \textbf{Step (II)}. Note that the last inequality
implies that $\sup_n \En  |z_n|_{L ^ {\tilde{q}}(0,T;X)} ^
\frac{ r}{\tilde{q}} <\infty$.

 Before we proceed further, we recall that there exist $\theta<1-\frac{q}{q_{\max}}$ and an UMD, type $p$ and separable Banach space $Y$ such that $D(A_Y^\theta) \subset
 X\subset Y$. To prove the second part we use the identity \footnote{Here $\Lambda^{-1}$ is defined through $A_Y$.} $z_n=\Lambda^{-1}
F_n(s, u_n(s))$  and Remark \ref{rem-M type-p} along with Corollary
\ref{C:comp}. But first we need to show that $ \sup_{n\in \mathbb{N}}\lvert F_n(\cdot,
u_n(\cdot))\rvert_{L^{\tilde{p}}(0,T; Y)}$ is bounded in probabilty for some $\tilde{p}\in (1, \frac{q_{\max}}{q})$.
Let us fix  $\tilde{p}\in (1, \frac{q_{\max}}{q})$. From Lemma \ref{Lem:a'priori} and the continuous embedding $E\subset X$ we infer that
$$ \lvert F_n(s,u_n(s))\rvert_X^{\tilde{p}}\le C(1+|e^{-s
A}x|^{\tilde{p}}+|v_n(s)|^{\tilde{p}q}+|z_n(s)|^{\tilde{p}q}).
$$   We easily obtain from this inequality that
\begin{equation}\label{est-Fn}
 \lvert F_n(\cdot, u_n(\cdot)) \rvert^{\tilde{p}}_{L^{\tilde{p}}(0,T;X)}\le C(1+
+\lvert v_n\rvert^{\tilde{p}q}_{L^{\tilde q}(0,T;X) }+ \lvert
z_n\rvert^{\tilde{p}q}_{L^{\infty}(0,T;X) }).
\end{equation}
Taking $\tilde{q}=\tilde{p}q\in (q,q_{\max})$ and rasing to the power $\frac{r}{\tilde{q}^2}$ both sides of  \eqref{est-Fn} implies that
\begin{equation*}%\label{est-Fn-2}
 \lvert F_n(\cdot, u_n(\cdot)) \rvert^{\frac{r\tilde{p}}{\tilde{q}^2}}_{L^{\tilde{p}}(0,T;X)}\le C(1+
+\lvert v_n\rvert^{\frac r{\tilde{q}}}_{L^{\tilde q}(0,T;X) }+ \lvert
z_n\rvert^{\frac r{\tilde{q}}}_{L^{\infty}(0,T;X) }).
\end{equation*}
Since $Y\subset X$ is continuous we obtain that
\begin{equation}\label{est-Fn-2}
 \lvert F_n(\cdot, u_n(\cdot)) \rvert^{\frac{r\tilde{p}}{\tilde{q}^2}}_{L^{\tilde{p}}(0,T;Y)}\le C(1+
+\lvert v_n\rvert^{\frac r{\tilde{q}}}_{L^{\tilde q}(0,T;X) }+ \lvert
z_n\rvert^{\frac r{\tilde{q}}}_{L^{\infty}(0,T;X) }).
\end{equation}
By Chebychev inequality we derive that for any $m\in \mathbb{N}$
\begin{equation*}
\begin{split}
 \mathbb{P}_n\biggl(\lvert F_n(s, u_n(s))
\rvert_{L^{\tilde{p}}(0,T;Y)}\ge  m \biggr)\le
\frac{1}{m^{\frac{r}{q\tilde{q}} }}\biggl(C +\En \lvert
v_n(s)\rvert^\frac{ r}{\tilde{ q}}_{L^{\tilde q}(0,T;X) }\\
+\En\lvert z_n(s)\rvert^\frac{ r}{ \tilde{q}}_{L^{\tilde
q}(0,T;X) } \biggr).
\end{split}
\end{equation*}
Since $\En \lvert v_n(s)\rvert^\frac {r}{\tilde{ q}}_{L^{\tilde
q}(0,T;X) }$ and $\En\lvert z_n(s)\rvert^\frac{ r}{\tilde{
q}}_{L^{\tilde q}(0,T;X) }$ are uniformly bounded w.r.t $n$, we
derive that $ \sup_{n\in \mathbb{N}}\lvert F_n(\cdot,
u_n(\cdot))\rvert_{L^{\tilde{p}}(0,T; Y)}$ is bounded in probabilty for any $\tilde{p}\in (1, \frac{q_{\max}}{q})$.

Now we choose $\tilde{p}\in (1,\frac{q_{\max}}{q})$ such that $\theta\in (0, 1-\frac{1}{\tilde{p}})$. Since $\theta <1-\frac{1}{\tilde{p}}$ and $D(A_Y^\theta) \subset
 X$, we can infer from Corollary
\ref{C:comp} that the family of laws of $\{z_n=\Lambda^{-1} F(\cdot, u_n(\cdot)); n\in \mathbb{N}\}$ is tight on $
{\mathcal{C}}([0,T]; D(A_Y^\theta))$ and hence  on ${\mathcal{C}}([0,T]; X)$.
%\addtocounter{zaehler}{1}
\end{proof}
\begin{remark}\label{rem-stoc-conv}
Let $q$ be a number in the interval $[p,\infty)$. It follows from Step (I) and Step (II) that for any  $\tilde{q} \in (q, q_{\max})$ and $r\in (1,p)$
\begin{equation*}
\sup_{n\ge 1}\En \Lve u_n\Rve^{\frac{r}{\tilde{q}}}_{L^{\tilde
q}(0,T;X)}\le C.
\end{equation*}

For $q<p$ the above inequality holds with $\tilde{q}=p$.
\end{remark}
\begin{remark}
 In \cite{Brz+Gat_1999} the first named author and Gatarek constructed an approximation of $F$ as follows.
 Let
$(F_n)_{n\in{\mathbb{N}}}$ be defined by
\begin{equation*}
F_n(s,x)=\begin{cases} F(s,x) \text{ if } \lvert x\rvert_X\le n,\\
F(s, \frac{n}{\lve x\rve_X}x)\text{ otherwise.}
\end{cases}
\end{equation*}
By \eqref{dar-2.3} $\lve F_n(s,y)\rve\le a(n)$, for all $s\ge0$,
$y\in E$. They solved the Problem \eqref{dar-2.66} driven by
Wiener noise on the random interval $[0,\tau_n\wedge T]$ where the
sequence of stopping times $\{\tau_n: n\ge 1\}$ is defined by
$$\tau_n=\inf\{t\in [0,T]: \lve u_n(t)\rve_X\ge n\}.$$ By proving that
$\sup_{t\in [0,T]}\lve u_n(t)\rve_X$ is uniformly bounded, which
implies that $\tau_n\uparrow T$ alsmost surely as $n\rightarrow
\infty$, and then using \eqref{dar-2.3} and Corollary \ref{C:comp}
they could show that the laws of $z_n$ is tight on $C([0,T];X)$. In
our framework we know a priori that $u_n$ is only \cadlag in $D(A^{\rho^\prime-1})$ for $\rho^\prime\in (0,\rho)$, hence $\tau_n$ will not be a well defined stopping time and we
will not be able to show that $\sup_{t\in [0,T]}\lve u_n(t)\rve_X$
is uniformly bounded.
\end{remark}
In order
to use the Skorokhod embedding theorem, we also need the following.
\begin{steps}
\item For $\rho^\prime\in (0,\rho)$ let $B_0=D(A^{\rho^\prime-1})$. The laws of the family  $\{  v_n:n\in{\mathbb{N}}\}$ are tight on
$L^{\tilde q}(0,T,E)\cap \mathbb{D}([0,T];B_0)$ and those of $\{ \eta_n:n\in{\mathbb{N}}\}$ are tight on
$M_I(Z\times [0,T])$. \addtocounter{lil}{1}
\end{steps}
\begin{proof}
By Lemma \ref{tight-sto-conv-1} and Lemma \ref{tight-sto-conv-2}
the laws of the family $\{ v_n:  n\in{\mathbb{N}}\}$ are tight on
$L^{\tilde q}(0,T,E)\cap  \mathbb{D}([0,T], B_0)$. Since ${\mathbb{E}} _
n\eta_n(A,I) = {\mathbb{E}} \eta(A,I)$ for all $A\in{\mathcal{B}}(Z)$ and
$I\in{\mathcal{B}}([0,T])$, the laws of the random variables
$\{\eta_n:n\in{\mathbb{N}}\}$ are identical on $M_{\mathbb{N}}(Z\times [0,T])$.
Hence we can deduce from \cite[Theorem 3.2]{para} that the laws of
the family $\{\eta_n:n\in{\mathbb{N}}\}$ are tight on $M_{\mathbb{N}}(Z\times
[0,T])$.

\end{proof}
Let $B_0$ as in Step (III). From Steps (I), (II), (III) and Prokhorov's theorem it follows that there exists a
subsequence of $\{ (z_n,v_n,\eta_n)  : n\in{\mathbb{N}}\}$, also denoted by
$\{ ( z_n,v_n,\eta_n):   n\in{\mathbb{N}}\}$ and a $ {\mathcal{C}}([0,T], X)\times
L^{\tilde q}(0,T,E)\cap \mathbb{D}([0,T], B_0)\times M_{\mathbb{N}}(Z\times
[0,T])$ -- valued random variable $(z_\ast, v _\ast,\eta_\ast)$,
such that the sequence of laws of $ (z_n,v_n,\eta_n)$ converges to
the law of $(z_\ast, v _\ast,\eta_\ast)$. Moreover, by Theorem
\ref{thm-Skorokhod}, there exists a probability space $(\hat
\Omega,\hat {\mathcal{F}},\hat{\mathbb{P}})$ and $ {\mathcal{C}}([0,T], X)\times L^{\tilde
q}(0,T,E)\cap \mathbb{D}([0,T], B_0)\times  M_{\mathbb{N}}(Z\times [0,T])$- valued
random variables $( z_\ast, v_\ast,\eta_\ast)$, $(\hat z_n,\hat
v_n,\hat \eta_n)$, $n\in{\mathbb{N}}$, such that $ \hat{\mathbb{P}}$-a.s.
\begin{eqnarray}\label{conv-aa}
(\hat z_n,\hat v_n,\hat \eta_n)\to(z_\ast, v_\ast,\eta_\ast)
\end{eqnarray}
on $ {\mathcal{C}}([0,T ]; X)\times L^{\tilde{q}}(0,T,E)\cap \mathbb{D}([0,T], B_0)\times  M_{\mathbb{N}}(Z\times [0,T])$,   $\hat \eta_n=\eta_\ast$,
and
$${\mathcal{L}}( (\hat z_n,\hat v_n,\hat \eta_n))={\mathcal{L}}( ( z_n,v_n,\eta_n)),$$ for all $n\in{\mathbb{N}}$.
{We define a filtration $\hat {\mathbb{F}} = (\hat {\mathcal{F}}_t)_{t\in [0,T]}$ on
$(\hat \Omega,\hat{\mathcal{F}})$ as  the one generated by $\eta_\ast$,
$z_\ast$, $v_\ast$ and the families  $\{z_n: n\in{\mathbb{N}}\}$ and
$\{v_n: n\in{\mathbb{N}}\}$.
}%
\\[0.1cm]
The next two steps imply that the following two It\^o integrals
over
the filtered probability space $(\hat \Omega,\hat {\mathcal{F}},\hat {\mathbb{F}} ,
\hat{\mathbb{P}})$
$$
\int_0 ^ t  \int_Z e ^ {-(t-s)A}\: G(s,\hat v_n(s)+\hat z_n(s)
,z))\, \tilde{\hat \eta}_n(dz,ds),\quad t\in [0,T],
$$
and
$$
\int_0 ^ t  \int_Z e ^ {-(t-s)A}\: G( s,v_\ast(s)+z_\ast (s),z))\,
\tilde{ \eta}_\ast(dz,ds),\quad t\in [0,T],
$$
are well defined.
\begin{steps}
\item The following  holds \addtocounter{lil}{1}
\begin{steps-2}
\item for all $n\in{\mathbb{N}}$, $\hat  \eta_n$ is a time homogeneous
Poisson random measure on ${\mathcal{B}}(Z)\times {\mathcal{B}}([0,T])$
over $(\hat \Omega,\hat {\mathcal{F}},\hat {\mathbb{F}},
\hat{\mathbb{P}})$ with intensity measure $\nu$; \item $\eta_\ast$ is a
time homogeneous Poisson random measure on ${\mathcal{B}}(Z)\times {\mathcal{B}}(0,T)$
over $(\hat \Omega,\hat {\mathcal{F}},\hat {\mathbb{F}}_t,
\hat{\mathbb{P}})$ with intensity measure $\nu$;
\end{steps-2}
\item The following  holds \addtocounter{lil}{1}
 \begin{steps-2} \item  for all $n\in{\mathbb{N}}$, the processes $\hat  v_n$ and $\hat z_n$ are,
 with respect
 to $(\hat  {\mathcal{F}}_t)_{t\ge 0}$, progressively  measurable;
\item the processes $ z_\ast$ and $v_\ast$ are,  with respect to
$(\hat {\mathcal{F}}_t)_{t\ge 0}$, progressively measurable.
\end{steps-2}  \end{steps}
The proofs of Step (IV) and (V) are the same as the proofs of Step
(IV) and (V) of Theorem \ref{Th:bound}.
Also as earlier in the proof of Theorem \ref{Th:bound} in order to
complete the proof of Theorem \ref{Th:general} we have to prove the following claim.

\begin{steps}
\item \addtocounter{lil}{1} Let  $u_\ast=e^{-\cdot A}u_0+
z_\ast+v_\ast$ and ${\mathcal{K}}$ be the mapping defined by\begin{eqnarray*} {\mathcal{K}}(u_0,u,
\eta)(t) &=&e ^ {-tA} u_0 +\int_0 ^ te ^ {-(t-s)A} F(s, {u}(s))\,ds
\\ && {}+ \int_0 ^ t \int_Z e ^ {-(t-s)A} G(s, {u}
(s);z)\, \tilde \eta (dz;ds), \nonumber \end{eqnarray*} for $t\in [0,T]$,
$u\in \mathbb{D}([0,T],B)$ and $\eta\in M_{\mathbb{N}}(Z\times [0,T])$. Then  $\hat{\mathbb{P}}$-a.s. \begin{eqnarray}\label{num-aa} u_\ast
(t)={\mathcal{K}}(u_0,u_\ast,\eta_\ast)(t), \quad \forall t\in[0,T]. \end{eqnarray}
\end{steps}
\begin{proof}
Since the process $u_n=e^{-\cdot A}u_0+z_n+v_n$ is a martingale
solution of problem \eqref{dar-2.66} and
$${\mathcal{L}}( (\hat z_n,\hat v_n,\hat \eta_n))={\mathcal{L}}( ( z_n,v_n,\eta_n)),$$ for all
$n\in{\mathbb{N}}$, we can argue as in Step (V) of the previous section and
prove that $\hat{\mathbb{P}}$-a.s.
\begin{equation}\label{huneqkn}
\hat{u}_n(t)=\mathcal{K}_n(u_0,\hat{u}_n, \hat{\eta}_n)(t),
\end{equation}
 for all
$t\in [0,T]$. Here $\mathcal{K}_n$ is obtained by replacing $F$
with $F_n$ in the definition of $\mathcal{K}$.
 Since the laws of $z_n$  and $\hat{z}_n$
are equal on $C([0,T];X)$, by \eqref{dar-2.66-z} and \eqref{EST-ZN-STEP2} we infer that
\begin{equation}\label{est-hzn}
\sup_{n\ge 1}\En \Lve \hat{z}_n\Rve^{\frac{r}{
\tilde{q}}}_{C([0,T];X)}<\infty.
\end{equation}
Since $v_n$ and $\hat{v}_n$ have equal laws on
$L^{\tilde{q}}(0,T;E)$, from \eqref{EST-VN-STEP1} we
deduce that
\begin{equation}\label{est-hvn}
\sup_{n\ge 1}\En \Lve
\hat{v}_n\Rve^{r}_{L^{\tilde{q}}(0,T;E)}<\infty.
\end{equation}

Invoking \eqref{conv-aa} we infer that, $\hat{\mathbb{P}}$-a.s., as $n\rightarrow \infty$,
\begin{equation*}
\Lve \hat{z}_n -z_\ast\Rve^{\frac{r}{2\tilde{q}}}_{C([0,T];X)}\rightarrow 0\;\; \mbox{ and } \;\;\Lve \hat{z}_n\Rve^{\frac{r}{2\tilde{q}}}_{C([0,T];X)} \rightarrow \Lve
z_\ast\Rve^{\frac{r}{2\tilde{q}}}_{C([0,T];X)}.
\end{equation*}
 Thanks to \eqref{est-hzn} the sequence  $\Lve
\hat{z}_n\Rve^{\frac{r}{2\tilde{q}}}_{C([0,T];X)}$ is $\hat{\mathbb{P}}$-uniformly integrable. Thus the applicability of the Lebesgue
Dominated Convergence Theorem implies that
\begin{equation*}\lim_{n\rightarrow \infty}
\En\Lve \hat{z}_n\Rve^{\frac{r}{2\tilde{q}}}_{C([0,T];X)} = \En\Lve
z_\ast\Rve^{\frac{r}{2\tilde{q}}}_{C([0,T];X)}.
\end{equation*}
 Thanks to \eqref{est-hzn} and this last
convergence  we can prove, by a similar argument used as above, that
\begin{equation}\label{conv-hzn}
\lim_{n\rightarrow \infty}  \En\Lve \hat{z}_n -z_\ast\Rve^{\frac{r}{2\tilde{q}}}_{C([0,T];X)}=
0.
\end{equation}
With similar argument we can also show that
\begin{equation}\label{conv-hvn-1}
\lim_{n\rightarrow \infty} \En \Lve \hat{v}_n
-v_\ast\Rve^{\frac{r}{2 \tilde{q}}}_{L^{\tilde{q}}(0,T;B_0)}= 0,
\end{equation}
and
\begin{equation}\label{conv-hvn-2}
\En \Lve \hat{v}_n
-v_\ast\Rve^{\frac{r}{2 \tilde{q}}}_{L^{\tilde{q}}(0,T;E)}\rightarrow 0,
\end{equation}
We derive from \eqref{conv-hzn}, \eqref{conv-hvn-2} and
\eqref{conv-hvn-2} that
\begin{equation*}
\lim_{n\rightarrow \infty}  \En \Lve \hat{u}_n
-u_\ast\Rve^{\frac{r}{2\tilde{q}}}_{L^{\tilde{q}}(0,T;B_0)}= 0.
\end{equation*}
and
\begin{equation*}
\lim_{n\rightarrow \infty}   \En \Lve \hat{u}_n
-u_\ast\Rve^{\frac{r}{2\tilde{q}}}_{L^{\tilde{q}}(0,T;X)}= 0,
\end{equation*}
which is correct because $E\subset X$.
In the other hand, by Proposition \ref{zero} we have
\begin{equation*}
\begin{split}
& \En \Lve
\mathcal{K}_n(x,\hat{u}_n,\hat{\eta}_n)-\mathcal{K}(x,u_\ast,\eta_\ast)\Rve^{\tilde
r}_{L^{p^\ast}(0,T;B_0)}\\ &= \En \Lve
\mathcal{K}_n(x,\hat{u}_n,\eta_\ast)-\mathcal{K}(x,u_\ast,\eta_\ast)\Rve^{\tilde
r}_{L^{p^\ast}(0,T;B_0)}\\
&\quad  \le C \En \biggl\Lve \int_0^t
e^{-(t-s)A}\left(F_n(s,\hat{u}_n(s))-F(s,u_\ast(s))\right)ds\biggr\Rve^{\tilde r}_{L^{p^\ast}(0,T;X)}\\
&\quad  + C \biggl\Lve \int_0^t \int_{\mathcal{Z}}
e^{-(t-s)A}\left(G(s,\hat{u}_n(s))-G(s,u_\ast(s))\right)\tilde{\eta}_\ast(dz,ds)\biggr\Rve^{\tilde{r}}_{\mathcal{M}^{p}(0,T;E)}\\
& \le \quad I^n_1+I^n_2,
\end{split}
\end{equation*}
where $\tilde{r}=\frac{r}{2\tilde{q}q}$ and
$p^\ast=\min(\frac{\tilde q}{q},p)$. Arguing as in Step (VI) of
the previous section we can show that $I^n_2\rightarrow 0$ as
$n\rightarrow 0$. To deal with $I^n_1$ we first use Assumption
\ref{assum-main-3} and \eqref{conv-aa} to derive that
\begin{align*}
\Lve F_n(\cdot,\hat{u}_n)-F(\cdot, u_\ast)\Rve^{\tilde
r}_{L^{p^\ast}(0,T;X)}\rightarrow 0,
\end{align*}
as $n\rightarrow \infty$. Since, by \eqref{est-Fn},
\eqref{est-hzn} and \eqref{est-hvn}, $$\sup_{n\ge 1}\En \Lve
F_n(\cdot, \hat{u}_n)\Rve^{2
\tilde{r}}_{L^{p^\ast}(0,T;X)}<\infty,$$ we can apply the
Lebesgue Dominated Convergence Theorem and derive that
$I^n_1\rightarrow 0$ as $n\rightarrow \infty$. Therefore
\begin{align*}
\En \Lve \hat{u}_n
-u_\ast\Rve^{\tilde{r}}_{L^{p^\ast}(0,T;B_0)}\rightarrow 0,\\
\En \Lve
\mathcal{K}_n(x,\hat{u}_n,\hat{\eta}_n)-\mathcal{K}(x,u_\ast,\eta_\ast)\Rve^{\tilde
r}_{L^{p^\ast}(0,T;B_0)}\rightarrow 0,
\end{align*}
as $n\rightarrow \infty$. These two facts along with
\eqref{huneqkn} implies that $\hat{\mathbb{P}}$-a.s. and for a.e. $t\in
[0,T]$
\begin{equation*}
u_\ast(t)=\mathcal{K}(x,u_\ast,\eta_\ast)(t).
\end{equation*}
Since $u_\ast$ and $\mathcal{K}(x,u_\ast,\eta_\ast)$ are \cadlag
functions taking values in $B_0$, the last equation holds
$\hat{\mathbb{P}}$-a.s. and for all $t\in [0,T]$.
\end{proof}
The paths of the stochastic process $u_\ast$ are \cadlag in $D(A^{\rho^\prime-1})$ for any $\rho^\prime\in (0,\rho)$.
However, since $-A$ is the infinitesimal generator of a contraction type $C_0$-semigroup on $D(A^{\rho-\frac1p})$ and $u_\ast \in L^{\tilde{q}}(0,T;X)$ almost surely,
then we easily infer from
Lemma \ref{lem-est-stoc-conv-2}-\eqref{lem-est-stoc-conv-2-dissipative} that the paths of $u_\ast$ are cadlag in $D(A^{\rho-\frac1p})$.
Similar calculations as done in Step (I) and Step (II) (see also
Remark \ref{rem-stoc-conv}) yield that for any $\tilde{q}\in (q, q_{\max})$ and $r\in (1,p)$
\begin{equation*}
\mathbb{E}\Lve u_\ast\Rve^{\frac{r}{\tilde{q}}}_{L^{\tilde{q}}(0,T;X)}<\infty.
\end{equation*}
This completes the proof of Theorem \ref{Th:general}

%%%%%%%%%%%%%%%%%%%%
%%%%%%%%%%%%%%%%%%%%%%%%%
%
\end{proof}

\usecounter{zaehler}
\setcounter{zaehler}{0} \addtocounter{zaehler}{1}

\appendix
%%
%%%%%%%%%%%%%%%%%%%%%%%%%

\section{Stochastic Preliminaries}\label{sec_stoch_prel}
In this paper we use the setting of Poisson random measures to deal
with equations of type \eqref{spde01}. Therefore,  in
the first part of this section we introduce this concept and in the
second part we point out the relationship between Poisson random
measures and L\'evy processes.
{Let   $1<p\le 2$ be fixed throughout the whole paper. Moreover,
throughout the whole paper, we assume that $(\Omega,{\mathcal{F}},{\mathbb{F}},\mathbb{P})$
is a complete filtered probability space. Here, for simplicity, we
denoted the filtration $\{{\mathcal{F}}_t\}_{t\ge 0}$ by ${\mathbb{F}}$. }
\subsection{Poisson random measures}
We start with the following definition.
\begin {definition}
(See Pisier \cite{0606.60008})
%
%Let $1\le p\le 2$.
A Banach space $E$ is of martingale type $p$, iff there exists a
constant $C=C(E,p)$ such that for any $E$-valued discrete
martingale $(M_0,M_1,M_2, \ldots)$
with $M_0=0$ the following inequality holds
$$
\sup_{n\ge 1} {\mathbb{E}}|M_n| ^p \le C\: \sum_{n\ge 1} {\mathbb{E}}|M_n-M_{n-1}|
^p.
$$
\end{definition}

The following definitions are presented here for the sake of
completeness because the notion of time homogeneous random measure
is introduced in many, not always equivalent ways.
\begin{definition}\label{def-Prm}(See \cite{Ikeda+Watanabe_1981}, Definition I.8.1)
Let us assume that  $(Z,{\mathcal{Z}})$  is a measurable space and $\nu\in {\mathcal{M}}_+(Z)$.
\\  A measurable function
\[\eta:
(\Omega,{\mathcal{F}})\to (M_I(Z\times {\mathbb{R}}_+),{\mathcal{M}}_I(Z\times {\mathbb{R}}_+)) \]
is called a
{\sl time homogenous Poisson random measure} $\eta$ on $(Z,{\mathcal{Z}})$
over $(\Omega,{\mathcal{F}},{\mathbb{F}},\mathbb{P})$  iff the following conditions are satisfied
\begin{enumerate}[(i)]
\item \label{prm-i} for each $B\in  {\mathcal{Z}} \otimes
\mathcal{B}({\mathbb{R}_+}) $,
 $\eta(B):=i_B\circ \eta : \Omega\to \bar{\mathbb{N}} $ is a Poisson random variable with parameter ${\mathbb{E}} \eta(B)$.\footnote{If ${\mathbb{E}} \eta(B)=\infty$, then $\eta(B)=\infty$}
\item \label{prm-ii} $\eta$ is independently scattered, i.e. if the sets $
B_j \in   {\mathcal{Z}}\otimes \mathcal{B}({\mathbb{R}_+})$, $j=1,\cdots,
n$, are  pair-wise disjoint, then the random variables $\eta(B_j)$,
$j=1,\cdots,n $, are pair-wise independent;
\item \label{prm-3}
for all $B\in  {\mathcal{Z}} $ and $I\in \mathcal{B}({\mathbb{R}_+})$, $\mathbb{E}\big[\eta (B\times I)\big]=\Leb(I)\nu(B)$;
\item \label{prm-4} for each $U\in {\mathcal{Z}}$, the $\bar{\mathbb{N}}$-valued
process $(N(t,U))_{t\ge 0}$  defined by
$$N(t,U):= \eta(U \times (0,t]), \;\; t\ge 0$$
is ${\mathbb{F}}$-adapted and its increments are independent of the past,
i.e.\ if $t>s\geq 0$, then $N(t,U)-N(s,U)=\eta(U \times (s,t])$ is
independent of $\mathcal{F}_s$.
\end{enumerate}

\end{definition}

Given a complete filtered  probability space $(\Omega,{\mathcal{F}},{\mathbb{F}}
,\mathbb{P})$, where the symbol ${\mathbb{F}}=\{\mathcal{F}_t\}_{t\geq 0}$ denotes the
filtration, the predictable random field ${\mathcal{P}} $ on $\Omega\times
{\mathbb{R}} _ +$ is the $\sigma$--field generated by all continuous
${\mathbb{F}}$--adapted  processes (see e.g.\ Kallenberg \cite[Chapter
25]{MR1876169}).
A real valued stochastic process $\{x(t):t\geq 0\}$,  defined on a
filtered probability space $(\Omega,{\mathcal{F}},{\mathbb{F}},\mathbb{P})$ is  called
predictable, if the mapping  $x:\Omega\times\Rb{+}\to\Rb{}$ is
${\mathcal{P}}/\mathcal{B}({\mathbb{R}})$-measurable. A random measure $\gamma$ on
${\mathcal{Z}}\times {\mathcal{B}}({\mathbb{R}} _ +)$ over $(\Omega;{\mathcal{F}},{\mathbb{F}},\mathbb{P})$ is called
predictable, iff for each $U\in {\mathcal{Z}}$, the ${\mathbb{R}}$--valued process
${\mathbb{R}}_+\ni t\mapsto \gamma( U\times(0,t])$ is predictable.

\begin{definition}\label{def-imPrm}
Assume that  $(Z,{\mathcal{Z}})$ is a measurable space and $\nu$ is a
non-negative measure on $(Z,{\mathcal{Z}})$. Assume that $\eta$ is a
  {time homogeneous} Poisson
random measure  {with  intensity measure $\nu$ on  $(Z,{\mathcal{Z}})$  }
over $(\Omega,{\mathcal{F}},{\mathbb{F}},\mathbb{P})$. The {\em compensator } of  $\eta$ is
the unique {predictable random measure}, denoted by $\gamma$, on $
{\mathcal{Z}} \times {\mathcal{B}}({\mathbb{R}}_+ ) $ over {$(\Omega,{\mathcal{F}},{\mathbb{F}},\mathbb{P})$} such that
for each $T<\infty$ and $A\in {\mathcal{Z}} $ with ${\mathbb{E}}\eta( A\times
(0,T])<\infty$, the $\mathbb{R}$-valued processes
$\{\tilde{N}(t,A)\}_{t\in (0,T]}$  defined by
$$
\tilde{N}(t,A):= \eta(  A\times (0,t] )-\gamma( A\times (0,t] ),
\quad   0< t\le T,
%t \in \mathbb{R}_+
$$
is a martingale on $(\Omega,{\mathcal{F}},{\mathbb{F}},\mathbb{P})$.
\end{definition}

\begin{remark}
Assume that $\eta$ is a time homogeneous Poisson random measure
with intensity $\nu$ on  $(Z,{\mathcal{Z}})$ over $(\Omega,{\mathcal{F}},{\mathbb{F}},\mathbb{P})$. It
turns out that the compensator $\gamma$ of $\eta$ is uniquely
determined and moreover
\[
\gamma(A,I)=  \nu(A)\times\Leb(I),\;\; (A,I)\in{\mathcal{Z}} \times {\mathcal{B}}({\mathbb{R}}_+).
\]
The  difference between a time homogeneous  Poisson random measure
$\eta$  and its compensator $\gamma$, i.e.  $\tilde
\eta=\eta-\gamma$, is called a  {\em compensated Poisson random
measure}.
\end{remark}
%
%
%
%
%%%%%%%%%%%%%%%%%
%%%%%%%%%%%%%%
We have the following important result.
\begin{proposition}\label{PROP-EQ-PRM}
Let us assume that  $(Z, {\mathcal{Z}})$ be a measurable space and  $\nu \in M(Z)$. Let $\eta$ and $ \mu$  be two measurable functions from $(\Omega, \mathcal{F}) $ to $ (M([0,T]\times Z), \mathcal{B}\times {\mathcal{Z}})$
 If the  laws $\eta$ and $ \mu$ are equal and  $\mu$ satisfies the conditions \eqref{prm-i}-\eqref{prm-3} from Definition \ref{def-Prm}
 then $\eta$ satisfies them as well.
\end{proposition}
\begin{proof}
For any $A\in {\mathcal{Z}}$ let $N_\mu$ be the $\bar{\mathbb{N}}$-valued process defined by
$$N_\mu(t,A):= \mu(A \times (0,t]), \;\; t\ge 0.
$$
In addition, we denote by  $(N_\mu(t))_{t\ge 0}$ the measure
valued process defined by $N_\mu(t) =\{ {\mathcal{Z}} \ni A \mapsto
N_\mu(t,A)\in\bar{\mathbb{N}}\}$, $t\in[0,T]$. Note that for each $t\ge 0$ and $A\in {\mathcal{Z}}$ the laws of the $\mathbb{N}$-valued random variables
$N_\mu(t,A)$ and $i_{[0,t]\times A} \circ \eta$ are equal.

Now, for each $t\ge 0$ and $A\in {\mathcal{Z}}$ we have
$$ \mathbb{P}\{N_\mu(t,A)=k\} = \mathbb{P}\{N_\mu(t,A)=k\}, \quad k \in \bar{\mathbb{N}},$$ because $\eta$ and $\mu$
have the same law. Hence if $\mu$ satisfies Definition \ref{def-Prm}\eqref{prm-i} and Definition \ref{def-Prm}\eqref{prm-3},
then $\eta$ satisfies them as well.

To prove Definition \ref{def-Prm}\eqref{prm-ii} it is sufficient to show that for
any two disjoint sets $A_1,A_2\in{\mathcal{Z}}$, $t\ge 0$ and $\theta \in \mathbb{R}$, we have
\begin{eqnarray}\label{KIK-FUNC} {\mathbb{E}}
e ^{i\left(\theta_1  N_{\eta} (t,A_1)  + \theta_2  N_{
\eta} (t,A_2) \right)} =
 {\mathbb{E}} e ^{i\;\theta_1  N_{ \eta} (t,A_1) }{\mathbb{E}} e ^{i\;\theta_2  N_{ \eta} (t,A_2) }.
\end{eqnarray}
Since, by assumption, identity \eqref{KIK-FUNC} is true for $\mu$, it is also true for $\eta$ because $\mu$ and $\eta$ have the same law and
 they have the same characteristic function. %$\ref{def-Prm}\eqref{prm-ii}

\end{proof}
\begin{remark}\label{REM-EQ-PRM}
 The conclusion of the proposition is still true even if $\eta$ and $\mu$ are defined on two different probability spaces.
\end{remark}

%%%%%%%%%%%%%%%%%%%%%%
%%%%%%%%%%%%%%%%%%%%%%

The classical It\^o stochastic integral has been generalised in
several directions, for example in Banach spaces of martingale
type $p$.
Since it would exceed the scope of the paper, we have decided not
to present a detailed introduction on this topic and restrict
ourselves  only to the necessary definitions. A short summary of
stochastic integration in Banach spaces of martingale type $p$ is
given in Brze{\'z}niak \cite{MR1313905}, Hausenblas \cite{levy2},
or Brze{\'z}niak and Hausenblas \cite{maxreg}.

We finish with the following version of the {\em Stochastic Fubini
Theorem} (see \cite{zhu}).

\begin{theorem}\label{stochasticfubini}
Assume that $E$ is a Banach space of martingale type $p$ and
$$\xi:[a,b]\times \Rb{+}\times Z \to  E$$
is a progressively measurable process. Then, for each $T\in
[a,b]$,  a.s.
 \begin{eqnarray*}  \int_a^b\big[ \int_0^T \int_Z
\xi(s,r;z)\,\tilde\eta( dz;dr)\big] \,ds &=& \int_0 ^T \int_Z
\big[ \int_a^b \xi(s,r;z)\,ds \big]\, \tilde \eta( dz;dr)  \end{eqnarray*}
\end{theorem}

\medskip

\subsection{L\'evy processes and Poisson random measures}\label{subsection-rel-l-prm}

From a L\'evy process one can construct a corresponding Poisson
random measure. Conversely, given a Poisson random measure, one
gets easily a corresponding L\'evy process.
To illustrate this fact, let us recall firstly the definition of a
 L\'evy process.
\begin{definition} Let $E$ be a Banach space.
A stochastic process $L=\{ L(t):t\ge 0\}$ over a filtered probability space
$(\Omega,{\mathcal{F}}, \mathbb{F},\mathbb{P})$, where ${\mathbb{F}}=\{{\mathcal{F}}_t\}_{t\ge 0}$ is a filtration, is called an $E$-valued L\'evy process if the
following conditions are satisfied.
\begin{numlistn}
\item $L(t)$ is $\mathcal{F}_t$-measurable for any $t\ge 0$;
\item the random variable $L(t)-L(s)$ is independent of $\mathcal{F}_s$ for any  $0\le s< t$;
\item $L_0=0$ a.s.;
\item For
all $0\le s<t$, the law of $L(t+s)-L(s)$ does not depend on $s$;
\item $L$ is stochastically continuous;
\item the trajectories of $L$ are a.s. c\`{a}dl\`{a}g in $E$.
\end{numlistn}
\end{definition}

If $E$ is a
Banach space  of type $p$ with dual $E^\ast$ and $L= \{ L(t):t\ge 0\}$ is an
$E$-valued L\'evy process,
 then there exist a nonnegative symmetric operator $Q:E ^\ast \to E$,
a non-negative measure $\nu$ concentrated on $E\setminus \{0\}$
satisfying  $\int_E(1\wedge |z|^p)\,\nu(dz)<\infty$, and $m\in E$ such that
\begin{equation*}
\begin{split}
  {\mathbb{E}} e^{i\la L(1),x^\ast\ra}=  \exp \Big(
\int_E \left( 1-e ^ {i\la y,x^\ast\ra }+1_{(-1,1)}(|y|)i\la y,x^\ast\ra
\right)\nu(dy)  \Big)\\ \times \exp\Big( i \la m,x^\ast \ra -\frac 12 \la Qx^\ast,x^\ast \ra\Big) ,\quad x^\ast\in E^\ast.
\end{split}
\end{equation*}
We refer, for instance, to
\cite{apple-int,araujo,Linde,Gine} for this result.\\
The measure $\nu$ is called {\em the characteristic measure} of the L\'evy process $L$. %=\{ L(t):t\ge 0\}$.
A L\'evy process is called of {\em pure jump type}  iff $Q=0$.
Moreover, the triplet $(Q,m,\nu)$ uniquely determines the law of
the  L\'evy process.

Now, starting with an $E$--valued L\'evy process over a filtered
probability space $(\Omega,{\mathcal{F}},{\mathbb{F}}, \mathbb{P})$ one can construct an
integer valued random measure as follows. For each $ (B,I)\in
{\mathcal{B}}({\mathbb{R}})\times {\mathcal{B}}({\mathbb{R}}_+)$ let
$$
 \eta_L(B\times I ) := \#\{s\in I \mid \Delta_s L
 \in B\}\in  \deln{{\mathbb{N}}_0\cup\{\infty\}}\bar{\mathbb{N}}.
$$
The {\em jump process} $\Delta X = \{\Delta _tX,\,0\le t<\infty\}$
of a process $X$ is defined by $ \Delta_t X(t) := X(t) - X(t-)=
X(t) - \lim_{s\uparrow  t} X(s) ,\quad t> 0 $ and $\Delta _0X:=0$.
If $E={\mathbb{R}} ^ d$, it can be shown that $\eta_L$ defined above is a
time homogeneous Poisson random measure (see Theorem 19.2
\cite[Chapter 4]{sato}).

\medskip

A measure $\nu$ on a Banach space $E$ is called symmetric if $\nu(A)=\nu(-A)$ for any $A\in \mathbb{B}(E)$.
For a Borel measure $\nu$ on $E$ we denote by $\tilde{\nu}$ its symmetrisation which is defined by $\tilde{\nu}(A)=\frac12 (\nu(A)+\nu(-A))$ for any $A\in \mathbb{B}(E)$.
\begin{definition}\label{def:levy}(See \cite[Chapter 5.4]{Linde})
Let $E$ be a separable Banach space with dual $E^\ast$.
A symmetric $\sigma$-finite Borel measure $\nu$ on $E$ is called a {\sl
symmetric L\'evy measure} if and only if
\begin{trivlist}
\item[(i)]
 $\nu(\{0\} )=0$, and
\item[(ii)]
 the
function
$$
E^\ast \ni a\mapsto  \exp \left( \int_E (\cos\langle x,a\rangle
-1) \; \nu(dx)\right)
$$
is a characteristic function of a Radon measure on $E$.
\end{trivlist}
A $\sigma$-finite Borel measure $\nu$ on $E$ is
called a L\'evy measure provided its symmetrisation part $\tilde{\nu}$ is a symmetric L\'evy measure.
\end{definition}
\begin{remark}
As remarked in \cite[Chapter 5.4]{Linde} we do not need to suppose
that the integral $\int_E (\cos\langle x,a\rangle -1) \;
\nu(dx)$ is finite. However, see ibid. Corollary 5.4.2, if
$\nu$ is a symmetric L\'evy measure, then, for each $a \in E^\ast$, the integral in question is finite.
\end{remark}
\begin{remark}(See e.g.\ \cite{0347.46016})
If $E$ is a Banach space of type $p$, then the following statements are equivalents
\begin{enumerate}
 \item a nonnegative measure $\nu$
is a L\'evy measure,
\item $\nu(\{0\})=0$ and
$\int_E|z|^p\nu(dz)<\infty$,
\item there exists an $E$-valued pure jump L\'evy process $L$ such that $\nu$ is its characteristic measure.
\end{enumerate}
\end{remark}

\subsection{The space time Poissonian white noise}\label{intr-stpn}
Analogously to the space time Gaussian white noise one can construct
a space time L\'evy  white noise or space time Poissonian white
noise.
But before doing this, let us recall the definition of  a Gaussian
white noise (see e.g.\ Dalang \cite{dalang}).

\begin{definition}\label{def-Gaussian white noise}
Let $(\Omega,{\mathcal{F}},\mathbb{P})$ be a complete probability space and let
$(S,{\mathcal{S}},\sigma)$ a  measure space. A \textbf{ Gaussian white noise
on $(S,{\mathcal{S}},\sigma)$} is an ${\mathcal{F}}\big\slash {\mathcal{M}}(S)$ measurable
mapping
$$
W:\Omega \to   M(S)
$$
such that
\begin{numlistn}
\item for every $A\in{\mathcal{S}}$ such that
$\sigma(A)<\infty$,  $W(A):=i_A\circ W$ is a $N(0, \sigma(A))$
Gaussian random variable; \item
 if the sets $A_1,A_2\in{\mathcal{S}} $ are disjoint, then the random variables $W(A_1)$ and $W(A_2)$ are independent
and $W(A_1\cup A_2)= W(A_1)+W(A_2)$.
\end{numlistn}
\end{definition}

The space time
Gaussian white noise can be defined as follows. Let $\mathcal{O}\subset
{\mathbb{R}} ^ d$ be a domain. Put $S=\mathcal{O}\times [0,\infty)$,
${\mathcal{S}}={\mathcal{B}}(\mathcal{O})\otimes {\mathcal{B}}(  [0,\infty))$ and let $\leb$ be the
Lebesgue measure on $S$.
The space time Gaussian white noise is an
$M(\mathcal{O})$-valued process
$ \left\{ W_{st}(t)
: t\geq 0\right\}$ defined by
$$W_{st}(t) =\{ {\mathcal{B}}(\mathcal{O}) \ni A\mapsto W(A\times [0,t))\in {\mathbb{R}}\}
,\quad t\ge 0.
$$
If $(\Omega,{\mathcal{F}},{\mathbb{F}},\mathbb{P})$ is a filtered probability space, then we
say that the Gaussian white noise $W$ on $(S,{\mathcal{S}},\sigma)$ is a
space time Gaussian white noise over $(\Omega,{\mathcal{F}},{\mathbb{F}},\mathbb{P})$, if
the process
$ \left\{ W_{st}(t)
: t\geq 0\right\}$ is ${\mathbb{F}}$--adapted.

Moreover, one can  show that the $M(\mathcal{O})$-valued process
$ \left\{ W_{st}(t)
: t\geq 0\right\}$

generates, in a unique way,  an  $L ^ 2(\mathcal{O})$-cylindrical Wiener
process $(\hat{W}_t)_{t\geq 0}$, see \cite[Definition
4.1]{{MR1880243}}. In particular, for any $A\in \mathcal{B}(\mathcal{O})$ such that
$\leb(A)<\infty$, and any $t\geq 0$,  $\hat{W}_t(1_A)=W(A
\times [0,t))=W_{st}(t,A)$.

\medskip
A L\'evy white noise and a space time L\'evy white noise can
be defined analogously.
\begin{definition}\label{def-stl}
Let $(\Omega,{\mathcal{F}},\mathbb{P})$ be a complete probability space,
$(S,{\mathcal{S}},\sigma)$ be a measurable space, $\gamma \in \mathbb{R}$
 and  $\nu$ a L\'evy measure on  ${\mathbb{R}}$, see Definition \ref{def:levy}. Then \textbf{ a  L\'evy white
noise on $(S,{\mathcal{S}},\sigma)$  with intensity jump size measure $\nu$}
is an ${\mathcal{F}}\big\slash {\mathcal{M}}(S)$-measurable mapping
$$
L:\Omega \to  M(S)
$$
such that
\begin{numlistn}
\item for all $A\in{\mathcal{S}}$ such that $\sigma(A)<\infty$,  $L(A):=
i_A\circ L$ is an infinite divisible random variables satisfying,  \text{ for all $\theta \in{\mathbb{R}}$ },
\begin{eqnarray*}\lefteqn{ {\mathbb{E}} e ^ {i \theta L(A)}} && \\
& =&\exp\left( \sigma(A)\, \left[ i\gamma \theta +\int_{\mathbb{R}} \left( 1 - e
^ {i\theta x}- i\theta x 1_{[-1,1]}(x) \right)\, \nu(dx)\right] \right).
 \end{eqnarray*} \item If the sets
$A_1,A_2\in{\mathcal{S}} $ are disjoint, then the random variables $L(A_1)$
and $L(A_2)$ are independent and $L(A_1\cup A_2)= L(A_1)+L(A_2)$.
\end{numlistn}
\end{definition}

\begin{definition}\label{def-st-1}
Let $(\Omega,{\mathcal{F}},\mathbb{P})$ be a complete probability space. Suppose
that $\mathcal{O}\subset {\mathbb{R}} ^ d$ is a domain and let $S=\mathcal{O}\times
[0,\infty)$, ${\mathcal{S}}={\mathcal{B}}(\mathcal{O})\times {\mathcal{B}}( [0,\infty))$ and $\leb$
the Lebesgue measure on $S$. %=\leb_{d+1}$.
Let $\nu$ be a L\'evy measure on ${\mathbb{R}}$.
If $L:\Omega \to M( S)$ is a L\'evy white noise on
$(S,{\mathcal{S}},\leb)$  with intensity jump size measure $\nu$,
then the $M(\mathcal{O})$-valued process $\left\{ L_{st} (t) :
t\ge 0\right\}$ defined by
$$
L_{st}(t)=\{{\mathcal{B}}(\mathcal{O}) \ni A\mapsto L(A\times [0,t))\in {\mathbb{R}} \},\quad
t\ge 0.
$$
is called a \textbf{ space time L\'evy noise  on $\mathcal{O}$ with jump
size characteristic $\nu$}. If $(\Omega,{\mathcal{F}},{\mathbb{F}},\mathbb{P})$ is a
filtered probability space, then we say $L$ is a \textbf{ space
time L\'evy white noise} over $(\Omega,{\mathcal{F}},{\mathbb{F}},\mathbb{P})$, iff the
corresponding  measure valued
process $ \left\{ L_{st}(t)
: t\geq 0\right\}$, is
${\mathbb{F}}$--adapted.
\end{definition}
\begin{remark}
If $L$ is a space time L\'evy white noise on $(S,{\mathcal{S}},\leb)$ as
given in Definition \ref{def-st-1}, then the  corresponding $M(\mathcal{O})$-valued
process $ \left\{ L_{st}(t)
: t\geq 0\right\}$
is a weakly cylindrical process on $L ^ 2(\mathcal{O})$, see Definition
3.2 \cite{markus}.
\end{remark}

In the very same way, one can first define a
Poissonian white noise, and then, a   space time Poissonian  white
noise.
\begin{definition}\label{def-stp}\label{st-stp}
Let us assume that $(\Omega,{\mathcal{F}},\mathbb{P})$ is a complete probability
space and $\nu$ be a L\'evy measure on ${\mathbb{R}}$.
 \begin{trivlist}
 \item[(a)]
A \textbf{  Poissonian white noise  with intensity jump size
measure $\nu$ on   a  measurable space  $(S,{\mathcal{S}},\leb)$ } is a
${\mathcal{F}}\big\slash{\mathcal{M}}(M_I(S\times{\mathbb{R}}))$ measurable mapping
\[
 \eta:\Omega\to M_I(S\times {\mathbb{R}})
\]
such that
\begin{numlistn}
\item  for all $A \in {\mathcal{S}}\otimes  {\mathcal{B}}({\mathbb{R}})$, $\eta(A):=i_{A}\circ
\eta$ is a Poisson  random variable with parameter  $(\leb
\otimes\nu)(A)$, provided $(\leb \otimes\nu)(A)<\infty$;
%and $\ nu\times \leb(B)<\infty$;
\item if the sets $A_1\in {\mathcal{S}}\otimes  {\mathcal{B}}({\mathbb{R}}) $
and $A_2\in {\mathcal{S}}\otimes {\mathcal{B}}({\mathbb{R}})$
are disjoint, then the random variables $\eta(A_1)$ and
$\eta(A_2)$ are independent and $\eta\left(A_1\cup A_2 \right)=
\eta(A_1)+\eta(A_2)$, almost surely.
\end{numlistn}

 \item[(b)] Let   $\mathcal{O}\subset {\mathbb{R}} ^ d$ be a domain. Then $\eta$ is called  a \textbf{space time Poissonian white noise on $\mathcal{O}$ with
intensity jump size measure $\nu$} iff $\eta$ is  a
\textit{Poissonian white noise  on $(S,{\mathcal{S}},\leb)$}  with
intensity jump size measure
$\nu$, where $S=\mathcal{O}\times [0,\infty)$, ${\mathcal{S}}={\mathcal{B}}(\mathcal{O})\otimes {\mathcal{B}}(  [0,\infty))$. \\
The corresponding measure-valued process $\{\Pi_t:t\ge 0\}$
defined by \begin{eqnarray}\label{mvp} \Pi_t: {\mathcal{B}}(\mathcal{O})\times {\mathcal{B}}({\mathbb{R}})\ni
(A,B) \mapsto
 \eta ( A\times  [0,t)\times B )\in \bar{\mathbb{N}}_0.
\end{eqnarray} is called   the \textbf{(homogeneous) space time Poissonian
white noise process}. \item[(c)] If $(\Omega,{\mathcal{F}},{\mathbb{F}},\mathbb{P})$ is a
filtered probability space, then  a \textit{space time Poissonian
white noise $\eta$ on $\mathcal{O}$} is called
 a
\textbf{ (homogeneous) space time Poissonian white noise over
$(\Omega,{\mathcal{F}},{\mathbb{F}},\mathbb{P})$}, iff  the measure valued process
$\{\Pi_t:t\ge 0\}$ defined above is ${\mathbb{F}}$--adapted. (Compare this
definition with \cite[Definition 7.2]{Peszat_Z_2007}.)
\end{trivlist}
\end{definition}

\begin{theorem}\label{thm-ZB}  Let us assume that
$(\Omega,{\mathcal{F}},\mathbb{P})$ is a complete probability space,  $\nu$ be a
L\'evy measure on ${\mathbb{R}}$ and
  $\mathcal{O}\subset {\mathbb{R}} ^ d$ be a bounded set and
\[\eta:  \Omega\to M_I(S\times {\mathbb{R}}) \] is
  a \textbf{space time Poissonian white noise on $\mathcal{O}$ with
intensity jump size measure $\nu$}, where $S=\mathcal{O}\times
[0,\infty)$, ${\mathcal{S}}={\mathcal{B}}(\mathcal{O})\otimes {\mathcal{B}}(  [0,\infty))$ and $\leb$
the Lebesgue measure on $S$. Define a random measure
\[\gimel :\Omega\to M_I(\mathcal{O} \times {\mathbb{R}} \times [0,\infty) )\]
by
\[ \gimel (\omega)(A\times B \times C)= \eta(\omega)(A \times C\times B), \;\; A \in {\mathcal{B}}(\mathcal{O}), B\in  {\mathcal{B}}({\mathbb{R}}), \in  {\mathcal{B}}([0,\infty).
\]
Then, $\gimel$ is a time homogeneous Poisson random measure on
$\mathcal{O} \times {\mathbb{R}}$ with the intensity measure $\daleth$ satisfying
\[
\daleth: {\mathcal{B}}( \mathcal{O} \times {\mathbb{R}})\ni B \mapsto   \int_{\mathbb{R}} \int _\mathcal{O}
1_B(\xi,\zeta) \, \nu(d\zeta)\, d\xi.
\]
\end{theorem}
\begin{proof}
The proof is similar to the proof of \cite[Proposition 7.21]{Peszat_Z_2007}.
\end{proof}

\begin{remark}
The compensator $\gamma$ of a  homogeneous space time Poisson
(white) noise with jump size intensity $\nu$ is a measure on
$\mathcal{O}\times {\mathbb{R}}\times [0,\infty)$ defined  by
$$
{\mathcal{B}}(\mathcal{O})\times {\mathcal{B}}({\mathbb{R}})\times {\mathcal{B}}([0,\infty))
 \ni (A,B ,I  ) \mapsto  \gamma(A\times B\times I ) =\leb(A)\, \nu(B) \, \leb(I).
$$
\end{remark}

%%%%%%%%%%%%%%%%%%%%%%%%%
%%%%%%%%%%%%%%%%%%%%%%%%%
%%%%%%%%%%%%%%%%%%%%%%%%%
%%%%%%%%%%%%%%%%%%%%%%%%%
\section{Stochastic Integration - a useful  result}

As seen in \cite{maxreg} the stochastic integral is a unique
bounded extension of the integral defined for the class of simple
functions to a set of all progressively measurable functions. In
the same way we will show the following Proposition.

\begin{proposition}\label{zero}
{Assume that  $(S,{\mathcal{S}})$ is  a measurable space, $\nu$
is a non-negative measure on $(S,{\mathcal{S}})$ and
$\mathfrak{P}=(\Omega,{\mathcal{F}},(\mathcal{F}_t)_{t\geq 0},\mathbb{P})$ is a
filtered probability space.} Assume also that  $\eta_1$ and
$\eta_2$ are two time homogeneous Poisson random measures  over
$\mathfrak{P}$,  with the intensity measure $\nu$.
 Assume that $p\in (1,2]$ and $E$ is a martingale type $p$ Banach space
and $\xi\in  \mathcal{M}^p(0,\infty, L^p(S,\nu;E))$. If
$\mathbb{P}$-a.s.\ $\eta_1=\eta_2$ on $M_{\mathbb{N}}(S,[0,T])$, then for any
$t\ge 0$ $\mathbb{P}$-a.s. \begin{eqnarray}\label{equal...} \int_0 ^ t \xi(s,z)
\,\tilde  \eta_1(dz,ds) = \int_0 ^ t \xi(s,z) \,\tilde
\eta_2(dz,ds). \end{eqnarray}
\end{proposition}

\begin{proof}[ Proof of Proposition \ref{zero}] In order to show the equality, we will go back to the definition
of the stochastic integral by step functions. First, put
$I=(a,b]$. We may suppose that
$f=\sum_{j=1,i=1}^{I,J}f_{ji}1_{A_{ji}\times B_i}$  with
$f_{ji}\in E$, $A_{ji}\in \mathcal{F}_{a}$ and $B_i\in{\mathcal{S}}$. For
fixed $i$ we may suppose that the finite family of sets
$\{{A_{ji}\times B_i}:j=1,\cdots, J\}$ are pairwise disjoint and
that the family of sets $\{B_i:i=1,\dots I\}$ are also  pair-wise
disjoint and $\nu(B_i)<\infty$. Let us notice that for
$i=1,\cdots, I$
$$
  \int_{B_i} f (x)\, \tilde{\eta_1}(dx,I)=\sum_{j=1}^J1_{A_{ji}} \, \tilde \eta_1({B_i\times I }) \, f_{ji},$$
and
$$
  \int_{B_i}  f (x)\, \tilde{\eta_2}(dx,I)=\sum_{i=1}^J 1_{A_{ji}} \,\tilde \eta_2({B_i\times I }) \,f_{ji}.$$
Since $\mathbb{P}$-a.s.\ $\eta_1=\eta_2$ on $M_{\mathbb{N}}(S,[0,T])$, it
follows that a.s.\ $ \int_S f (x)\tilde{\eta_1}(dx,I)=\int_S f
(x)\tilde{\eta_2}(dx,I)$. It remains to investigate the limit as
done in \cite[Appendix C]{maxreg}. But if $x_n=y_n$ for all
$n\in{\mathbb{N}}$ and $x_n\to x$, $y_n\to y$ in a Banach space $X$, then
$x=y$.
\end{proof}
%%%%%%%%%%%%%%%%%%%%%%%%%%%%%%%%%%%%%%%%%%%%%%%%%%%%%%%%%%%%%%%%%%%%%%%%%%%%%%%%%%%%%%%%%%%%%%

%%
\section{Besov spaces and their properties}\label{besov}

An introduction to Besov spaces are given in Runst and Sickel
\cite[p.\ 8, Def.\ 2]{873.35001}. We are interested in the continuity of the
mapping $G$ described in \ref{assum-main-2}. To be  precise, we
will prove the following result.
\begin{proposition}\label{prop-delta}
Let $p,p^\ast\in (1,\infty)$ satisfy $\frac1 p+ \frac1 {p^\ast}=1$.  Then for every  $f\in {\mathcal{S}}({\mathbb{R}}^ {d})$  and $a\in \mathbb{R}^d$ the tempered
distribution $f\delta_a$ belongs to the Besov space $B^{-\frac
d{p^\ast}}_{p,\infty}({\mathbb{R}} ^ d)$ and
\begin{eqnarray}\label{eqn-delta}
 \int_{{\mathbb{R}} ^d} \left| f\delta_a\right|_{B^{-\frac d{p^\ast}}_{p,\infty}({\mathbb{R}} ^ {d}) }^p da = (2\pi)^{-\frac d2}2^{-\frac d {p^\ast}}\left|f\right| ^p _{L ^p({\mathbb{R}} ^ {d})}. \end{eqnarray}
In particular, there exists a unique bounded linear map
 $$\Lambda: L^p(\mathbb{R}^d) \to L^p({\mathbb{R}} ^ {d}, B^{-\frac d{p^\ast}}_{p,\infty}({\mathbb{R}} ^ {d}))
$$
 such that $[\Lambda(f)](a)=f \delta_a$, $f\in {\mathcal{S}}({\mathbb{R}} ^ {d})$, $a \in {\mathbb{R}} ^ {d}$.
\dela{Moreover, for any $s\ge 0$ there exists a constant $C$ such that
\begin{eqnarray}\label{zweites}
 \int_{{\mathbb{R}} ^d} \left| f\delta_a\right|_{B^{-\frac d{p^\ast}-s}_{p,\infty}({\mathbb{R}} ^ {d}) }^ p \: da \le C\: \left|f\right|  _{W ^{-s}_p ({\mathbb{R}} ^ {d})} ^ p .
\end{eqnarray}
}
\end{proposition}
In the following we denote the value of  $\Lambda(f)$ at $a$,
where $f \in L^p({\mathbb{R}} ^ {d})$ by $f \delta_a$. Note that $\lb f\delta_a, \phi \rb=f(a) \phi(a)$, ${\mathcal{S}}({\mathbb{R}}^ {d})$ so that $f\delta_a= f(a) \delta_a$. Hence, in order to prove \eqref{eqn-delta} it is sufficient to prove it for $f=1$, i.e. that
\begin{equation}\label{eqn-delta_a}
\left| \delta_a\right|_{B^{-\frac d{p^\ast}}_{p,\infty}({\mathbb{R}} ^ {d}) }^p=(2\pi)^{-\frac d2}2^{-\frac d {p^\ast}}.
\end{equation}

Let us recall the
definition of the Besov spaces  as  given in \cite[Definition 2,
pp. 7-8]{873.35001}. First we choose a function
$\psi\in{\mathcal{S}}(\mathbb{R}^d)$ such that $0\le \psi(x)\le 1$, $x\in
\mathbb{R}^d$ and
$$
\psi(x) = \left\{ \begin{array}{rcl} 1,&\mbox{ if } & |x|\le 1,\\
0&\mbox{ if } & |x|\ge \frac 32.
\end{array}\right.
$$
Then put \begin{eqnarray*}
 \left\{ \begin{array}{rcl}\phi_0(x) &=&\psi(x), \; x\in \mathbb{R}^d,
\\
\phi_1(x) &=&\psi(\frac x2)-\psi(x), \; x\in \mathbb{R}^d,
\\
\phi_j(x) &=&\phi_1(2 ^{-j+1} x),\; x\in \mathbb{R}^d, \quad
j=2,3,\ldots.
\end{array} \right.
\end{eqnarray*}
We will use
the definition of the Fourier transform ${\mathcal{F}}={\mathcal{F}}^{+1}$ and its inverse ${\mathcal{F}}^{-1}$ as in \cite[p.\
6]{873.35001}. In particular,  with  $\lb \cdot,\cdot \rb$ being
the scalar product in ${\mathbb{R}} ^ {d}$, we put
\begin{equation*}({\mathcal{F}}^{\pm 1} f)(\xi):= (2\pi)^{-d/2} \int_{{\mathbb{R}} ^ {d}} e^{\mp i\lb x,\xi\rb}f(x)\, dx,\;
\; f \in {\mathcal{S}}({\mathbb{R}} ^ {d}),\xi\in {\mathbb{R}} ^ {d}.
\end{equation*}

With the choice of $\phi=\{\phi_j\}_{j=0} ^\infty$ as above and
${\mathcal{F}} $ and ${\mathcal{F}} ^{-1}$ being the Fourier  and the inverse Fourier
transformations (acting in the space ${\mathcal{S}}^\prime(\mathbb{R}^d)$ of
Schwartz distributions)  we have the following definition.
\begin{definition}
Let $s\in{\mathbb{R}}$,  $0<p\le \infty$ and and  $f \in {\mathcal{S}}^\prime(\mathbb{R}^d)$. If $0<q <
\infty$  we put
\begin{eqnarray*} \left| f\right|_{B^{s}_{p,q}} &=& \left( \sum_{ j=0} ^\infty  2
^{sjq}\left| {\mathcal{F}} ^{-1} \left[ \phi_j{\mathcal{F}} f\right] \right|_{L ^p} ^q\right)
^\frac 1q = \Vert\Big( 2 ^{sj}\left| {\mathcal{F}}
^{-1} \left[ \phi_j{\mathcal{F}} f\right]\right|_{L ^p}\Big)_{j\in \mathbb{N}} \Vert_{l^q}.
\end{eqnarray*}
If  $q=\infty$ we put
\begin{eqnarray*}
\left| f\right|_{B^{s}_{p,\infty}} &=& \sup_{j\in{\mathbb{N}}}  2 ^{sj}\left| {\mathcal{F}}
^{-1} \left[ \phi_j{\mathcal{F}} f\right]\right|_{L ^p}=\Vert \Big( 2 ^{sj}\left| {\mathcal{F}}
^{-1} \left[ \phi_j{\mathcal{F}} f\right]\right|_{L ^p}\Big)_{j\in \mathbb{N}} \Vert_{l^\infty} . \end{eqnarray*} We denote by
$B^{s}_{p,q}(\mathbb{R}^d)$ the space of all $f \in
{\mathcal{S}}^\prime(\mathbb{R}^d)$ for which $\left| f\right|_{B^{s}_{p,q}}$ is
finite.
\end{definition}
\dela{\begin{lemma}\label{lemma_con}
If $h\in {\mathcal{S}}(\mathbb{R}^d)$, $g\in {\mathcal{S}}(\mathbb{R}^d)$ $a\in
\mathbb{R}^d$, then
$$
\left|  (h \delta_a) \ast g\right|_{L ^p(\mathbb{R}^d)} = |h(a)|
|g|_{L^p(\mathbb{R}^d)}.
$$
\end{lemma}
\begin{proof} Let us recall that $h\delta_b\in {\mathcal{S}}^\prime({\mathbb{R}} ^ {d})$  assigns to
a function $g\in {\mathcal{S}}({\mathbb{R}} ^ {d})$ a value $\delta_a(hg)=h(a)g(a)$.
Since by the  definition of a convolution of a distribution with a
test
function, where $\check{g}(\cdot)=g(-\, \cdot)$ %, see e.g.\ \cite[section 6.29]{MR1157815}
$$[(h\delta_a)\ast g](x)=(h\delta_a)(\tau_x\check{g})=h(a)(\tau_x\check{g})(a)=h(a)\check{g}(a-x)=h(a){g}(x-a),
\; x\in {\mathbb{R}} ^ {d},
$$
we have
\begin{eqnarray*} \left|  (h\delta_a)\ast g\right|_{L ^p(\mathbb{R}^d)}^p &=&
\int_{{\mathbb{R}} ^d }\left| g(x-a) h(a) \right| ^p \,dx
%\\
=  %| h(a)| ^p \: \int_{{\mathbb{R}} ^d }\left| g(x)  \right| ^p dx = | h(a)| ^p
\:\left| g  \right|_{L ^p(\mathbb{R}^d)} ^p. \end{eqnarray*}
\end{proof}}

\begin{lemma}\label{lemma_con1}
If $\vp\in {\mathcal{S}}(\mathbb{R}^d)$, $\lambda>0$ and $g(x):=\vp(\lambda
x)$, $x \in {\mathbb{R}} ^ {d}$, then
$$
|  {\mathcal{F}}^{-1}g|_{L ^p(\mathbb{R}^d)} =\lambda^{d(\frac1p-1)}
|{\mathcal{F}}^{-1}\vp|_{L^p(\mathbb{R}^d)}.
$$
\end{lemma}
\begin{proof}The proof follows from simple calculations  so it is omitted.
\end{proof}
\begin{proof}[Proof of Proposition \ref{prop-delta}] As remarked earlier it is enough to show equality \eqref{eqn-delta_a}. Since   ${\mathcal{F}}^{-1}(\vp
u)=(2\pi)^{-d/2}({\mathcal{F}}^{-1}\vp)\ast ({\mathcal{F}}^{-1} u)$, $\vp \in {\mathcal{S}}$, $u
\in {\mathcal{S}}^\prime$ we infer that for $j\in \mathbb{N}^\ast$,
\begin{eqnarray*} |  {\mathcal{F}}^{-1}[\phi_j{\mathcal{F}}(\delta_a)|_{L ^p(\mathbb{R}^d)} &=&
(2\pi)^{-d/2}|  ({\mathcal{F}}^{-1}\phi_j)\ast \delta_a|_{L
^p(\mathbb{R}^d)} \\
&=&(2\pi)^{-d/2} |{\mathcal{F}}^{-1}\phi_j|_{L ^p(\mathbb{R}^d)}
\\&=&(2\pi)^{-d/2}2^{d(\frac1p-1)}2^{-jd(\frac1p-1)}  |{\mathcal{F}}^{-1}\phi_1|_{L
^p(\mathbb{R}^d)}.
 \end{eqnarray*}
Hence, $\delta_a $ belongs to the Besov space
$B_{p,\infty}^{d(\frac1p-1)}(\mathbb{R}^d)$ as requested and the
equality (\ref{eqn-delta_a}) follows immediately.
\end{proof}
\begin{corollary}\label{cor-delta}Assume also that
 $\mathcal{O}$ is a bounded and open subset of $\mathbb{R}^d$ with
boundary $\partial \mathcal{O}$ of class $\mathcal{C}^\infty$. Let $
r, q\in (1,\infty)$ with $q\ge r$ then there exists a unique
bounded linear map
 \begin{equation}\label{eqn-delta-O}
 \Lambda: L^q(\mathcal{O}) \to L^q(\mathcal{O}, B^{-\frac d{r}}_{r,\infty}(\mathcal{O}))
\end{equation}
 such that $[\Lambda(f)](a)=f \delta_a$, $f\in L^q(\mathcal{O})$, $a \in
 \mathcal{O}$. In particular, there exists a constant $C$ such that  for any $f\in L^q(\mathcal{O})$
\begin{eqnarray}\label{zweites-bounded}
 \int_{\mathcal{O}} \left| f\delta_a\right|_{B^{-d(1-\frac1r)}_{r,\infty}(\mathcal{O}) }^ q \: da \le C\: \left|f\right|_{L^q(\mathcal{O})}^{q} .
\end{eqnarray}
\end{corollary}
\begin{proof} It is enough to prove \eqref{eqn-delta-O} for any $f\in C_0^\infty(\mathcal{O})$ as that set is dense in $L^q(\mathcal{O})$.
As before, we first need to show
 the following version of  \eqref{eqn-delta_a}
\begin{equation}\label{eqn-delta_a-O}
\sup_{a \in \mathcal{O}}\left| \delta_a\right|_{B^{-\frac
d{r^\ast}}_{r,\infty}(\mathcal{O}) }<C(r,d),
\end{equation}
for a constant $C(r,d)>0$ depending only on $r$ and $d$.
For that aim let us fix $a\in \mathcal{O}$ and let us recall that
according to Definition 4.2.1 from \cite{Triebel_1995}, $\left|
\delta_a\right|_{B^{-\frac d{r^\ast}}_{r,\infty}(\mathcal{O}) }$ is equal to
infimum of $\left|  u \right|_{B^{-\frac
d{r^\ast}}_{r,\infty}(\mathbb{R}^d) }$ over all $u \in B^{-\frac
d{r^\ast}}_{r,\infty}(\mathbb{R}^d) $ such that $u_{\vert
\mathcal{O}}=\delta_a$. Thus $\left| \delta_a\right|_{B^{-\frac
d{r^\ast}}_{r,\infty}(\mathcal{O}) } \leq \left| \delta_a\right|_{B^{-\frac
d{r^\ast}}_{r,\infty}(\mathbb{R}^d) }$ and the result follows by
applying \eqref{eqn-delta_a}.

Second, let $\Lambda $ be the linear map defined on $L^q(\mathcal{O})$ by $\Lambda f=f\delta_a$ for $f\in L^q(\mathcal{O})$ and $a\in \mathcal{O}$.
Since, by the assumption $q\ge r$, $L^q(\mathcal{O})\subset L^r(\mathcal{O})$ it follows from the first part of the
proof that
\begin{align*}
 \int_\mathcal{O} \vert \Lambda f(a)\vert^q_{B^{-\frac
d{r^\ast}}_{r,\infty}(\mathcal{O}) } da \le & \int_\mathcal{O} \vert f(a) \vert^q \vert \delta_a\vert_{B^{-\frac
d{r^\ast}}_{r,\infty}(\mathcal{O}) } da\\
\le & \sup_{a\in \mathcal{O}} \vert \delta_a\vert_{B^{-\frac
d{r^\ast}}_{r,\infty}(\mathcal{O}) } \int_\mathcal{O} \vert f(a) \vert^q da\\
\le & C(r,d)^q \lvert f\rvert^q_{L^q(\mathcal{O})}.
\end{align*}
The last inequality completes the proof of Corollary \ref{cor-delta}.
\end{proof}

\section{A modified version of the Skorokhod embedding Theorem}

Within the proof we are considering the limit of pairs of random
variables. For us it was important that certain propetries of the
pairs are conserved by the Skorokhod embedding Theorem. Therefore,
we had to use a modified version wich is stated below.
\begin{theorem}\label{thm-Skorokhod}
Let $(\Omega,{\mathcal{F}},\mathbb{P})$ be a probability space and $E_1,E_2$ be two
separable metric spaces. Let $\chi_n:\Omega\to E_1\times E_2$,
$n\in{\mathbb{N}}$, be a family of random variables, such that the sequence
$\{\Law(\chi_n): {n\in{\mathbb{N}}}\}$ is weakly convergent on $E_1\times
E_2$.

For $i=1,2$ let $\pi_i:E_1\times E_2$ be the projection onto
$E_i$, i.e.
$$
E_1\times E_2\ni \chi=(\chi ^ 1, \chi ^ 2) \mapsto \pi_i(\chi)=
\chi  ^ i \in E_i.
$$
Finally let us assume  that
there exists a random variable $\rho:\Omega\to E_1$
 such that
$\Law(\pi_1\circ \chi_n) = \Law(\rho)$, $\forall n\in{\mathbb{N}}$.

Then, there exists a probability space $(\bar \Omega,\bar {\mathcal{F}},\bar{\mathbb{P}})$, a family of $E_1\times E_2$--valued random variables
$\{\bar \chi_n: n\in{\mathbb{N}}\}$, on $(\bar \Omega,\bar {\mathcal{F}},\bar{\mathbb{P}})$
and a  random variable $\chi_\ast:\bar \Omega\to E_1\times E_2$
such that
\begin{numlist}
\item $\Law(\bar \chi_n) = \Law(\chi_n)$, $\forall n\in{\mathbb{N}}$;
\item $\bar\chi_n\to \chi_\ast$ in $ E_1\times E_2$ $\bar{\mathbb{P}}$--a.s.\ \item
$\pi_1\circ \bar\chi_n(\bar \omega)= \pi_1 \circ \chi_\ast (\bar
\omega)$ for all $\bar\omega\in\bar \Omega$.
%$\bar{\mathbb{P}}$ a.s.\ \Red{I have the impression, we can even omit a.s.}
\end{numlist}
\end{theorem}

\begin{proof}[Proof of Theorem \ref{thm-Skorokhod}.]
The proof is a modification of the proof of \cite[Chapter 2,
Theorem 2.4]{761.60052}. For simplicity, let us put  ${\mathcal{P}}MU_n :=
\Law(\chi_n)$, ${\mathcal{P}}MU_n ^ 1:= \Law(\pi _1 \circ \chi_n)$,
$n\in{\mathbb{N}}$, and ${\mathcal{P}}MU_\infty :=\lim_{n\to\infty} {\mathcal{L}}(\chi_n)$. We
will generate  families of partitions of $E_1$ and $E_2$. To start
with let $\{ x_i, i\in{\mathbb{N}}\}$ and $\{ y_i,i\in{\mathbb{N}}\}$ be dense
subsets in $E_1$ and $E_2$, respectively, and let
$\{r_n,n\in{\mathbb{N}}\}$ be  a sequence of natural numbers convergent to
zero. Some additional condition on the sequence will be given
below.

Now, let $A_1 ^ 1:= B(x_1,r_1)\footnote{For $r>0$ and $x$ let
$B(x,r):=\{y, |y|\le r\}$.} $, $A_k ^ 1 := B(x_k,r_1)\setminus
\left(\cup_{i=1} ^ {k-1} A_i ^ 1\right) $ for $k\ge 2$. Similarly, $C_1
^ 1:= B(y_1,r_1)$, $C_k ^ 1 := B(y_k,r_1)\setminus \left( \cup_{i=1}
^ {k-1} C_i ^ 1\right)$ for $k\ge 2$. Inductively, we put \begin{eqnarray*}
A_{i_1,\cdots, 1} ^ k &:=& A_{i_1,\cdots, i_{k-1}} ^ {k-1}\cap
B(x_{i_k},r _ k),
\\
A_{i_1,\cdots, i_k} ^ k &:=& A_{i_1,\cdots, i_{k-1}} ^ {k-1}\cap
B(x_{i_k},r _ k) \setminus \left( \cup_{j=1} ^ {i_k-1}
A_{i_1,\cdots,i_{k-1}, j} ^ k\right) ,\quad k\ge 2, \end{eqnarray*} and
similarly, where we will replace $"A"$ by $"C"$ \begin{eqnarray*}
C_{i_1,\cdots, 1} ^ k &:=& C_{i_1,\cdots, i_{k-1}} ^ {k-1}\cap
B(y_{i_k},r _ k),
\\
C_{i_1,\cdots, i_k} ^ k &:= &C_{i_1,\cdots, i_{k-1}} ^ {k-1}\cap
B(y_{i_k},r _ k) \setminus \left(  \cup_{j=1} ^ {i_k-1}
C_{i_1,\cdots,i_{k-1}, j} ^ k\right),\quad k \ge 2 . \end{eqnarray*} For
simplicity, we enumerate for any $k\in{\mathbb{N}}$ these families  and
call them $(A^ k_i)_{ i\in{\mathbb{N}}}$, and  $(C^ k_j)_{j\in{\mathbb{N}}}$.

Let $\bar \Omega:= [0,1)\times [0,1)$ and $\leb$ %, the probability measure on $\bar \Omega$,
be the Lebesgue measure on $[0,1)\times [0,1)$. In the first step,
we will construct a family of partition consisting of rectangles
in $\bar \Omega$.
\begin{definition}\label{def-part}
Suppose that $\mu$ is a Borel probability measure on $E=E_1\times
E_2$ and $\mu^1$ is the marginal of $\mu$ on $E_1$, i.e.
$\mu^1(A):=\mu(A\times E_2),$ $A\in\mathcal{B}(E_1)$. Assume that
$(A_i)_{i\in\mathbb{N}}$ and $(C_i)_{i\in\mathbb{N}}$ are
partitions of  $E_1$ and $E_2$, respectively. Define the following
partition of the square $[0,1)\times [0,1)$. For $i,j\in
\mathbb{N}$ we put
\begin{eqnarray}
\nonumber \label{egn-part}
I_{ij}&:=&\left[\mu_1\Big(\bigcup_{\alpha=1}^{i-1} A_\alpha\Big),
\mu_1\Big(\bigcup_{\alpha=1}^{i} A_\alpha\Big) \right)
\\
 && \times
\left[\frac1{\mu_1(A_i)}\mu\Big(A_i \times
\bigcup_{\alpha=1}^{j-1} C_\alpha\Big), \frac1{\mu_1(A_i)}
\mu\Big(A_i \times \bigcup_{\alpha=1}^{j} C_\alpha\Big) \right)\;
.
\end{eqnarray}
\end{definition}
\begin{remark}
 Obviously, if $ \mu_1(A_i)=0$ for some
$i\in\mathbb{N}$, then $I_{ij}=\emptyset$ for all
$j\in\mathbb{N}$.
\end{remark}
Next for fixed $l\in\mathbb{N}$ and $n\in\deln{{\mathbb{N}}_0\cup\{\infty\}}\bar{\mathbb{N}}$, we
will define a partition
$\Big(I_{ij}^{l,n}\Big)_{i,j\in\mathbb{N}}$
 of $\bar\Omega=[0,1)\times [0,1)$ corresponding to partitions $(A_i^l)_{i\in\mathbb{N}}$
 and $(C_i^l)_{i\in\mathbb{N}}$
 of the spaces $E_1$ and $E_2$, respectively.

We denote by  $\mu(A|C)$ the conditional probability of $A$ under
$C$.
Then, we have for $n\in\deln{{\mathbb{N}}_0\cup\{\infty\}}\bar{\mathbb{N}}$ \begin{eqnarray*} && I_{1,1} ^{1,
n} := \left[0, {\mathcal{P}}MU ^ 1_n(A_1 ^ 1)\right)\times \left[
0, {\mathcal{P}}MU _ n(E_1\times C_1 ^ 1\mid A_1 ^ 1\times E_2)\right), \\
&& I_{2,1} ^{1, n} := \left[ {\mathcal{P}}MU ^ 1_n(A_1 ^ 1),{\mathcal{P}}MU ^ 1_n(A_1 ^
1)+{\mathcal{P}}MU ^ 1_n(A_2 ^ 1)   \right)\times \left[0, {\mathcal{P}}MU _n(E_1\times C_1
^ 1\mid A_2 ^ 1\times E_2)\right)
\\
\; & &\ldots \;\ldots \, \end{eqnarray*} \begin{eqnarray*} &&
 I_{k,1} ^{1, n} := \left[\sum_{m=1} ^ {k-1}  {\mathcal{P}}MU ^ 1_n(A_k ^ 1),\sum_{m=1} ^ {k}  {\mathcal{P}}MU ^ 1_n(A_k ^ 1)   \right)\times
\left[0, a_{n.k}\right),
\,\,\, k\ge 2, \end{eqnarray*} and \begin{eqnarray*} &&I_{1,2} ^{1, n} := \left[0, {\mathcal{P}}MU ^
1_n(A_1 ^ 1)\right)\times \left[a_{n,1}, b_n\right),
\\ && I_{2,2} ^{1, n} := \left[ {\mathcal{P}}MU ^ 1_n(A_1 ^ 1),{\mathcal{P}}MU ^ 1_n(A_1 ^ 1)+{\mathcal{P}}MU ^ 1_n(A_2 ^ 1)   \right)\times \left[ a_{n,1}, b_n\right),
\end{eqnarray*}
where $a_{n,k}={\mathcal{P}}MU _n(E_1\times C_1 ^ 1\mid A_k ^ 1\times E_2)$ and
$$b_n={\mathcal{P}}MU _n(E_1\times C_1 ^ 1\mid A_1 ^ 1\times E_2)+ {\mathcal{P}}MU _n(E_1\times C_2 ^ 1\mid A_1 ^ 1\times E_2).$$
More generally, for $k\in{\mathbb{N}}$ \begin{eqnarray*} && I_{k,2} ^{1, n} :=
\left[\sum_{m=1} ^ {k-1}  {\mathcal{P}}MU ^ 1_n(A_k ^ 1),\sum_{m=1} ^ {k}
{\mathcal{P}}MU ^ 1_n(A_k ^ 1)   \right)\times \left[a_n, b_n \right),
\end{eqnarray*}
and, for $k,r\in{\mathbb{N}}$ %in general,
\begin{eqnarray*} &&I_{k,r} ^{1, n} :=  \left[\sum_{m=1} ^ {k-1}  {\mathcal{P}}MU ^
1_n(A_m ^ 1),\sum_{m=1} ^ {k}  {\mathcal{P}}MU ^ 1_n(A_m ^ 1)   \right)\times
\\ &&\quad
 \left[\sum_{m=1} ^ {r-1}  {\mathcal{P}}MU_n(E_1\times C_m ^ 1\mid A_k ^ 1\times E_2),\sum_{m=1} ^ {r}  {\mathcal{P}}MU _n
(E_1\times C_r ^ 1\mid A_k ^ 1 \times E_2)   \right). \end{eqnarray*} Finally,
for  $k,l,r\ge 2$ \begin{eqnarray}&&\label{relation-a}\\\nonumber&& I_{k,r}
^{l, n} :=  \left[\sum_{m=1} ^ {k-1}  {\mathcal{P}}MU ^ 1_n(A_m ^
l),\sum_{m=1} ^ {k}  {\mathcal{P}}MU ^ 1_n(A_m ^ l)   \right)\times
\\ &&\quad\quad
 \left[\sum_{m=1} ^ {r-1}  {\mathcal{P}}MU_n(E_1\times C_m ^ l\mid A_k ^ l\times E_2),\sum_{m=1} ^ {r}  {\mathcal{P}}MU _n\nonumber
(E_1\times C_m ^ l\mid A_k ^ l \times E_2)   \right).
%\\&&\nonumber \quad
\end{eqnarray} Let us observe that  for fixed $l\in{\mathbb{N}}$, the rectangles
$\{ I_{k,r} ^ l,k,r\in{\mathbb{N}}\}$ are pairwise disjoint and  the family
$\{ I_{k,r} ^ l:k,r\in{\mathbb{N}}\}$ is a covering
of  $\bar \Omega$.
Therefore, we conclude that for any $n\in{\mathbb{N}}\cup\{\infty\}$ we
have  ${\mathcal{P}}MU_n(E_1\times E_2)=1$ and $\sum_{m\in{\mathbb{N}}} {\mathcal{P}}MU_n(
E_1\times C_m ^ l\mid A_k ^ l)=1$. Hence, as a consequence,  it
follows that  for fixed $l,n\in{\mathbb{N}}$ the family of sets $\{
I_{k,r}: k,r\in{\mathbb{N}}\}$ is a covering of $[0,1)\times [0,1)$ and
consists of disjoint sets.

The next step is to construct the random variables $\bar
\chi_n:\bar \Omega \to E_1\times E_2$, such that $\Law(\bar
\chi_n)=\Law( \chi_n)$. We assume that $r_m$ is chosen in such a
way, that the measure of the boundaries of the covering
$(A_j)_{j\in{\mathbb{N}}}$
and $(C_j)_{j\in{\mathbb{N}}}$ are zero.
Now, in each non-empty sets  $\inte(A ^ {m}_j)$ and $\inte(C ^
{m}_j)$ we choose points $x_j ^ m$ and $y_j ^ m$, respectively,
from the dense subsets $\{x_i,i\in{\mathbb{N}}\}$ and $\{y_i,i\in{\mathbb{N}}\}$ and
define the following random variables. First, we put for $m\in{\mathbb{N}}$
%First, let
\begin{eqnarray*} Z ^ 1_{n,m}(\bar \omega) &=& x_k ^ m \mbox{  if  }  \bar
\omega\in I ^ {m,n}_{k,r},
\\ Z ^ 2_{n,m}(\bar \omega)  &=&  y_r ^ m \mbox{  if  } \bar \omega\in I ^ {m,n}_{k,r}, \quad n\in{\mathbb{N}}\cup\{\infty\},\end{eqnarray*}
and then, for $n\in{\mathbb{N}}\cup\{\infty\}$ \begin{eqnarray*} \bar \chi_n ^ 1(\bar
\omega)  &=& \lim_{m\to\infty} Z ^ 1_{n,m}(\bar \omega), \quad
\\\bar \chi_n ^ 2(\bar \omega)  &=& \lim_{m\to\infty} Z ^ 2_{n,m}(\bar \omega).
\end{eqnarray*} Due to the construction of the partition, the limits above
exist. To be precise, for any $n\in{\mathbb{N}}\cup\{\infty\}$ and $\bar
\omega\in\bar \Omega$, we have \begin{eqnarray}\label{c-s-z} |Z
^i_{n,m}(\bar \omega)-Z ^ i_{n,k}(\bar \omega)|\le r_m, \quad k\ge
m,\; i=1,2, \end{eqnarray} and therefore $(Z_{n,m}(\bar \omega))_{m\ge 1}
$ is  a Cauchy-sequence for all $\bar\omega\in\bar
\Omega=[0,1)\times [0,1)$. Hence, for $i,m$  $Z ^ i_n(\bar
\omega)$ is well defined.
Furthermore, $\chi_n$ is measurable, since $Z ^ i_{n,m}$ are
simple functions, hence measurable. Therefore, for $i,m$  $Z ^
i_n(\bar \omega)$ is a random variable.

Let $\bar \chi_n:= (\bar \chi _n^ 1,\bar \chi_n ^ 2
)$.

Finally we have to proof that the random variables $\bar \chi$ and
$\bar\chi_n$ have the following three properties:
\begin{numlist}
\item $\Law(\bar \chi_n) = \Law(\chi_n)$, $\forall n\in{\mathbb{N}}$,
\item $\chi_n\to \chi$ a.s.\ in $ E_1\times E_2$,
\item $\pi_1\circ \chi_n(\omega)= \pi_1 \circ \chi_\ast(\omega)$. % a.s.\
\end{numlist}
%\\[0.4cm]
\noindent {\sl Proof of (i):} The following identity  for the
rectangle $ A_k ^ l\times C_k ^ l$ holds
\begin{eqnarray*}
\leb (I_{k,r} ^ {l,n}) & = &\leb\left( \bar\chi_n\in A_k ^ l\times C_k ^ l\right)\\
&
=& {\mathcal{P}}MU_n ^ 1( A_k ^ l)\times {\mathcal{P}}MU_n( E_1 \times C_r ^ l\mid A_k
^ l\times E_2)
\\
& = & {\mathcal{P}}MU_n ( A_k ^ l\times E_2 )\times {\mathcal{P}}MU_n( E_1 \times
C_r ^ l\mid A_k ^ l\times E_2)\\ &=& {\mathcal{P}}MU_n( u ^ 1\in A_k ^ l
\mbox{ and } u ^ 2\in C_r ^ l)
\\& = & {\mathcal{P}}MU_1((u^1,u^2)\in A_k^l\times C_r^l).
\end{eqnarray*}
Now, one has to use that the set of rectangles of
$[0,1)\times [0,1)$ form a $\pi$ system in  $ {\mathcal{B}}([0,1)\times
[0,1))$. Moreover, $\leb$ and $ {\mathcal{P}}MU_1$ are identical on the
set of rectangles. Therefore, by Lemma 1.17 \cite[Chapter
1]{MR1876169}, $\leb$ and ${\mathcal{P}}MU$ are equal on $
{\mathcal{B}}([0,1)\times [0,1))$.
\\[0.4cm]
\noindent {\sl Proof of (ii):} We will first prove that there
exists a random variable $\chi=( \chi  ^ 1, \chi ^ 2)$ such that
$\bar \chi_n ^ 1\to
 \chi ^ 1$ and $\bar \chi_n ^ 2\to
 \chi ^ 2$ $\leb$-a.s.\ for $n\to\infty$.
For this it is enough to  show that the sequences $\{\bar \chi_n ^
1,n\in{\mathbb{N}}\}$ and $\{\bar \chi_n ^ 2,n\in{\mathbb{N}}\}$ are
$\leb$--a.s.\ Cauchy sequences. From the triangle inequality we
infer that  for all $ n,m,j\in{\mathbb{N}}$, $i=1,2$
$$
\left|\bar \chi_n  ^ i- \bar \chi ^ i_m \right|\le \left| \bar\chi ^
i_n-Z_{n,j} ^ i\right|+ \left|Z_{n,j} ^ i-Z_{m,j} ^ i\right| + \left|Z_{m,j}
^ i-  \bar \chi ^ i_m \right|.
$$

Let us
first observe that by \eqref{c-s-z},  for any $n\in{\mathbb{N}}$, the
sequences $\{Z_{n,j} ^ 1,j\in{\mathbb{N}}\}$ and $\{Z_{n,j} ^ 2,j\in{\mathbb{N}}\}$
converge uniformly on $\bar \Omega$   to $\bar \chi_n ^ 1$ and
$\bar \chi_n ^ 2$, respectively Therefore, it suffices to show
that for all $\eps>0$ there exists a number $n_0$ such that (see
\cite[Lemma 9.2.4]{dudley})
$$\leb\left( \left\{\omega\in [0,1)\times [0,1)\mid Z_{n,l}(\omega)\not=Z_{m,l}(\Omega), n,m\ge n_0 \right\}\right)\le \eps
$$

Since $\{\mu_n,n\in{\mathbb{N}}\}$ converges weakly,  for any $\delta>0$
there exists a number $n_0\in{\mathbb{N}}$ such that\footnote{$\rho$
denotes the Prokhorov metric on measurable space $(E,{\mathcal{E}})$, i.e.
$\rho(F,G):=\inf\{\eps>0: F(A ^ \eps)\leq G(A)+\eps,   A\in {\mathcal{E}}$\}}.
 $\rho (\mu_n,\mu_m)\leq \eps$ for all
$n,m\ge n_0$. %is a Cauchy sequence.
Hence,  for any $\delta>0$ we can find  a number $n_0\in{\mathbb{N}}$ such
that
\begin{eqnarray}\label{estt}
\\
\nonumber \left|a_{k,r} ^ {i,l,n}-a_{k,r} ^ {i,l,m}\right|,\left|b_{k,r}
^ {i,l,n}-b_{k,r} ^ {i,l,m}\right|\le \delta, \quad k,r,l\in{\mathbb{N}},
i=1,2, \, n,m\ge n_0, \end{eqnarray} where $$I ^ {l,n}_{k,r}=[a_{k,r} ^
{1,l,n},b_{k,r} ^ {1,l,n} )\times [a_{k,r} ^ {2,l,n},b_{k,r} ^
{2,l,n})$$ and $$I ^ {l,m}_{k,r}=[a_{k,r} ^ {1,l,m},b_{k,r} ^
{1,l,m} )\times [a_{k,r} ^ {2,l,m},b_{k,r} ^ {2,l,m}).$$
In fact, by the construction of $I ^ {l,n}_{k,r}$ and $I ^
{l}_{k,r}$ we have \begin{eqnarray*} a_{k,r} ^ {1,l,n} &=& \mu_n ^ 1\left(
\cup_{j=1} ^ {k-1} A_j^ l \right), \; n,m\in{\mathbb{N}},  r,k\in{\mathbb{N}},
\\
a_{k,r} ^ {1,l,m} &=& \mu ^ 1\left( \cup_{j=1} ^ {k-1} A_j ^ l
\right),\; n,m\in{\mathbb{N}},  r,k\in{\mathbb{N}}. \end{eqnarray*} Let us fix $\delta>0$ and let
us choose $n_0$ such that $\rho (\mu_n,\mu_m)\le \delta$, $ n,m\ge
n_0$. Then
$$
 \mu_n ^ 1\left( \cup_{j=1} ^ {k-1} A_j ^ l \right)\le  \mu_n ^ 1\left(\left( \cup_{j=1} ^ {k-1} A_j ^ l\right) ^ \delta \right)
\le \mu _m^ 1\left( \cup_{j=1} ^ {k-1} A_j ^ l \right)+\delta,\quad
r,k,l\in{\mathbb{N}}.
$$
On the other hand, by  symmetricity of  the Prokhorov metric,
$$
 \mu^ 1_m\left( \cup_{j=1} ^ {k-1} A_j ^ l \right)\le  \mu _m^ 1\left(\left( \cup_{j=1} ^ {k-1} A_j ^ l\right) ^ \delta \right)
\le \mu ^ 1_n\left( \cup_{j=1} ^ {k-1} A_j ^ l \right)+\delta,\quad
r,k\in{\mathbb{N}}.
$$
 Hence, we infer that
\begin{eqnarray}\label{estt-1}
 \left|a_{k,r} ^ {1,l,n}-a_{k,r} ^ {1,l,m}\right|\le \delta,\quad   r,k,l\in{\mathbb{N}}, \quad n,m\ge n_0.
\end{eqnarray} The second inequality in \eqref{estt} can be shown
similarly.

Since the sequence  $\{\mu_n,n\in{\mathbb{N}}\}$ is tight on $E_1\times
E_2$ we can find  a compact set $K_1\times K_2$ such that
$$
\sup_n\left(\mu_n \left((E_1\times E_2)\setminus (K_1\times
K_2)\right)\right)\le \frac \eps2.
$$
Let us fix $l\in{\mathbb{N}}$. Since the set $K_1 \times K_2$ is compact,
from the  covering  $(A_k ^ l\times C_r^l )_{k\in{\mathbb{N}},r\in{\mathbb{N}}}$ of
$E_1\times E_2$ there exists a finite covering $(A_k ^ l\times C_r
^ l)_{k=1, \cdots, K, j=1,\cdots,R}$ of $K_1 \times K_2$.
Next, observe that the estimate \eqref{estt}
is uniformly for all $n,m\ge n_0$.
Therefore we can use estimate \eqref{estt} with $\delta=\eps/2(KR)$
and  infer  that
$$
\sum_{k,r=1} ^ {K,R} \leb\left( I ^ {l,n}_{k,r}\bigtriangleup I ^
{l,m}_{k,r},n,m\ge n_0 \right)\le \frac \eps2.
$$
Moreover, since
\begin{eqnarray*} \lefteqn{ \leb\left( \left\{\omega\in [0,1)\times
[0,1): Z_{n,l}(\omega)\not=Z_{m,l}(\omega),\, n,m\ge n_0 \right\}\right)
\le}
%% %subseteq
&&
\\ &&
\le
\sum_{k,r=1}  \leb\left( I ^ {l,n}_{k,r}\bigtriangleup I ^
{l,m}_{k,r},\, n,m\ge n_0 \right)+ \leb\left((E_1\times
E_2)\setminus (K_1\times K_2)\right),
\end{eqnarray*} it follows that \begin{eqnarray*}
\leb\left( \left\{\omega\in [0,1)\times [0,1):
Z_{n,l}(\omega)\not=Z_{m,l}(\omega),\, n,m\ge n_0 \right\}\right)  &\le
& \frac\eps 2 +\frac \eps 2.
\end{eqnarray*}
Summarizing, we proved that for any $\eps>0$ there exists a number
$n_0\in{\mathbb{N}}$ such that $\left\{\omega\in [0,1)\times [0,1)\mid
Z_{n,l}(\omega)\not=Z_{m,l}(\Omega), n,m\ge n_0 \right\}\le \eps$.
Applying  then \cite[Lemma 9.2.4]{dudley} we infer (ii).
\\[0.4cm]

\noindent {\sl Proof of (iii):} Let us denote \begin{eqnarray}
\label{relation-aaa}
 J_{k} ^{l, n} := \left[\sum_{m=1} ^ {k-1}  {\mathcal{P}}MU ^ 1_n(A_k ^ l),\sum_{m=1} ^ {k}  {\mathcal{P}}MU ^ 1_n(A_k ^ l)   \right), \quad k=1,\cdots, N_l ^ 1.
\end{eqnarray}
Since the laws of $\pi_1\circ \chi_n$ and $\pi_1\circ
\chi_n$,  $n,m\in\deln{{\mathbb{N}}_0\cup\{\infty\}}\bar{\mathbb{N}}$ are  equal, $ J_{k} ^{l, n}=
J_{k} ^{l, m}$, $n,m\in\deln{{\mathbb{N}}_0\cup\{\infty\}}\bar{\mathbb{N}}$. Let us denote these
(equal) sets  $J_k ^ l$. Since  for each $J_k ^ l$ we can find
a set $A_k ^ l$ satisfying \eqref{relation-aaa} for
all $n\in\deln{{\mathbb{N}}_0\cup\{\infty\}}\bar{\mathbb{N}}$, we infer  that for any $m\in{\mathbb{N}}$,
$$ Z_{n,m} ^ 1 (\bar \omega)=  Z_{1,m} ^ 1 (\bar \omega), \quad  \bar \omega\in\Omega,\;  n\in\deln{{\mathbb{N}}_0\cup\{\infty\}}\bar{\mathbb{N}}.
$$
Let $n\in{\mathbb{N}}$ be fixed. Considering the limit for the sequence
$\{ Z_{n,m} ^ 1 (\bar \omega),m\in{\mathbb{N}}\}$ for  $m\to\infty$ and
keeping in mind  that $\{Z_{1,m} ^ 1 (\bar \omega),m\in{\mathbb{N}}\}$ is a
Cauchy sequence implies the assertion (iii).

\end{proof}
\section{A Tightness criteria in $\mathbb{D}([0,T];Y)$}
Let $Y$ be a separable and complete metric space and $T>0$. The
space $\mathbb{D}([0,T];Y)$ denotes the space of all right continuous
functions $x:[0,T]\to Y$ with left limits.
The space of continuous function is usually equipped with the
uniform topology. But, since $\mathbb{D}([0,T];Y)$ is complete but not
separable in the uniform topology, we equip $\mathbb{D}([0,T];Y)$ with the
Skorohod topology in which $\mathbb{D}([0,T];Y)$ is both separable and
complete. For more information about Skorokhod space and topology
we refer to Billingsley's book \cite{billingsley} or Ethier and
Kurtz \cite{MR838085}. In this appendix we only state the
following tightness criterium which is necessary for our work. For
this we denote by ${\mathcal{P}}\left( \mathbb{D}([0,T];Y)\right)$ the space of Borel
probability measures on $\mathbb{D}([0,T];Y)$.

\begin{theorem}\footnote{Compare with \cite[Chapter III, Theorem 13.5, p. 142]{billingsley}.}\label{billingley-comp}
A subset ${\mathcal{C}} \subset {\mathcal{P}}\left( \mathbb{D}([0,T];Y)\right)$ is tight, iff
\begin{enumerate}
\item  for any $\eps>0$ there exists a compact set $K_\eps\subset Y$
such that
$$ F\left( \left\{ x \in \mathbb{D}([0,T];Y): x(t)\in K_\eps\, \,\forall\, t\in [0,T]\right\} \right) \ge 1-\eps
,\,\quad  \forall F\in{\mathcal{C}}  ;
$$
\item there exist two real numbers $\gamma>0$, $\alpha >\frac 12$
and a nondecreasing continuous function $g:[0,T]\to{\mathbb{R}}_+$ such
that for all $t_1\le t\le t_2$, $ n\ge 0$ and $\lambda>0$
\begin{eqnarray*}\lefteqn{ F\left( \left\{  x \in \mathbb{D}([0,T];Y):\, |x(t)-x(t_1)|\ge
\lambda, \, |x(t)-x(t_2)|\ge \lambda\right\} \right)
} &&\\
&& \phantom{mmmmmmmm}\le \frac{1}{ \lambda ^ {2\gamma}}\, \left[
g(t_2)-g(t_1)\right],\quad \forall F\in{\mathcal{C}}.
%\hspace{3cm}
\phantom{mmmmmmmm} \end{eqnarray*}
\end{enumerate}
\end{theorem}

\begin{corollary}\label{comp-2}
Let $\{ x_n:n\in {\mathbb{N}}\}$ be a sequence of  c\`{a}dl\`{a}g processes,
each of the process defined on a probability space $
(\Omega_n,{\mathcal{F}}_n,\mathbb{P}_n)$.
Then the sequence of laws of  $\{ x_n:n\in {\mathbb{N}}\}$ %the sequence
is tight on  $\mathbb{D}([0,T];Y)$ if
\begin{enumerate}[(a)]
\item  for any $\eps>0$ there exists a compact set $K_\eps\subset Y$
such that
$$ \mathbb{P}_n\left( x_n(t)\in K_\eps, \, t\in [0,T]\right) \ge 1-\eps
,\, \forall n\in{\mathbb{N}} ;
$$
\item there exist two constants $c>0$ and $\gamma> 0$ and a real
number $r>0$ such that for all $\theta>0$, $t\in[0,T-\theta]$, and
$ n\ge 0$
$$
 {\mathbb{E}}_n\sup_{t\le s\le t+\theta} |x_n(t)-x_n(s )| ^ r\le c\, \theta ^ \gamma.
$$
\end{enumerate}
\end{corollary}
\begin{proof}
The inequality \ref{comp-2}-(a) and the
Chebyscheff  inequality imply Inequality \ref{billingley-comp}-(a).
Now fix $t_1\le t\le t_2$. Then \begin{eqnarray*}\lefteqn{
 \mathbb{P} _ n\left( |x_n(t)-x_n(t_1)|\ge \lambda, \, |x_n(t)-x_n(t_2)|\ge \lambda\right) }
&&\\&& \le \mathbb{P}_n\left( \sup_{t_1\le s\le t_2} |x_n(s)-x_n(t_1)|\ge
\lambda\right). \end{eqnarray*} Estimating the RHS by the Chebyshev
inequality and using inequality \ref{comp-2}-(b) leads to inequality
\ref{billingley-comp}-(b). Thus the assertion follows.
\end{proof}

\section{An inequality}\label{sec-inequality}

Let $Y$ is a Banach space with norm $\lve \cdot \rve$, $T>0$  and
$f:(0,T]\to Y$ is a Bochner integrable function.
$$\int_0^T |f(s)|^p\,ds <\infty.$$
For fixed $n$ let $I_k=(\frac{k} {2^n},\frac{k+1} {2^n}]$ and
$\hat{f}_n:(0,T]\to Y$ be the function defined by $\hat{f}_n(s)=0$
for $s\in I_0$ and $\hat{f}_n(s)={2^n} \int_{I_k}\, f(t) \,dt $
for $s\in I_{k+1}$, $k=0,1,2,\cdots$.
We have the following facts.
\begin{proposition}\label{convergen}
%\phantom{mm} \\ \phantom{mm}  %The following holds
\begin{itemize}
\item[(i)] If $f$ belongs to $L^p(0,T;Y)$, then so is $\hat{f}_n$
for each $n$. \item[(ii)] Let $\alpha \in (0,\frac1p)$ and $f\in
W^{\alpha,p}(0,T;Y)$. Then there exists $C>0$ such that
\begin{eqnarray}\label{aninequality} \left\|
f(s)-\hat{f}_n(s)\right\|_{L^p(0,T;Y)} &\le & C 2^{-n\alpha}\,
\left\|f\right\|_{W^{p,\alpha}(0,T;Y)},
\end{eqnarray}
for all $n$.
\end{itemize}
\end{proposition}
\noindent
\begin{proof}
Without loss of generality we take $T=1$ and set
$t^n_j=\frac{j}{2^n}$. We have
\begin{align*}
\int_0^1 \lve \hat{f}_n(s)\rve^p ds=&\sum_{j=1}^{2^n-1}
\int_{t^n_j}^{t_{j+1}^n}\lve \hat{f}_n(s)\rve^p ds \\
=&\sum_{j=1}^{2^n-1}\frac{1}{2^n} \lve \hat{f}_n(t^n_{j})\rve^p.
\end{align*}
From this last equation and the definition of $\hat{f}$ we derive
\begin{align*}
\int_0^1 \lve \hat{f}_n(s)\rve^p ds\le&
\sum_{j=1}^{2^n-1}\frac{1}{2^n} \biggl\lve 2^n
\int_{t^n_{j-1}}^{t^n_{j}} f(s) ds\biggr\rve^p\\
\le & \sum_{j=1}^{2^n-1} \int_{t^n_{j-1}}^{t^n_{j}}\lve f(s)\rve^p
ds.
\end{align*}
Therefore $$ \int_0^1 \lve \hat{f}_n(s)\rve^p ds\le
\int_0^{1-\frac{1}{2^n}} \lve f(s)\rve^p ds,$$ which ends the
proof of (i).

Now we prove item (ii). Let $s\in [0,1]$, $\alpha\in (0,\frac
1p)$ and $f\in W^{\alpha,p}(0,T;Y)$. Since the intervals $I_k$,
$k=0,1,...,2^n-1,$ form a partition of $[0,1]$ then either $s\in
I_0$ or $s\in (2^{-n},1]$. In one hand if $s\in (2^{-n},1]$ then
there exists $k\ge1$ such that $s\in I_k$. In this case we have by
H\"older's inequality that
\begin{align*}
\lve f(s)-\hat{f}_n(s)\rve\le  2^n \biggl\lve \int_{I_{k-1}}
\frac{\left(f(s)-f(r)\right)}{\lve s-r\rve^{\frac1p+\alpha}}
\times \lve s-r\rve^{\frac 1p+\alpha} dr\biggr\rve\\
\le  2^n\biggl[\int_{I_{k-1}} \lve
s-r\rve^{\frac{p}{p-1}(1+\alpha)-1}
dr\biggr]^\frac{p-1}{p}\biggl[\int_{I_{k-1}}\frac{\lve
f(s)-f(r)\rve^p}{\lve s-r\rve^{1+p\alpha}}dr\biggr]^\frac1p\\
\le  2^{-np\alpha}\biggl[\int_{I_{k-1}}\frac{\lve
f(s)-f(r)\rve^p}{\lve s-r\rve^{1+p\alpha}}dr\biggr]^\frac1p.
\end{align*}
Therefore
\begin{eqnarray}
\int_{2^{-n}}^1 \lve f(s)-\hat{f}_n(s)\rve^p ds &\le &
2^{-np\alpha}\int_{2^{-n}}^1
\biggl[\sum_{k=1}^{2^n-1}\int_{I_{k-1}}\frac{\lve
f(s)-f(r)\rve^p}{\lve s-r\rve^{1+p\alpha}}\,dr\biggr]ds\nonumber\\
&&\hspace{-3truecm}\lefteqn{\le  2^{-np\alpha}\int_{0}^1 \int_0^1\frac{\lve
f(s)-f(r)\rve^p}{\lve s-r\rve^{1+p\alpha}}dr\,ds
\leq  2^{-np\alpha} \Lve
f\Rve^p_{W^{\alpha,p}(0,T;Y)}.}
\label{ineq-fhat-1}
\end{eqnarray}
On the other hand, by making use of H\"older's inequality we
obtain
\begin{eqnarray*}
\int_0^{2^{-n}}\lve f(s)-\hat{f}_n(s)\rve^p \,ds&=&
\int_0^{2^{-n}}\lve f(s)\rve^p\,ds\\
&\leq & 2^{-np\alpha}\biggl[\int_0^{2^{-n}} \lve
f(s)\rve^qds\biggr]^\frac pq,
\end{eqnarray*}
where $q=\frac{p}{1-p\alpha}$. Since $L^r(I_0)\subset
W^{\alpha,p}(I_0)$ for any $r\in [p,\frac{p}{1-p\alpha}]$ we
infer from the last inequality that there exists $C>0$ such that
 \begin{equation}
\int_0^{2^{-n}}\lve f(s)-\hat{f}_n(s)\rve^p ds\le C
2^{-np\alpha}\Lve
f\Rve^p_{W^{\alpha,p}(0,T;Y)}.\label{ineq-fhat-2}
 \end{equation}
 Now  inequality \eqref{aninequality} follows from inequalities
 \eqref{ineq-fhat-1} and \eqref{ineq-fhat-2}. This completes
 the proof of our proposition.
\end{proof}

\end{document}